\documentclass[twoside]{article}
\usepackage[a4paper]{geometry}
\usepackage{a4wide}
\usepackage[latin1]{inputenc} 
\usepackage[T1]{fontenc} 
\usepackage{RR}
\usepackage{hyperref}

\usepackage{float,leqno,epsf}
\usepackage{times}
\usepackage{graphicx,epsfig}
\usepackage{amsmath,amssymb}

\RRNo{8442}
\RRdate{December 2013}
\RRauthor{
Adimurthi\thanks[sfn]{TIFR-CAM, PB 6503, Chikkabommasandra, Sharadanagar, GKVK PO Bangalore-560065, India}\thanks{\texttt{aditi@math.tifrbng.res.in}}%
  \and
G. D. Veerappa Gowda\thanksref{sfn}\thanks{\texttt{gowda@math.tifrbng.res.in}}%
\and  J\'er\^ome Jaffr\'e\thanks{INRIA, BP 105, 78153 Le Chesnay Cedex, France}\thanks{\texttt{Jerome.Jaffre@inria.fr}}
}
\RRtitle{Le flux DFLU pour les syst\`emes de lois de conservation}
\RRetitle{The DFLU flux for systems of conservation laws\thanks{This work was partially supported by the Indo-French collaboration project IFCPAR/CEFIPRA 3401-2.}}
\titlehead{The DFLU flux for systems of conservation laws}
\RRnote{This paper is published in the Journal of Computational and Applied Mathematics 247 (2013) 102-123.}

\RRresume{Le flux num\'erique DFLU a \'et\'e introduit pour r\'esoudre des lois de conservation scalaires hyperboliques dont la fonction de flux est discontinues en espace. Nous montrons comment ce flux peut \^etre utilis\'e pour r\'esoudre une certaine classe de syst\`emes de lois de conservation tels que les syst\`emes mod\'elisant l'injection de polym\`eres en ing\'eni\'erie de r\'eservoirs p\'etroliers. En outre, ces r\'esultats s'\'etendent au cas de fonctions de flux discontinus par rapport \`a la variable d'espace. Une telle situation appara\^\i t par exemple quand on consid\`ere des r\'eservoirs p\'etroliers qui sont h\'et\'erog\`enes. Des exp\'eriences num\'eriques sont pr\'esent\'ees pour illustrer l'efficacit\'e de ce nouveau sch\'ema  compar\'e \`a d'autres sh\'emas standard tels que les sch\'emas Mobilit\'es Amont, Lax-Friedrichs et Force.}
\RRabstract{
The DFLU numerical flux was introduced in order to solve hyperbolic scalar conservation laws with a flux function discontinuous in space. We show how this flux  can be used to solve certain class of  systems of conservation laws such as systems  modeling polymer flooding in oil reservoir engineering. Furthermore, these results are extended to  the case where the flux function is discontinuous in the space variable. Such a situation arises for example while dealing  with oil reservoirs which are heterogeneous. Numerical experiments are presented to illustrate the efficiency of this  new scheme compared to other standard schemes like
upstream mobility, Lax-Friedrichs and Force schemes.}

\RRmotcle{Volumes finis, diff\'erences finies, solveurs de Riemann, syst\`emes de lois de conservation, \'ecoulements en milieu poreux, injection de polym\`eres.}
\RRkeyword{ Finite volumes, finite differences, Riemann solvers, system of conservation laws, flow in porous media, polymer flooding.}

\RRprojet{Pomdapi}

\usepackage{float,leqno,epsf}
\usepackage{times}
\usepackage{graphicx,epsfig}
\usepackage{amsmath,amssymb}
\usepackage{amssymb}
\newtheorem{lemma}{\bf  Lemma }[section]
\newcommand{\be} {\begin{equation}}
\newcommand{\ee} {\end{equation}}
\newcommand{\bea} {\begin{eqnarray}}
\newcommand{\eea} {\end{eqnarray}}
\newcommand{\Bea} {\begin{eqnarray*}}
\newcommand{\Eea} {\end{eqnarray*}}
\newcommand{\p} {\partial}

\newcommand{\la} {\lambda}

\newcommand{\no} {\nonumber}
\newcommand{\noi} {\noindent}

\newcommand{\var} {\varphi}
\newcommand{\R}{{\mathbb R}} 
\newcommand{\Z}{{\mathbb Z}} 
\newcommand{\cqfd}{\hfill \rule{1.5mm}{3mm}}
\RCParis 
\begin{document}

\makeRR

\section{Introduction}

 The main difficulty in the numerical solution of systems of conservation laws is the complexity of constructing the Riemann solvers. One way to overcome this difficulty is to consider centered schemes as in \cite{LaxWen60,NesTad90,Toro99,Toro06,AdiGowJaf09}. However, in general these schemes are more diffusive than Godunov type methods based on exact or approximate Riemann solvers when this alternative is available. Therefore in this paper we will consider Godunov type methods. Most often the numerical solution requires the calculation of eigenvalues or eigenvectors of the Jacobian matrix of the system. This is even more complicated when the system  is non-strictly hyperbolic, i.e. eigenvectors are not linearly independent. In this paper we present an  approach like in \cite{KarMisRis09a} and \cite{KarMisRis09c} which do not require, detailed information about the eigenstructure of the full system.

Let us consider a system of conservation laws in conservative form
\[ \mathbf{U}_t + (\mathbf{F}(\mathbf{U}))_x=0, \quad \mathbf{U}=(u^1,\cdots,u^J), \quad \mathbf{F}=(f^1,\cdots,f^J). \]
A conservative finite volume method reads
\[ \dfrac{\mathbf{U^{n+1}_i}-\mathbf{U^n_i}}{\Delta t} + \dfrac{\mathbf{F}^{n}_{i+1/2} - \mathbf{F}^{n}_{i-1/2}}{h} =0 \]
where $\mathbf{F}^{n}_{i+1/2}$ is a numerical flux calculated using an exact or approximate Riemann solver. 
In a first order scheme this numerical flux is calculated using the left and right values 
$\mathbf{U}_i^n$ and $\mathbf{U}_{i+1}^n$. If we solve the equation field by field the $j$-th equation reads
\[ \dfrac{u^{j,n+1}_i-u^{j,n}_i}{\Delta t} + \dfrac{F^{j,n}_{i+1/2} - F^{j,n}_{i-1/2}}{h} =0\]
where the $j$-th numerical flux is a function of $\mathbf{U}^n_i$ and $\mathbf{U}^n_{i+1}$:
\[ F^{j,n}_{i+1/2} = F^j(u^{1,n}_i,\cdots, u^{j,n}_i,\cdots,u^{J,n}_i,u^{1,n}_{i+1},\cdots, u^{j,n}_{i+1},\cdots,u^{J,n}_{i+1}), \quad j=1,\cdots,J. \]
This flux function  can be calculated by solving the scalar Riemann problem for $t>t_n$:
\begin{align}
\label{dconsl}
&u^j_t + (\tilde{f}^{j,n}(u^j,x))_x = 0,\\ 
&u^j(x,t_n)=u^{j,n}_i \mbox{ if } x <x_{i+1/2}, \; u^j(x,t_n)=u^{j,n}_{i+1} \mbox{ if } x >x_{i+1/2}, \nonumber
\end{align}
where the flux function $\tilde{f}^j$, discontinuous at the point $x=x_{i+1/2}$, is defined by
\be \begin{array}{l} 
\tilde{f}^{j,n}(u^j,x) \equiv \tilde{f}^{j,n}_L(u^j) = f^j(u^{1,n}_i,\cdots, u^{j-1,n}_i,u^{j},u^{j+1,n}_i,\cdots,u^{J,n}_i) \mbox{ if } x <x_{i+1/2},\\
\tilde{f}^{j,n}(u^j,x) \equiv  \tilde{f}^{j,n}_R(u^j) = f^j(u^{1,n}_{i+1},\cdots, u^{j-1,n}_{i+1},u^{j},u^{j+1,n}_{i+1},\cdots,u^{J,n}_{i+1}) \mbox{ if } x >x_{i+1/2}
\label{fluj}
  \end{array} \ee
(L and R refer to left and right of the point $x_{i+1/2}$).

Scalar conservation laws like equation (\ref{dconsl}) with a flux function discontinuous in space have been the object of many studies \cite{chacohjaf87,mochen87,lantvewin92,gimris92,jaf96,kaa99,Tow00,Tow01,BurKarRisTow03,KarRisTow03,SegVov03,AdiJafGow04, SidJaf09}. In particular, in \cite{AdiJafGow04} a Godunov type finite volume scheme was proposed and convergence to a proper entropy condition was proved, provided that the left and right flux functions have exactly one local maximum and the same end points (the case where the flux functions has exactly one local minimum can be treated by symmetry).   At the discontinuity the interface flux, that we call the DFLU flux, is given by the formula 
\be \begin{array}{l} F^n_{i+1/2}(u_L,u_R) =  \min\{ f_L(\min\{u_L,\theta_L\}), f_R(\max\{u_R,\theta_R\})\}, 
\label{dflu} \end{array}\ee
 if $f$ denotes the scalar flux function and $\theta_L=$argmax$(f_L)$, $\theta_R=$argmax$(f_R)$. When $f_L \equiv f_R$ this formula is equivalent to the Godunov flux so formula (\ref{dflu}) can be seen as an extension of the Godunov flux to the case of a flux function discontinuous in space.
In the case of systems formula (\ref{dflu}) can be applied to the fluxes $\tilde{f}^{j,n}_L$ and $\tilde{f}^{j,n}_R$ provided both agrees at the end points of the domain for all $j$, like in the case of scalar laws with a flux function discontinuous in space. In the case of an uncoupled triangular system, a similar scheme  is used in \cite{KarMisRis08,KarMisRis09a,KarMisRis09b} and its convergence analysis is studied. Also in  \cite{KarMisRis09c}, the  idea of discontinuous flux is used to study a coupled system arising in three-phase flows in porous media and shown its successfulness.

To illustrate the method we consider the system of conservation laws arising for polymer flooding in reservoir simulation which is described in section \ref{polymer}.
This system, or similar systems of equations, is nonstrictly hyperbolic and is studied in several papers \cite{Tem82, JohWin88, JohTveWin89,IssTem90}. For example in \cite{JohWin88} the authors solve Riemann problems associated to this system when gravity is neglected and therefore the fractional flow function is an increasing function of the unknown.  In this case, the eigenvalues of the corresponding Jacobian matrix  are positive and hence it is less difficult to construct Godunov type schemes which turn out to be upwind schemes. When  the above model with gravity effects is considered, then the flux function is not necessarily monotone and hence the eigenvalues can change sign. This makes the construction of Godunov type schemes more difficult as it involves exact solutions of Riemann problems with a non monotonous fractional flow function.
Therefore in section \ref{Riemann} we solve the  Riemann problems in the general case when gravity terms are taken into account so the flux function is not anymore monotone. This will allow to compare our method with that using an exact Riemann solver.
In section \ref{finite difference} we consider Godunov type finite volume schemes. We present the DFLU scheme for the system of polymer flooding and compare it to the Godunov scheme whose flux is given by the exact solution of the Riemann problem. We also present several other possible numerical fluxes, centered like Lax-Friedrichs or FORCE, or upstream like the upstream mobility flux commonly used in reservoir engineering \cite{AziSet79, BreJaf91,SidJaf09}.
In section \ref{numericalresults} we compare numerically the DFLU method with these fluxes.
Finally in section \ref{dispolymer} we considered the case where the flux function is discontinuous in the space variable and its corresponding Riemann problem  is discussed in appendix.
\section{A system of conservation laws modeling polymer flooding}
\label{polymer}
A polymer flooding model for enhanced oil recovery in petroleum engineering was introduced in \cite{Pop80} as the following $2\times 2$ system of conservation laws
\be
\begin{array}{rrll}
 s_t + f(s, c)_x & =& 0 \\
 (sc+a(c))_t+(cf(s,c))_x &=&0
 \label{polyeq1}
\end{array}
\ee
where $t > 0$ and $x \in \R$, $(s,c) \in I \times I$ with $I=[0,1]$. $s=s(x,t)$ denotes the saturation of the wetting phase, so $1-s$ is  the saturation of the oil phase. $c=c(x,t)$ denotes the concentration of the polymer in the wetting phase which we have normalized. Here the porosity was set to 1 to simplify notations. The flux function $f$ is the Darcy velocity of the wetting phase  $\var_1$ and is determined by the relative permeabilities and the mobilities of the wetting and oil phases, and by the influence of gravity:
\be
\label{fg-twophase}
f(s,c)=  \var_1 = \dfrac{\lambda_1(s,c)}{\lambda_1(s,c) + \lambda_2(s,c)} [ \var + (g_1-g_2)\lambda_2(s,c) ].
\ee
The quantities $\lambda_\ell, \ell=1,2$  are the mobilities of the two phases, with $\ell =1$ referring to the wetting phase and  $\ell =2$ referring to the oil phase:
\[ \lambda_\ell(s,c) = \dfrac{K kr_\ell(s)}{\mu_l(c)}, \ell=1,2, \]
where $K$ is the absolute permeability, and $kr_\ell$ and $\mu_\ell$ are respectively the relative permeability and the viscosity of the phase $\ell$.
$kr_1$ is an increasing function of $s$ such that $kr_1(0)=0$ while $kr_2$ is a decreasing function of $s$ such that $kr_2(1)=0$. Therefore $\lambda_\ell, \ell=1,2$ satisfy
\begin{equation}
\label{eq3}
\begin{array}{l}
\lambda_1=\lambda_1(s,c) \mbox{is an  increasing functions of} \, s ,\;
\lambda_1(0,c) = 0\;  \forall c \in [0,1],\\
\lambda_2= \lambda_2(s,c) \, \mbox{is a decreasing functions of} \, s ,\;
\lambda_2(1,c)= 0\;  \forall c \in [0,1].
\end{array}
\end{equation}
The idea of polymer flooding is to dissolve a polymer in the injected water in order to increase the viscosity of the injected wetting phase. Thus the injected wetting phase will not be able to bypass oil so one obtains a better displacement of the oil by the injected phase. Therefore $\mu_1(c)$ is increasing with $c$ while $\mu_2$ will be taken as a constant assuming there is no chemical reaction between the polymer and the oil.  Therefore $f$ will decrease with respect to $c$.
The function $a=a(c)$ models the adsorption of the polymer by the rock and is increasing with $c$. 

$\var$ is the total Darcy velocity, that is the sum of the Darcy velocities of the
two phases $\var_1$ and $\var_2$:
\[ \var = \var_1+\var_2, \quad \var_1 = \dfrac{\lambda_1}{\lambda_1 + \lambda_2} [ \var + (g_1-g_2)\lambda_2], \quad
\var_2 = \dfrac{\lambda_2}{\lambda_1 + \lambda_2} [ \var + (g_2-g_1)\lambda_1].\]
$\var $ is a constant in space since we assume that the flow is incompressible.
The gravity constants $g_1, g_2$ of the phases are proportional to their density.

To equation (\ref{polyeq1})  we add the initial condition
\begin{equation}
(s(x,0), c(x, 0)) = (s_0 (x),c_0(x)). \label{inicond}
\end{equation}

\begin{figure}[htbp]
 \includegraphics[width=6.2cm]{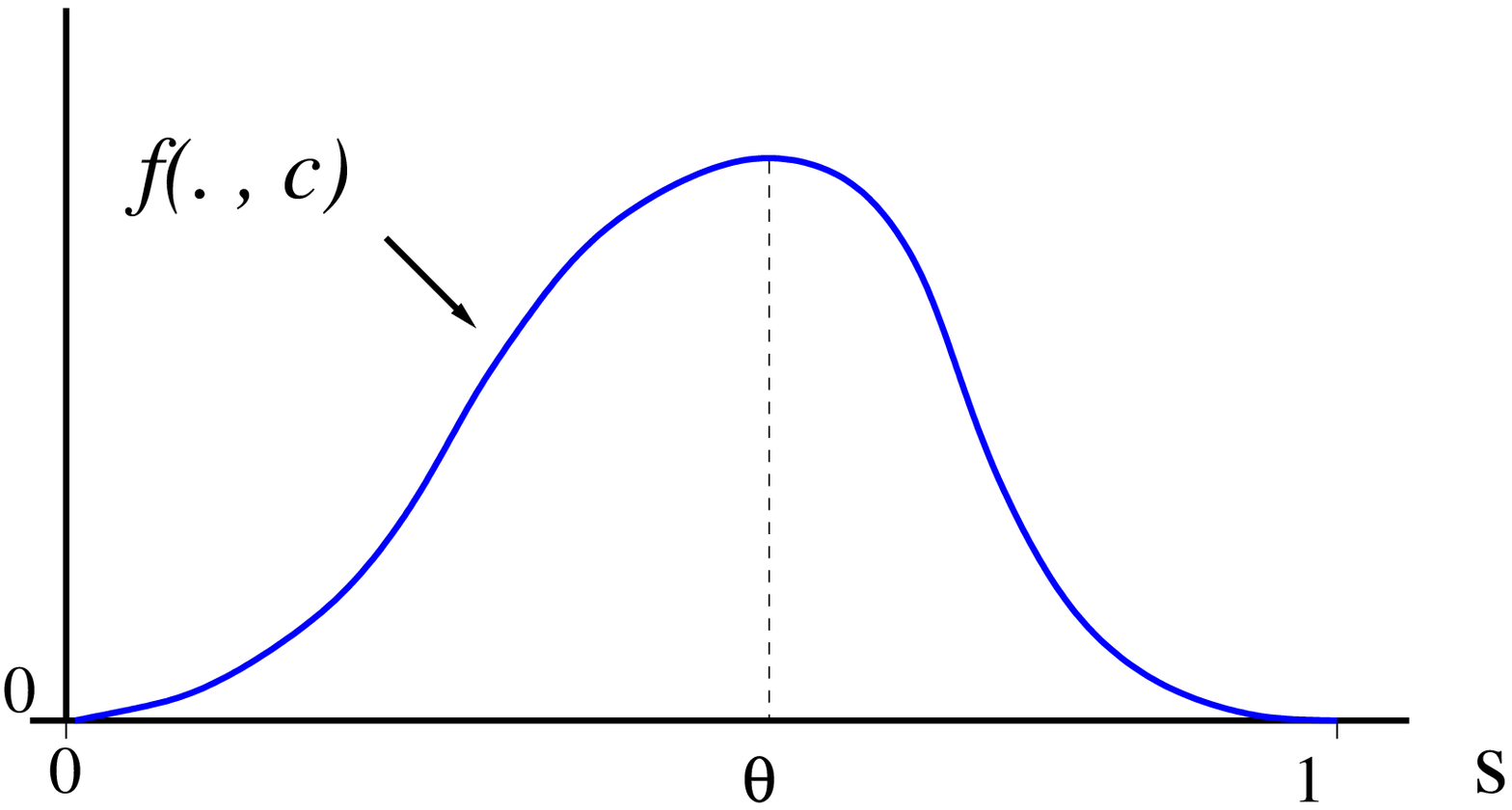} \hspace*{0.3cm}
 \includegraphics[width=4.2cm]{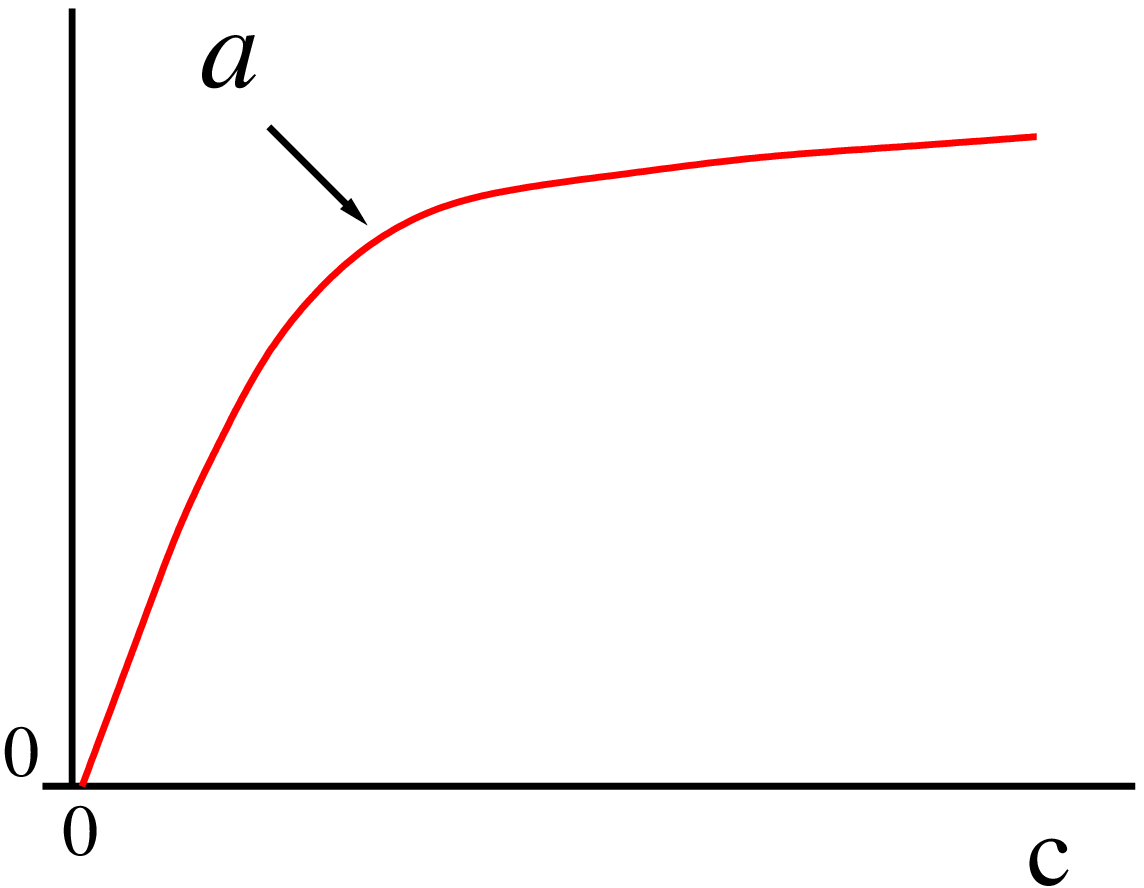}
\caption{Shapes of flux function $s \rightarrow f(s,c)$ (left) and adsorption function $c \rightarrow a(c)$ (right).}
\label{shape-f-a}
\end{figure}

Since the case when $f$ is monotone was already studied in \cite{JohWin88,JohTveWin89}, we concentrate on the nonmonotone case which is more complicated and corresponds to taking into account gravity. Here  we assume that $\var = 0$  for the
nonlinearities of the system (\ref{polyeq1}). We will assume also that phase 1 is heavier than phase 2 ($g_1>g_2$) so  we can assume the following properties:
\begin{enumerate}
\item[(i)] $f(s,c) \ge 0, f(0,c)= f(1,c)=0$ for all $c \in I$.
\item[(ii)] The function $s \rightarrow f(s,c)$ has exactly one global
 maximum in $I$ and no other local minima in the interior of $I$ with $\theta=$argmax$(f)$.
\item[(iii)] $f_c(s,c) < 0  \,\,\forall \,\,s \in (0,1)$ and for all $ c \in I$
\item[(iv)] The adsorption term $a=a(c)$ satisfies\\
   $a(0)=0,\,\,\,\, h(c)=\dfrac{da}{dc}(c)>0, \quad \dfrac{d^2a}{dc^2}(c) < 0$
 for all $c \in I$.  
 \end{enumerate}
 Typical shapes of functions $f$ and $a$ are shown in Fig. \ref{shape-f-a}.

We expand the derivatives in equations  (\ref{polyeq1}) and we plug the resulting first equation into the second one. Then we obtain the system in nonconservative form
 \[
\begin{array}{rrll}
 s_t + f_s(s, c) s_x + f_c(s,c) c_x & =& 0, \\
 (s +a'(c)) c_t+f(s,c) c_x &=&0.
\end{array}
\]
 Let $U$ denote the state vector $U=(s,c)$ and introduce the upper triangular matrix
  \[  A(U) =  \left( \begin{array}{cc}
  f_s & f_c \\[0.2cm] 0 & \dfrac{f}{s + a'(c)}\end{array} \right) \]
and the system (\ref{polyeq1}) can be read in matrix form as
$$ U_t+ A(U)\,U_x\,=\,0.$$   
  
 The eigenvalues of $A$ are $\lambda^s=f_s$ and $\lambda^c=\dfrac{f}{s+a'}$, with corresponding eigenvectors $e^s=(1,0), e^c=(f_c,\lambda^c-\lambda^s) $ if $0<s<1$   and $e^c=(0,1)$ if $s=0,1$. The
 eigenvalue  $\lambda^s$ may change sign whereas the eigenvalue $\lambda^c$ 
 is always positive. One can observe that for each $c \in I$ there exists a unique $s^*=s^*(c) \in (0,1)$ such that
$$ \lambda^c(s^*,c)=\lambda^s(s^*,c)$$
(see Fig.\ref{rs0}). For this couple $(s^*,c)$, $\lambda^c=\lambda^s$, hence eigenvectors are not linearly independent and the problem is nonstrictly hyperbolic.

 Any weak solution of (\ref{polyeq1}) has to satisfy the Rankine-Hugoniot
 jump  conditions given by
\be
\begin{array}{rrll}
 f(s_R,c_R)-f(s_L,c_L) & =& \sigma(s_R -s_L),\\
 c_R f(s_R,c_R)-c_L f(s_L,c_L)&=&\sigma (s_R c_R +a(c_R) -s_L c_L - a(c_L)),
\label{RHcondition1}
\end{array}
\ee
where $(s_L,c_L), (s_R,c_R)$ denote the left and right values of the couple $(s,c)$ at a certain point of discontinuity.

When $c_R=c_L$, the second equation reduces to the first equation and the speed of the discontinuity $\sigma$ is given by the first equation only. Now we are interested in the case $c_R \neq c_L$.
By combining the two equations (\ref{RHcondition1}) we may write
$$ (c_R - c_L)f(s_L,c_L)=\sigma(c_R - c_L)s_L + \sigma( a(c_R) - a(c_L)) $$
where
\[ \sigma=\frac{f(s_L,c_L)}{s_L+\bar{a}_L(c_R)}, \quad \bar{a}_L(c)=\left \{\begin{array}{lll}
\dfrac{a(c)-a(c_L)}{c-c_L} &\mbox{if} &c \neq c_L,\\
a'(c) &\mbox{if}& c=c_L. \end{array} \right. 
\]
Plugging this into first equation of (\ref{RHcondition1}), we obtain
$$\sigma(s_R + \bar{a}_L(c_R))=\sigma(s_L + \bar{a}_L(c_R))+f(s_R,c_R)-f(s_L,c_L)=f(s_R,c_R).$$
Hence when $c_L \neq c_R$ the Rankine-Hugoniot condition ({\ref{RHcondition1}})
reduces to
\begin{equation}
\dfrac{f(s_R,c_R)}{s_R+\bar{a}_L(c_R)}=\dfrac{f(s_L,c_L)}{s_L +\bar{a}_L(c_R)}=\sigma.
\label{rhcondition2}
\end{equation}

       In the absence of the adsorption term, i.e. $a=a(c)=0$, equation 
(\ref{polyeq1}) is studied in \cite{KliRis01} by using the equivalence of the Euler and Lagrangian 
 formulations and converting it into a scalar conservation law with a discontinuous flux function. In the presence of the adsorption term, this transformation fails to convert it into a scalar conservation law with a discontinuous flux function.

\section{Riemann problem}
\label{Riemann}
In this section we  solve the Riemann problems associated with our system, that we solve system (\ref{polyeq1}) with the initial condition
\be s(x,0) = \left\{ \begin{array}{lll}
s_L &\mbox{if}& x<0, \\ s_R &\mbox{if}& x>0 \end{array} \right.  , \quad
c(x,0) = \left\{ \begin{array}{lll}
c_L &\mbox{if}& x<0, \\ c_R &\mbox{if}& x>0 \end{array} \right.  .
\label{riemannpb}
\ee
Solution to (\ref{riemannpb}) is constructed by using elementary waves associated with the system. There are two families of waves, refered to as the $s$ and $c$ families. $s$ waves consist of rarefaction and shocks (or contact discontinuity) across which $s$ changes continuously and discontinuously respectively, but across which $c$ remains constant. $c$ waves consist solely of contact discontinuities, across which both $s$ and $c$ changes such that $\dfrac{f(s,c)}{s+a'(c)}$ remains constant in the sense of (\ref{rhcondition2}).

\begin{figure}[H]
\begin{center}
 \begin{picture}(200,100)(0,-10)
 \put(0,0){\includegraphics[height=3cm]{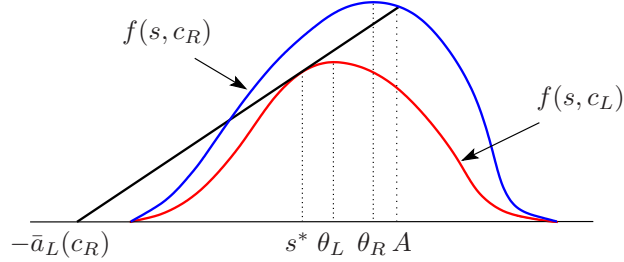}}
 \put(35,70){$f(s,c_R)$}
 \put(190,45){$f(s,c_L)$}
 \put(-7,-10){$-\bar{a}_L(c_R)$}
 \put(96,-10){$s^*$}
  \put(123,-10){$\theta_R$}
   \put(108,-10){$\theta_L$}
   \put(136,-10){$A$}
 \end{picture}
 \end{center}
 \caption{Two flux functions $f(s,c_L)$ and  $f(s,c_R)$ with $c_L > c_R$.}
 \label{rs0}
 \end{figure}

We will restrict to the case $c_L > c_R$. The case  $c_L <  c_R$ can be treated similarly. When  $c_L > c_R$ the flux functions for the first equation (\ref{polyeq1}) $s \rightarrow f(s,c_L)$ and $s \rightarrow f(s,c_R)$ are as represented in 
Fig. \ref{rs0}, that is $f(s,c_L) \leq f(s,c_R) \; \forall s \in (0,1)$. Let 
$\theta_L$ and $\theta_R$ be the points at which $f(.,c_L)$ and $f(.,c_R)$ reach their
maxima respectively.

 Let $s^* \in (0,1)$ be a point at which  
$ f_s(s^*,c_L)=\dfrac{f(s^*,c_L)}{s^*+\bar{a}_L(c_R)}$. Now draw a line through the points
$(-\bar{a}_L(c_R),0)$ and $(s^*,f(s^*,c_L))$ which intersects the curve $f(s,c_R)$ at a point $A \geq s^*$ (see Fig. \ref{rs0}).

 Our study of Riemann problems separates into two cases $s_L < s^* $ and $s_L \geq s^*$ which themselves separate into several subcases.

 \begin{itemize}
 \item {\bf Case 1:} $ s_L < s^*$.

    Draw a line passing through the points $(s_L,f(s_L,c_L))$ and $(-\bar{a}_L(c_R),0)$. This line intersects the curve $f(s,c_R)$ at points $\overline{s}$ and $B$ (see Fig. \ref{figR1a} ). Now we divide this into two subcases.  
 
 \item {\bf Case 1a:} $s_R < B$\\
  (a)  Connect $(s_L,c_L)$ to $(\overline{s},c_R)$ by $c$-wave with a speed 
 $$\sigma_c=\frac{f(s_L,c_L)}{s_L+\bar{a}_L(c_R)}=\frac{f(\overline{s},c_R)}{\overline{s}+\bar{a}_L(c_R)}. $$
  (b) Next connect  $(\overline{s},c_R)$ to $(s_R,c_R)$ by a $s$-wave, along the curve
     $f(s,c_R)$ (see Fig. \ref{figR1a}).

   For example if $ s_R \geq \overline{s}$ and $f(s,c_L)$ and $f(s,c_R)$ are concave functions then the solution of the Riemann problem is given by
\be (s(x,t),c(x,t)) = \left\{ \begin{array}{lll}
(s_L,c_L) &\mbox{if}& x<\sigma_c t, \\ (\overline{s},c_R) &\mbox{if}& \sigma_c t < x < \sigma_s t, \\ (s_R,c_R) &\mbox{if}& x > \sigma_s t, \end{array} \right. 
\ee
where
 $$\sigma_c=\frac{f(s_L,c_L)}{s_L+\bar{a}_L(c_R)}=\frac{f(\overline{s},c_R)}{\overline{s}+\bar{a}_L(c_R)},\quad \sigma_s=\frac{f(\overline{s},c_R) - f(s_R,c_R)}{\overline{s}-s_R}.$$

\noi Note that $ 0 < \sigma_c < \sigma _s$.     
\begin{figure}[H]
 \begin{picture}(320,100)(-40,-10)
 \put(0,0){\includegraphics[height=3cm]{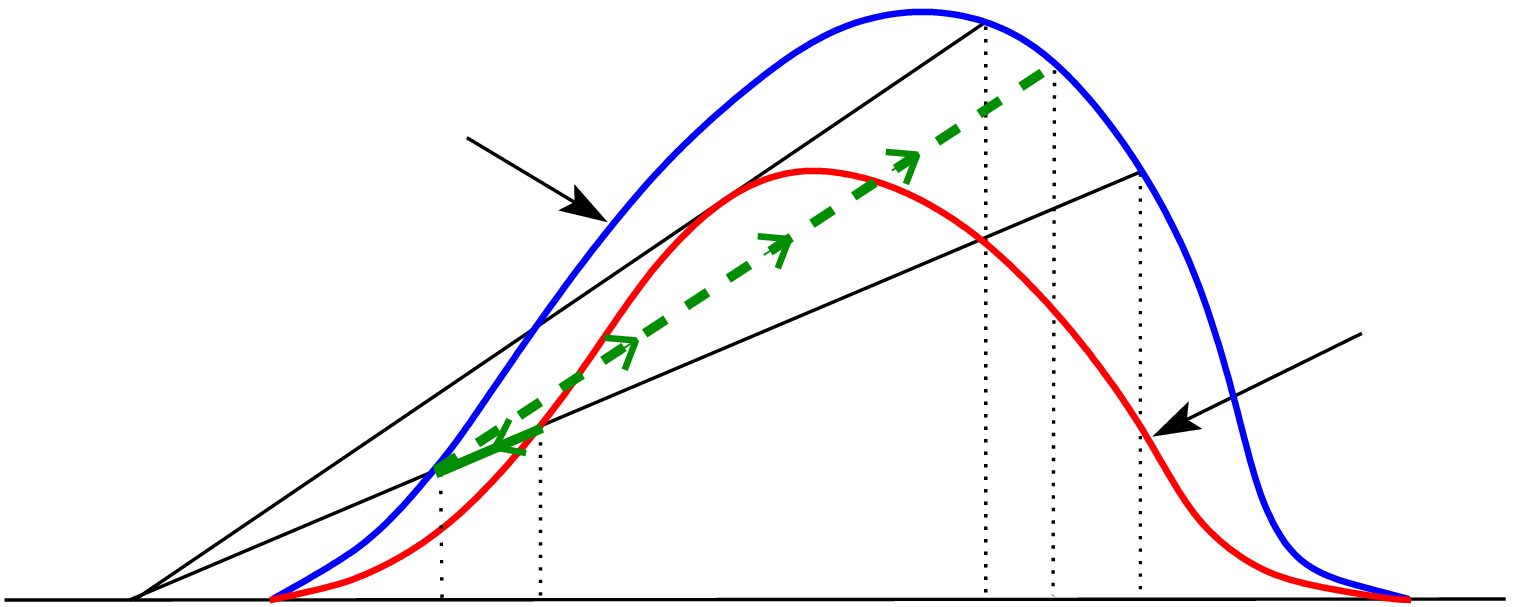}}
 \put(40,72){$f(s,c_R)$}
\put(96,-10){$s^*$}
 \put(190,45){$f(s,c_L)$}
 \put(60,-10){$\overline{s}$}
 \put(70,-10){$s_L$}
 \put(132,-10){$A$}
 \put(144,-10){$s_R$}
 \put(157,-10){$B$}

 \put(-7,-10){$-\bar{a}_L(c_R)$}
 \put(220,0){\includegraphics[height=1.8cm]{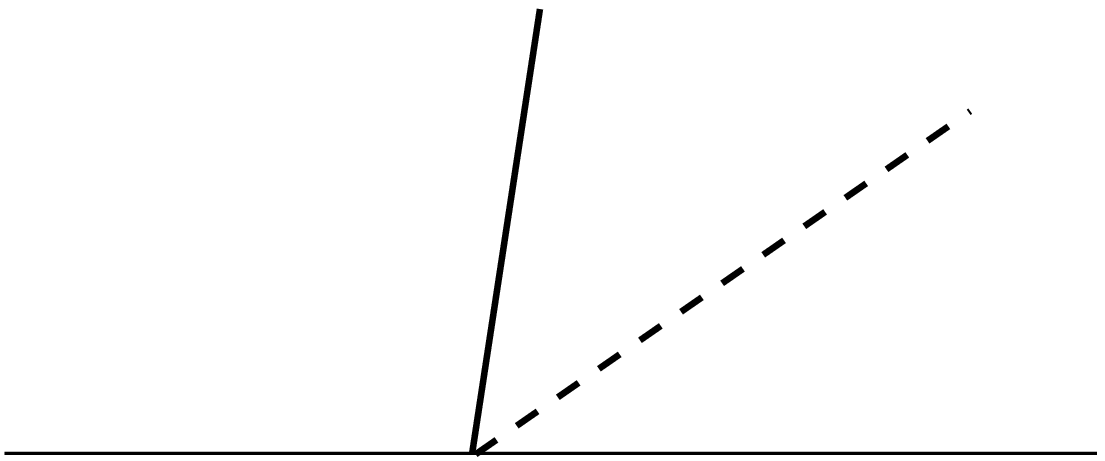}}
  \put(232,15){$(s_L,c_L)$}
 \put(228,-10){$(s_L,c_L)$}
 \put(271,-10){0}
  \put(300,-10){$(s_R,c_R)$}
  \put(310,13){$(s_R,c_R)$}
  \put(283,35){$(\overline{s},c_R)$}
  \put(255,57){$x=\sigma_c t$}
   \put(320,44){$x=\sigma_s t$}
 \end{picture}
 \caption{Solution of Riemann problem (\ref{riemannpb}) with $s_L <  s^*$ and $ s_R <  B$.}
 \label{figR1a}
 \end{figure}
 
 \item {\bf Case 1b:} $s_R \geq B $.

  Draw a line passing through the points $(s_R,f(s_R,c_R))$ and $(-\bar{a}_L(c_R),0)$. This line intersects the curve $f(s,c_L)$ at a point $\overline{s}$ (see Fig. \ref{rs1b}).

 (a)  Connect $(s_L,c_L)$ to $(\overline{s},c_L)$ by a $s$-wave along the curve $f(s,c_L)$.\\
 (b) Next connect  $(\overline{s},c_L)$ to $(s_R,c_R)$ by a $c$-wave with a speed $$\sigma_c=\frac{f(s_R,c_R)}{s_R+\bar{a}_L(c_R)}=\frac{f(\overline{s},c_L)}{\overline{s}+\bar{a}_L(c_R)}. $$
   
For example if $f(s,c_L)$ and $f(s,c_R)$ are concave functions then the solution is given by
\be (s(x,t),c(x,t)) = \left\{ \begin{array}{lll}
(s_L,c_L) &\mbox{if}& x<\sigma_s t, \\ (\overline{s},c_L) &\mbox{if}& \sigma_s
t < x < \sigma_c t \\ (s_R,c_R) &\mbox{if}& x > \sigma_c t \end{array} \right.
\ee
where
 $$\sigma_c=\dfrac{f(s_R,c_R)}{s_R+\bar{a}_L(c_R)}=\dfrac{f(\overline{s},c_L)}{\overline{s}+\bar{a}_L(c_R)}, \quad \sigma_s=\frac{f(\overline{s},c_L) - f(s_L,c_L)}{\overline{s}-s_L}.$$
\noi Note that $  \sigma_s < \sigma _c$ and $(s_L,c_L)$ is connected to
$(\overline{s},c_L)$ by a $s$-shock wave and $(\overline{s},c_L)$ is connected
to $(s_R,c_R)$ by a $c$-shock wave.

\begin{figure}[H]
 \begin{picture}(320,100)(-40,-10)
 \put(0,0){\includegraphics[height=3cm]{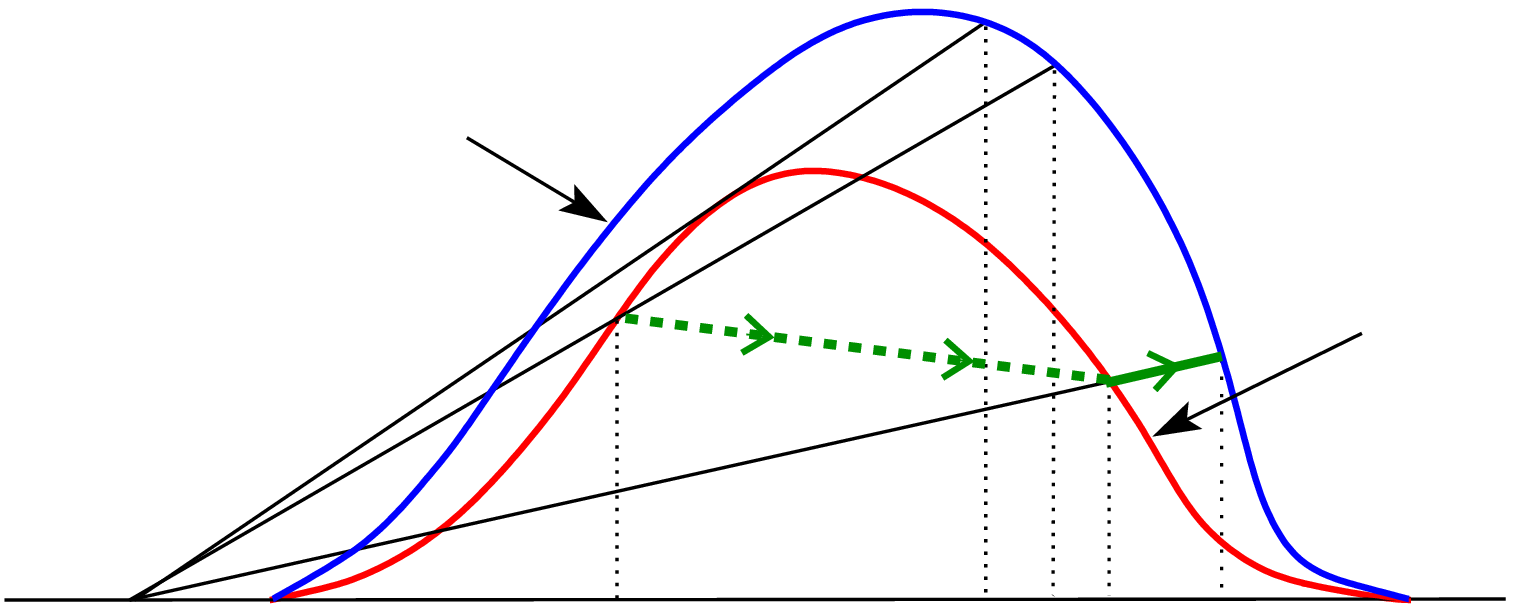}}
 \put(40,72){$f(s,c_R)$}
 \put(190,45){$f(s,c_L)$}
 \put(83,-10){$s_L$}
\put(96,-10){$s^*$}
 \put(132,-10){$A$}
 \put(144,-10){$B$}
 \put(157,-10){$\overline{s}$}
 \put(165,-10){$s_R$}
 \put(-7,-10){$-\bar{a}_L(c_R)$}
 \put(220,0){\includegraphics[height=1.8cm]{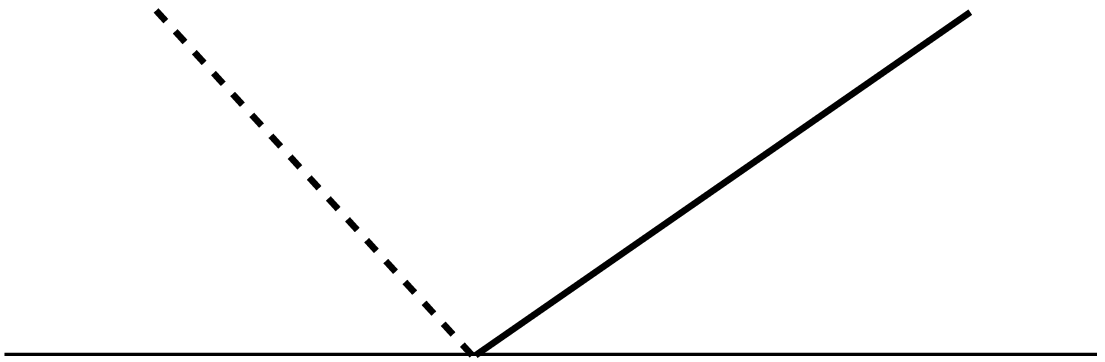}}
 \put(225,15){$(s_L,c_L)$}
 \put(230,-10){$(s_L,c_L)$}
 \put(285,-10){0}
  \put(315,-10){$(s_R,c_R)$}
  \put(325,15){$(s_R,c_R)$}
  \put(283,28){$(\overline{s},c_L)$}
  \put(243,55){$x=\sigma_s t$}
   \put(340,55){$x=\sigma_c t$}
 \end{picture}
 \caption{Solution of Riemann problem (\ref{riemannpb}) with $s_L <  s^*$ and
 $s_R \geq B $.}
 \label{rs1b}
 \end{figure} 
\end{itemize} 
\begin{itemize}
 \item {\bf Case 2:} $ s_L \geq s*$.\\

  \item {\bf Case 2a:} $s_R \leq A$ .\\
(a) Connect $(s_L,c_L)$ to $(s^*,c_L)$ by a $s$-wave along the curve $f(s,c_L)$.\\
(b) Connect $(s^*,c_L)$ to $(\overline{s},c_R)$ by a $c$-wave.\\ 
(c) Connect $(\overline{s},c_R)$ to $(s_R,c_R)$ by a $s$-wave along  the curve $f(s,c_R)$ (see Fig. \ref{rs2a}).

    For example if $s_R \leq \overline{s}$ and $f(s,c_L)$ and $f(s,c_R)$ are concave functions, then the solution is given by 

\[ (s(x,t),c(x,t)) = \left\{ \begin{array}{lll}
(s_L,c_L) &\mbox{if}& x<\sigma_1 t, \\ ((f_s)^{-1}(\frac{x}{t},c_L),c_L) &\mbox{if}& \sigma_1 t < x < \sigma_2 t ,\\ (\overline{s},c_R) &\mbox{if}& \sigma_2 t <x < \sigma_3 t,\\
((f_s)^{-1}(\frac{x}{t},c_R),c_R) &\mbox{if}&  \sigma_3 t < x < \sigma_4 t,\\
(s_R,c_R) &\mbox{if}&  x > \sigma_4 t,
 \end{array} \right.
\]
where
\[ \sigma_1=f_s(s_L,c_L), \quad \sigma_2 = f_s(s^*,c_L)=\frac{f(s^*,c_L)}{s^*+\bar{a}_L(c_R)}, \quad
 \sigma_3 = f_s(\overline{s},c_R), \quad \sigma_4= f_s(s_R,c_R).\]

   Here $ (s_L,c_L)$ is connected to $(s^*,c_L)$ by a $s$-rarefaction wave, $(s^*,c_L)$ is connected to $(\overline{s},c_R)$ by a $c$-shock wave and $(\overline{s},c_R)$ is connected to $(s_R,c_R)$ by a by a rarefaction wave(see Fig. \ref{rs2a}). If $s_R > \overline{s}$ then 
 $(\overline{s},c_R)$ would be connected to $(s_R,c_R)$ by a $s$-chock wave.
\begin{figure}[H]
 \begin{picture}(320,100)(-40,-10)
 \put(0,0){\includegraphics[height=3cm]{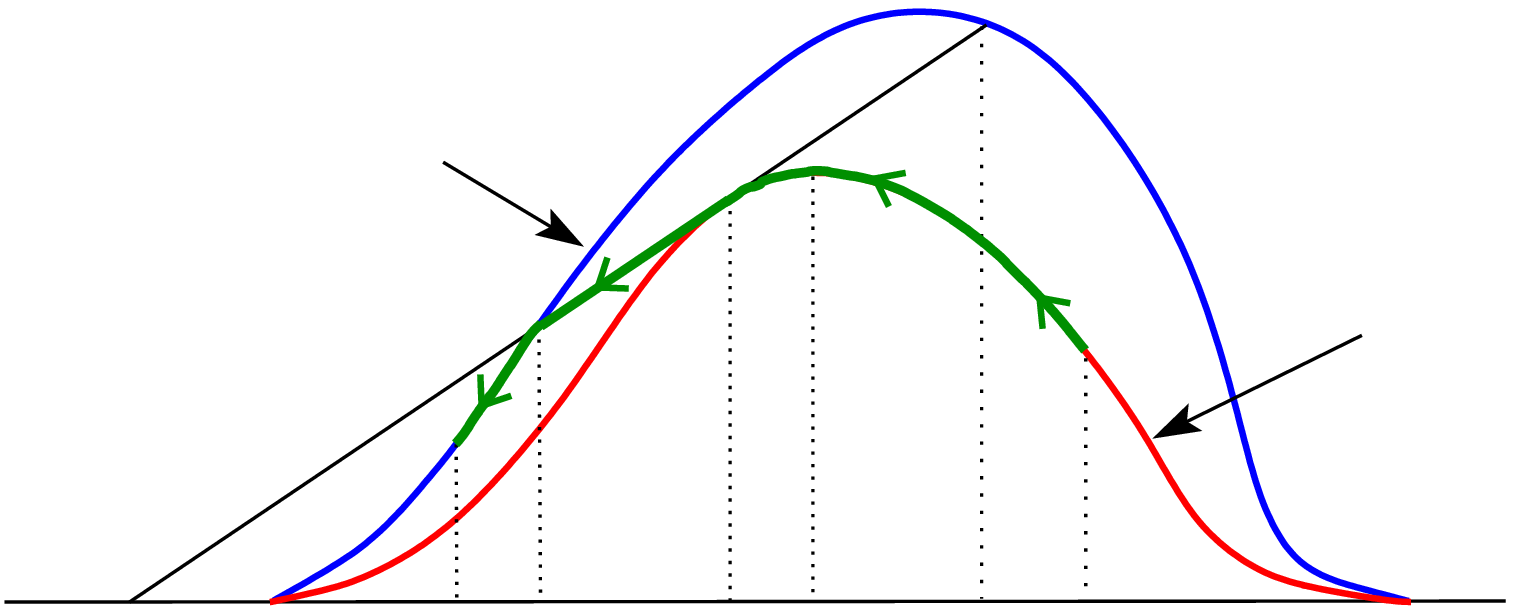}}
 \put(35,70){$f(s,c_R)$}
 \put(190,45){$f(s,c_L)$}
 \put(-8,-10){$-\bar{a}_L(c_R)$}
 \put(58,-10){$s_R$}
  \put(75,-10){$\overline{s}$}
 \put(96,-10){$s^*$}
   \put(111,-10){$\theta_L$}
   \put(220,0){\includegraphics[height=1.8cm]{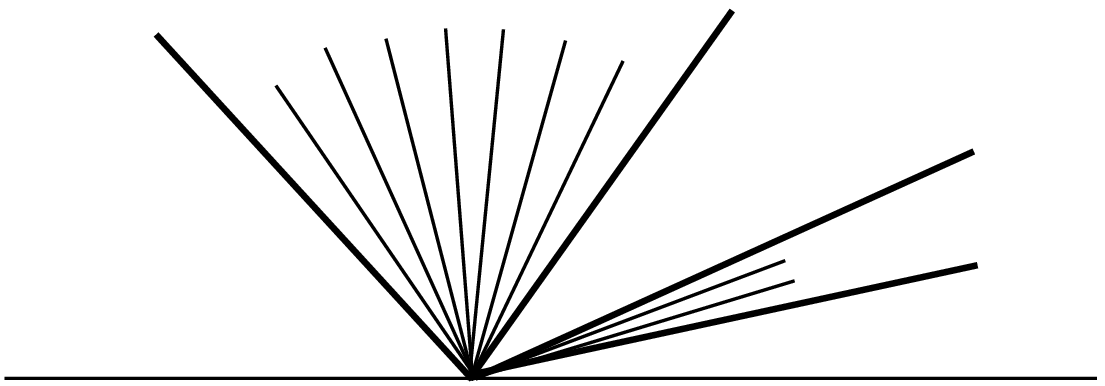}}
 \put(225,10){$(s_L,c_L)$}
 \put(230,-10){$(s_L,c_L)$}
 \put(282,-10){0}
  \put(310,-10){$(s_R,c_R)$}
  \put(330,5){$(s_R,c_R)$}
  \put(310,31){$(\overline{s},c_R)$}
  \put(235,55){$x=\sigma_1 t$}
   \put(300,55){$x=\sigma_2 t$}
   \put(350,35){$x=\sigma_3 t$}
\put(350,17){$x=\sigma_4 t$}
   \put(133,-10){$A$}
    \put(150,-10){$s_L$}
 \end{picture}
 \caption{Solution of Riemann problem (\ref{riemannpb}) with $s_L \geq s^*$ and
 $s_R < A$.}
 \label{rs2a}
 \end{figure}

\item {\bf Case 2b:} $s_R \geq A$\\

 Draw a line passing through the points $(s_R,f(s_R,c_R))$ and $(-\bar{a}_L(c_R),0)$. This line intersects the curve $f(s,c_L)$ at a point $\overline{s}$ (see Fig. \ref{rs2b}).
\begin{figure}[H]
 \begin{picture}(320,100)(-40,-10)
 \put(0,0){\includegraphics[height=3cm]{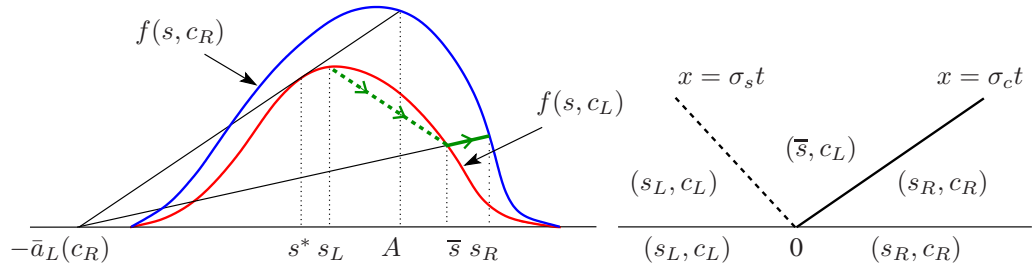}}
\put(40,72){$f(s,c_R)$}
 \put(190,45){$f(s,c_L)$}
 \put(108,-10){$s_L$}
\put(97,-10){$s^*$}
 \put(132,-10){$A$} 
 \put(157,-10){$\overline{s}$}
 \put(165,-10){$s_R$}
 \put(-7,-10){$-\bar{a}_L(c_R)$}
 \put(220,0){\includegraphics[height=1.8cm]{rs1bb.eps}}
 \put(225,15){$(s_L,c_L)$}
 \put(230,-10){$(s_L,c_L)$}
 \put(285,-10){0}
  \put(315,-10){$(s_R,c_R)$}
  \put(325,15){$(s_R,c_R)$}
  \put(283,28){$(\overline{s},c_L)$}
  \put(243,55){$x=\sigma_s t$}
   \put(340,55){$x=\sigma_c t$}
 \end{picture}
 \caption{Solution of Riemann problem (\ref{riemannpb}) with $s_L <  s^*$ and
 $s_R \geq A $.}
 \label{rs2b}
 \end{figure}
                                                                                
 (a)  Connect $(s_L,c_L)$ to $(\overline{s},c_L)$ by a $s$-wave along the curve $f(s,c_L)$,\\
 (b) Next connect  $(\overline{s},c_L)$ to $(s_R,c_R)$ by a $c$-wave with a speed
$$\sigma_c=\frac{f(s_R,c_R)}{s_R+\bar{a}_L(c_R)}=\frac{f(\overline{s},c_L)}{\overline{s}+\bar{a}_L(c_R)}. $$

    For example if $s_L < \overline{s}$ and $f(s,c_L)$ and $f(s,c_R)$ are concave functions, then the   solution is given by
\be 
(s(x,t),c(x,t)) = \left\{ \begin{array}{lll}
(s_L,c_L) &\mbox{if}& x<\sigma_s t, \\ (\overline{s},c_L) &\mbox{if}& \sigma_s
t < x < \sigma_c t, \\ (s_R,c_R) &\mbox{if}& x > \sigma_c t, \end{array} \right.
\ee
where
 $$\sigma_c=\frac{f(s_R,c_R)}{s_R+\bar{a}_L(c_R)}=\frac{f(\overline{s},c_L)}{\overline{s}+\bar{a}_L(c_R)}, \quad \sigma_s=\frac{f(\overline{s},c_L) - f(s_L,c_L)}{\overline{s}-s_L}.$$
\noi Note that $  \sigma_s < \sigma _c$ and $(s_L,c_L)$ is connected to
$(\overline{s},c_L)$ by a $s$-shock wave and $(\overline{s},c_L)$ is connected
to $(s_R,c_R)$ by a $c$-shock wave.
                                                                                
\end{itemize}

\section{Conservative finite volume schemes for the system of polymer flooding}
\label{finite difference}
    Let $h > 0$ and define the space grid points $x_{i+1/2}=ih,i \in \Z$ and
  for $\Delta t >0$ define the time discretization points  $t_n=n \Delta  t $ for all non-negative integer $n$. Let   
$\lambda=\frac{\Delta t}{h}$. A numerical scheme which is in conservative form for equation (\ref{polyeq1}) is given by

\begin{equation}
\begin{array}{l}
 (s_i^{n+1} - s_i^{n}) + \lambda ( F^n_{i+1/2} -  F^n_{i-1/2} )= 0,\\
(c_i^{n+1} s_i^{n+1} + a(c_i^{n+1})- c_i^n s_i^{n} - a(c_i^{n})) + 
\lambda (G^n_{i+1/2} -  G^n_{i-1/2})= 0
\label{finitevolumescheme}
\end{array}
 \end{equation}
 where the numerical flux $F^n_{i+1/2}$ and $ G^n_{i+1/2} $ are associated with the flux functions $f(s,c)$  and $g(s,c)=c f(s,c)$, and are functions of the left and right values of the saturation $s$ and the concentration $c$ at $x_{i+1/2}$:
 \[ F^n_{i+1/2} = F(s_i^n, c_i^n, s_{i+1}^n, c_{i+1}^n), \quad 
 G^n_{i+1/2} = G(s_i^n, c_i^n, s_{i+1}^n, c_{i+1}^n). \]
 The choice of the functions $F$ and $G$ determines the numerical scheme. To recover $c_i^{n+1}$ from the second equation of (\ref{finitevolumescheme}) one has to use an iterative method, like Newton-Raphson. We first present the new flux that we call DFLU, which is constructed as presented in the introduction. We compare it with the exact Riemann solver and show $L^\infty$ estimates for the associate scheme. Then we recall three other schemes to which to compare:  the upstream mobility flux and two centered schemes, Lax-Friedrichs's  and FORCE.
\subsection{The DFLU numerical flux}
\begin{figure}[htbp]   \begin{center}
\begin{picture}(350,70)(0,-5)
\thicklines
\put(70.,20){\parbox{100pt}{$ s_t + f(s, c_i^n)_x = 0$\\ $s(t_n) =
s_i^n$}}
\put(185.,20){\parbox{100pt}{$ s_t + f(s, c_{i+1}^n)_x = 0$\\ $s(t_n) =
s_{i+1}^n$}}
\put(0.,0.){\line(1,0){350}}
\put(0.,50.){\line(1,0){350}}
\put(175.,0.){\line(0,1){50}}
\put(60.,0.){\line(0,1){50}}
\put(300.,0.){\line(0,1){50}}
\put(160.,-10){$x_{i+1/2}$}
\put(45.,-10){$x_{i-1/2}$}
\put(285.,-10){$x_{i+3/2}$}
\put(0.,5){$t=t_n$}
\put(0.,55){$t=t_{n+1}$}
\end{picture}
\end{center}
\caption{The flux functions $f(\cdot,c)$ is discontinuous in c at the
discretization points.}
\label{fluxdiscr}
\end{figure}
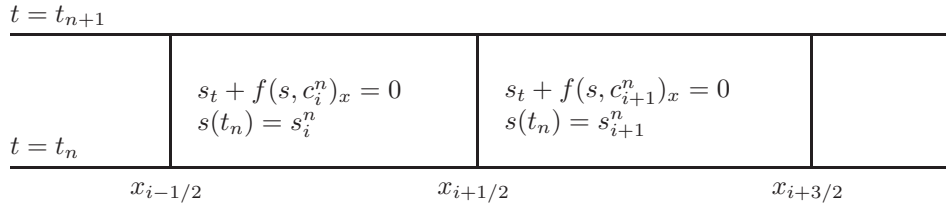

The DFLU flux is an extension of the Godunov scheme that we proposed and analyze in \cite{AdiJafGow04} for scalar conservations laws with a flux function discontinuous in space.
As the second eigenvalue $\lambda^c$ of the system is always non-negative we define
\begin{equation}
G^n_{i+1/2}=c^n_i\,F^n_{i+1/2}. 
 \label{DFLUg} 
\end{equation} 
Now the choice of the numerical scheme depends on the choice
of $F^n_{i+1/2}$. To do so we treat $c(x,t)$
in $f(s,c)$  as  a known function which may be discontinuous at the space discretization points.
Therefore on the border of each rectangle $(x_{i-1/2},x_{i+1/2}) \times (t_n,t_{n+1})$, we
 consider the conservation law:
\begin{equation}
s_t+f(s,c_i^n)_x=0 
\end{equation}
 with initial condition $s(x,0)=s_i^0$ for $x_{i-1/2}< x< x_{i+1/2}$(see 
Fig.\ref{fluxdiscr}).  

Extending the idea of \cite{AdiJafGow04},we define the DFLU flux  as
\be
\begin{array}{lcll}
 F_{i+1/2}^{n} & = & F^{DFLU}(s_i^n, c_i^n, s_{i+1}^n, c_{i+1}^n) \\
             & = & \min\{ f(\min\{s_i^n,\theta_i^n\},c_i^n),f(\max\{s_{i+1}^n,\theta_{i+1}^n\},c^n_{i+1})\},
\end{array} 
\label{DFLU}
\ee
where
$ \theta_i^n =\mbox{ argmax } f(\cdot,c_i^n) $.\\
                                                                         
\noindent {\bf Remarks:} \\
1) Suppose $c_i^n=c_0$, a constant for all $i$,then it is easy to see that $c_i^{n+1}=c_0$ for all $i$.\\
2) Suppose $s \rightarrow f(s,c)$ is an increasing function (case without gravity)  then $\theta_i^n=1$ for all $i$ and from (\ref{DFLU}) we have $ F^n_{i+1/2} = f(s_i^n,c_i^n)$  and the finite difference scheme (\ref{finitevolumescheme}) becomes
\be
\begin{array}{rll}
 s_i^{n+1} &= & s_i^n - \lambda (f(s_i^n,c_i^n) -f(s_{i-1}^n,c_{i-1}^n ))\\
c_i^{n+1} s_i^{n+1} + a(c_i^{n+1})&=& c_i^n s_i^{n} + a(c_i^{n}) -
 \lambda (c_i^n\,f(s_i^n,c_i^n) -c_{i-1}^n\,f(s_{i-1}^n,c_{i-1}^n ))
\label{finitediffeqn2}
\end{array}
\ee
 which  is nothing but the standard upwind scheme.

\subsection{Comparison of the  DFLU flux with the flux given by an exact Riemann solver}
\label{comparison}

     Now we would like to compare the exact Godunov flux $F_{i+1/2}^G $ with our DFLU flux $F^{DFLU}_{i+1/2}$ defined by (\ref{DFLU}). For sake of brevity we considered only the case $c^n_i \geq c^n_{i+1}$. The opposite case can be considered similarly. We discuss the cases considered in section \ref{Riemann}.

\noi {\bf Case 1a:} $ s_i < s^*, s_{i+1} <  B$. See  Fig. \ref{figR1a}. 
   In this case $F_{i+1/2}^G=f(s_i,c_i)=F^{DFLU}_{i+1/2}$.\\

\noi {\bf Case 1b:} $ s_i < s^*, s_{i+1} \geq B$. See  Fig. \ref{rs1b}. \\
   Then 
$F_{i+1/2}^G= \left\{ \begin{array}{lll}
   f(\overline{s},c_i) &\mbox{if}& \sigma_s < 0 \\ f(s_i,c_i) &\mbox{if}& \sigma_s  \geq 0 
\end{array} \right. $
where $ \sigma_s=\dfrac{f(\overline{s},c_i) - f(s_i,c_i)}{\overline{s}-s_i}$. On the other hand the DFLU flux gives   $ F^{DFLU}_{i+1/2}=\min \{f(s_i,c_i),f(\max \{s_{i+1},\theta_{i+1} \},c_{i+1}) \} $.
Therefore in this case the Godunov flux may not be  same as the DFLU flux.\\

\noi {\bf Case 2a:} $s_i \geq s^*,s_{i+1} \leq A $. See Fig.\ref{rs2a}. Then
\[
F_{i+1/2}^G = \left\{ \begin{array}{lll}
   f(\theta_i,c_i) &\mbox{if} &  s_i > \theta_i   \\ f(s_i,c_i) &\mbox{if}& s_i \leq \theta_i
\end{array} \right.
 =f(\min \{s_i,\theta_i \},c_i) =F_{i+1/2}^{DFLU}.
\]
\noi {\bf Case 2b:}$s_i \geq s^*,s_{i+1} > A $. See Fig.\ref{rs2b}. \\
Then $F_{i+1/2}^G= \left\{ \begin{array}{lll}
   f(\overline{s},c_i) &\mbox{if}&  \sigma_s < 0   \\ f(s_i,c_i) &\mbox{if}& \sigma_s \geq 0
\end{array} \right. $
where $ \sigma_s=\dfrac{f(\overline{s},c_i) - f(s_i,c_i)}{\overline{s}-s_i}$. \\
The DFLU flux is $F^{DFLU}_{i+1/2}= \min\{ f(\min\{s_i,\theta_i\},c_i),f(\max\{s_{i+1},\theta_{i+1}\},c_{i+1})\}$.
In this case these two fluxes are not equal, for example when $\sigma_s <0$.

One can actually observe that the Godunov flux can actually be calculated with  the following compact formula:

\noi {\bf Case 1:  $s_i < s^*_i$.}
\[
F_{i+1/2}^G = \!\! \left\{ \begin{array}{lllll}
  \!\! f(s_i,c_i) &\! \mbox{if}\, f_s(s_{i+1},c_{i+1})\geq 0 
  \mbox{ or }\, \dfrac{f(s_{i+1},c_{i+1})}{s_{i+1}+\bar{a}_L(c_{i+1})} \geq \dfrac{f(s_i,c_i)}{s_i+\bar{a}_L(c_{i+1})},\\
   \!\! \min(f(s_i,c_i),f(\overline{s}_i,c_i))&  \!\mbox{otherwise},
\end{array} \right.
\]
where $\overline{s}_i$ is given by
$\dfrac{f(s_{i+1},c_{i+1})}{s_{i+1}+\bar{a}_L(c_{i+1})}=\dfrac{f(\overline{s}_i,c_i)}{\overline{s_i}+\bar{a}_L(c_{i+1})}$.

\noi {\bf Case 2:  $s_i \ge s^*_i$.}
\[
F_{i+1/2}^G = \left\{ \begin{array}{llll}
   \!\! f(\min(s_i,\theta_i),c_i) &\! \mbox{if} \, f_s(s_{i+1},c_{i+1})\geq 0 
   \mbox{ or }\, \dfrac{f(s_{i+1},c_{i+1})}{s_{i+1}+\bar{a}_L(c_{i+1})} \geq \dfrac{f(s^*_i,c_i)}{s^*_i+\bar{a}_L(c_{i+1})},&\\
  \!\!  \min(f(s_i,c_i),f(\overline{s}_i,c_i)) & \! \mbox{otherwise},
\end{array} \right.
\]
where $\overline{s}_i$ is given by
$\dfrac{f(s_{i+1},c_{i+1})}{s_{i+1}+\bar{a}_L(c_{i+1})}=\dfrac{f(\overline{s}_i,c_i)}{\overline{s}_i+\bar{a}_L(c_{i+1})}$.

\subsection{$L^{\infty}$, TV bounds and convergence analysis for the DFLU scheme}
We show first $L^{\infty}$ bounds, and TVD bounds will follow immediately. Let 
$\displaystyle{M=\sup_{s,c}\{|f_s(s,c)|,\frac{f(s,c)}{s+a'(c)}\}}$. 
\begin{lemma}
Let $s_0$ and $c_0 \in  L^\infty (\R, [0,1])$ be the initial data and let
$\{s_i^n\}$ and $\{c_i^n\}$ be the corresponding solution calculated
by the finite volume scheme (\ref{finitevolumescheme}) using the DFLU flux (\ref{DFLUg}),  (\ref{DFLU}).
When  $\la M \leq 1$  then
\be
\begin{array}{l}
0 \leq s_i^n \leq 1 \; \mbox{ for all } i, n,\\
||c^n||_{\infty}  \leq    ||c^{n-1}||_{\infty} \mbox{ where } ||c^n||_{\infty}=\sup_i|c_i^n|.
\end{array}
\label{scinfinity}
\ee
\end{lemma}
{\bf Proof:} Since $0 \leq s_0 \leq 1$ and hence for all $i, \; 0 \leq s_i^0
\leq 1.$ By induction, assume that (\ref{scinfinity}) holds for all $n$. Let 
$$\begin{array}{rll}
 s_i^{n+1} &=& s_i^{n} - \lambda ( F^n_{i+1/2} -  F^n_{i-1/2} )\\
           & =& H(s_{i-1}^n,s_i^n,s_{i+1}^n,c_{i-1}^n,c_i^n,c_{i+1}^n)
\end{array}
$$
By (\ref{DFLU}),it is easy to check that if $\lambda M \leq 1$, then
$H=H(s_1,s_2,s_3,c_1,c_2,c_3)$ is an increasing function in $s_1,s_2,s_3$ and by the hypothesis on $f$, $H(0,0,0,c_1,c_2,c_3)=0, H(1,1,1,c_1,c_2,c_3)=1$. Therefore

$$\begin{array}{rll}
0 &=& H(0,0,0,c_{i-1}^n,c_i^n,c_{i+1}^n) \\
&\leq & H(s_{i-1}^n,s_i^n,s_{i+1}^n,c_{i-1}^n,c_i^n,c_{i+1}^n)=s_i^{n+1}\\
& \leq &  H(1,1,1,c_{i-1}^n,c_i^n,c_{i+1}^n)=1.
\end{array}
$$
This proves $0  \leq s_i^{n+1}  \leq 1$.

     To prove bounds for $c$, consider   
$$ (c_i^{n+1} s_i^{n+1} + a(c_i^{n+1})- c_i^n s_i^{n} - a(c_i^{n})) +
\lambda (G^n_{i+1/2} -  G^n_{i-1/2})= 0.$$
Add and subtract the term $c_i^ns_i^{n+1}$ to the above equation,then we have

$$c_i^{n+1}(s_i^{n+1}+a'(\xi_i^{n+1/2}))-c_i^n (s_i^{n+1}+a'(\xi_i^{n+1/2}))+c_i^n(s_i^{n+1}-s_i^n)+ \lambda (G^n_{i+1/2} -  G^n_{i-1/2})= 0.$$
where $ a(c_i^{n+1})- a(c_i^{n})=a'(\xi_i^{n+1/2}) (c_i^{n+1}-c_i^n)$ for some $\xi_i^{n+1/2}$ between $c_i^{n+1}$ and $c_i^n$. 
Then  substituting for $(s_i^{n+1}-s_i^n)$ from the first equation of (\ref{finitevolumescheme}),since $c_i^nF_{i+1/2}^n=G_{i+1/2}$, we have
 
$$c_i^{n+1}(s_i^{n+1}+a'(\xi_i^{n+1/2}))-c_i^n (s_i^{n+1}+a'(\xi_i^{n+1/2}))+
\lambda  F^n_{i-1/2} (c_i^n-c_{i-1}^n)=0$$.
This is equivalent to 
\be
\begin{array}{llll}
 c_i^{n+1}&=&c_i^n - \lambda \frac{ F^n_{i-1/2}}{(s_i^{n+1}+a'(\xi_i^{n+1/2}))} (c_i^n-c_{i-1}^n)
\end{array}
\label{nonconservative}
\ee
which is the scheme written in the non-conservative form.
Let $b_i^n=\lambda \dfrac{ F^n_{i-1/2}}{(s_i^{n+1}+a'(\xi_i^{n+1/2}))}$ then
\[
 c_i^{n+1} = (1-b_i^n)c_i^n+b_i^nc_{i-1}^n
          \leq \max\{c_i^n,c_{i-1}^n\} \mbox{ if } b_i^n \leq 1.
\]
This proves the second inequality.  \cqfd\\

Since $c_i^{n+1}$ is a convex combination
of $c_i^n$ and $c_{i-1}^n$ if $\lambda M \leq 1$, then we obtain the following total variation diminishing property for $c_i^{n}$:
\begin{lemma}
\label{ctvd}
 Let $\{c_i^n\}$ be the solution  calculated
by the finite volume scheme (\ref{finitevolumescheme}), (\ref{DFLUg}), (\ref{DFLU}).
When  $\la M \leq 1$  then
\[
\sum_i{|c_i^{n+1} - c_{i-1}^{n+1}|}  \leq  \sum_i|c_i^{n} -c_{i-1}^{n}|\,\, \mbox{ for all n}.
\]
\end{lemma}

Also we have from (\ref{nonconservative}) for $\lambda M \leq 1$, 
\be 
\sum_i{|c_i^{n+1} - c_i^n|}  \leq  \sum_i|c_i^{n} -c_{i-1}^{n}|\,\, \mbox{ for all n}.
\label{contraction}
\ee

Note that the saturation $s$  need not be of total variation bounded  because of $f=f(s,c)$ and $c=c(x,t)$ is dicontinuous(see \cite{Adsg11}). The singular mapping technique as in \cite{AdiJafGow04}  to prove the convergence of $\{s_i^n\}$ looks very  difficult to apply. However by using the method of compensated compactness, Karlsen,Mishra and Risebro \cite{KarMisRis09a} showed the convergence of an approximated solution in the case of a triangular system. Now we use their results to prove the convergence of $\{s_i^n\}$. Their method of proof of compensated compactness shows that actually they have proved the following.

   Assume that the flux $f(v,\alpha)$ and the function $k(x,t)$ satisfies the following hypothesis:

\begin{enumerate}
\item $f(0,\alpha)=f(1,\alpha)=0$ for all $\alpha$ in $I$.
\item $f_{vv}(v,\alpha) \neq 0$ for all $\alpha$ in  $I$ and a.e $v$ in $I$
\item There exists $M>0$ and a discretization $\{k_i^n\}$ of $k(x,t)$ exist such that for a subsequence $h$
\begin{enumerate}
\item $\{k_i^n\} \rightarrow k $ in $L^1_{loc}$ as $h \rightarrow 0,$
\item \noi $ \sum_i{|k_i^{n+1} - k_{i-1}^{n+1}|}  \leq  M\;\;\;\mbox{ for all n},$

\item \noi $ \sum_i{|k_i^{n+1} - k_{i}^{n}|}  \leq  M\;\;\;\mbox{ for all n}.$
\end{enumerate} 
\end{enumerate}

\noi Next we describe the discretisation $\{v_i^n\}$ of $v$ corresponding to $\{k_i^n\}$ as follows:
 
   Let $v_{\Delta}^n(x,t)$ be a function defined on the strip $\R \times (n\Delta t,(n+1)\Delta t)$ such that 
\be \left\{\begin{array}{llll}
   (v_{\Delta}^n)_t + f(v_{\Delta}^n,k_i^n)_x=0, \;\;\;(x,t)\;\in\;(x_{i-\frac{1}{2}},x_{i+\frac{1}{2}}) \times  (n\Delta t,(n+1)\Delta t), &&\\ 
v_{\Delta}^n(x,n \Delta t) = v_i^n \;\; \mbox{if}\;\; x\;\in\;(x_{i-\frac{1}{2}},x_{i+\frac{1}{2}}),
\end{array} \right.\ee
$$ f(v_{\Delta}^n(x_{i+\frac{1}{2}}^-,t))= f(v_{\Delta}^n(x_{i+\frac{1}{2}}^+,t    ))\;\; \mbox{for}\; t\;\in\;(n \Delta t,(n+1) \Delta t)$$

and
$$ v_i^{n+1}=\frac{1}{h}\int_{x_{i-\frac{1}{2}}}^{x_{i+\frac{1}{2}}}\,v_{\Delta}^n(\xi,(n+1)\Delta t)\,d\xi. $$

Then we have the following result from  \cite{KarMisRis09a}(see section 5.2).

\begin{lemma}
\label{kmb}
Assume that $v_\Delta^n$ satisfies
\begin{enumerate}
\item $0 \leq \sup_i |v_i^n| \leq 1.$
\item $v_\Delta^n$ satisfies "minimal jump condition" at each interface $x_{i+\frac{1}{2}}$. 
\end{enumerate}
Then there exists  subsequences of $\{k_i^n\}$ and $\{v_i^n\}$ converges   respectively to $k$ and $v\; a.e$  and these limits are the  solution of

\be \left\{ \begin{array}{llll}
   v_t + f(v,k)_x=0&& \\
   v(x,0)=v_0(x)&&
\end{array}
\right.\ee 
\end{lemma}
{ Proof of convergence of $\{s_i^n\}$:} Assume further that $c_0$ and $f$ satisfies the following.

\begin{enumerate}
\item[(i)] $c(x,0)=c_0(x)$ is of bounded variation.
\item[(ii)] $c \rightarrow f(s,c)$ is a non-increasing function.
\item[(iii)] $f_{ss}(s,c)\neq 0$ for all $c$ and a.e $s$.
\end{enumerate}

    Let $\{c_i^n\}$ be as in Lemma \ref{ctvd} and $s_i^n$ be the  corresponding  solution obtained from DFLU flux (\ref{DFLU}). Then it follows from the above hypothesis (ii), $s_i^n$ satisfies  the "minimal jump condition" across the interface. Hence by taking
$$c_i^n=k_i^n\;\;\mbox{and}\;\;s_i^n=v_i^n,$$
it follows from (\ref{contraction}) and Lemmas \ref{ctvd},\ref{kmb}, there exists  subsequences
 of $c_i^n$ and $ s_i^n$ converges respectively to $c$ and $s$. Further  $s$ satisfies
$$     s_t+f(s,c)_x=0$$.

\noi {\bf Remark:} As equation (\ref{nonconservative}) for $c$ is in non-conservative form, though the sequence $\{c_i^n\}$ is $L^{\infty}$ stable and TVD, it is difficult to prove the convergence$\{s_i^nc_i^n\}$ to a weak solution of $(sc+a(c))_t+(cf(s,c))_x=0$ unless, like in \cite{TveWin91b,TveWin91a}, the concentration $c$ is Lipschitz continuous or like in \cite{TveWin90} fluxes are in the special form. In the presence of viscosity, the convergence of the Lax-Friedrichs scheme for the polymer flooding model was proved in \cite{Tveito90}.

\subsection{The upstream mobility flux}

 Petroleum engineers have designed, from physical
considerations, another numerical flux
called the upstream mobility flux. It is an ad-hoc flux for two-phase
flow in porous media which
corresponds to an approximate solution to the Riemann problem. For this flux $G^n_{i+1/2}$ is given again by (\ref{DFLUg}) and $F_{i+\frac{1}{2}}^n$ is given by
\[
\begin{array}{l}
F_{i+\frac{1}{2}}^n=F^{UM}(s_i^n,c_i^n,s_{i+1}^n,c_{i+1}^n)= \displaystyle{
         \frac{\lambda_1^*}{\lambda_1^* + \lambda_2^*}
         [ \var + (g_1-g_2)\lambda_2^* ]} ,\\
\lambda^*_\ell=
\left\{ \begin{array}{ll} \lambda_\ell(s^n_i,c_i^n) & \mbox{if }
         \var +(g_\ell-g_i)\lambda_i^* >0, \;i=1,2, i\neq \ell,\\[3mm]
        \lambda_\ell(s_{i+1}^n,c_{i+1}^n) & \mbox{if }
        \var +(g_\ell-g_i)\lambda_i^* \leq0, \; i=1,2, i\neq \ell,
\end{array}
\right. \ell=1,2.
\end{array}
\]                                                               

\subsection{The Lax-Friedrichs flux}
In this case fluxes are given by
\[
\begin{array}{llll}
F^n_{i+1/2}&=& \frac{1}{2} [f(s^n_{i+1},c^n_{i+1})+f(s^n_i,c^n_i)-\dfrac{(s^n_{i+1}-s^n_i)}{\lambda} ]\\
G^n_{i+1/2}&=&\frac{1}{2} [c^n_{i+1}f(s^n_{i+1},c^n_{i+1})+c^n_i
f(s^n_i,c^n_i)-\dfrac {(c^n_{i+1} s^n_{i+1}+a(c^n_{i+1})-c^n_is^n_i-a(c^n_i))}
{\lambda} ]
\end{array}
\]

\subsection{The FORCE  flux}
This flux \cite{Toro99,AdiGowJaf09}, introduced by E. F. Toro, is an average of the Lax-Friedrichs and Lax-Wendroff flux. It is defined by
\[
\begin{array}{llll}
F^n_{i+1/2}&=& \frac{1}{4} [f(s^n_{i+1},c^n_{i+1})+f(s^n_i,c^n_i)+2 f(s^{n+1/2}_i,c^{n+1/2}_i)-\dfrac{(s^n_{i+1}-s^n_i)}{\lambda} ]\\[0.3cm]
G^n_{i+1/2}&=&\frac{1}{4} [c^n_{i+1}f(s^n_{i+1},c^n_{i+1})+c^n_i
f(s^n_i,c^n_i)+2 c^{n+1/2}_i f(s^{n+1/2}_i,c^{n+1/2}_i)\\[0.2cm]
&& -\dfrac {(c^n_{i+1} s^n_{i+1}+a(c^n_{i+1})-c^n_i s^n_i-a(c^n_i))}{\lambda}]\\
\end{array}
\]
 where
$$s^{n+1/2}_i=\frac{(s^n_{i+1}+s^n_i)}{2}-\frac {\lambda}{2}(f(s^n_{i+1},c^n_{i+1})-f(s^n_i,c^n_i))$$
and
\[
\begin{array}{llll}
 s^{n+1/2}_ic^{n+1/2}_i+a(c^{n+1/2}_i)&=&\dfrac{(s^n_{i+1}c^n_{i+1}+s^n_ic^n_i)}{2}+\dfrac{(a(c_{i+1}^n)+a(c_i^n))}{2}\\[0.3cm]
&&\qquad \qquad -\frac{\lambda}{2}(c^n_{i+1}f(s^n_{i+1},c^n_{i+1})-c^n_if(s^n_i,c^n_i)).
\end{array}
\]

\section{Numerical experiments}
\label{numericalresults}
  To evaluate the performance of the DFLU scheme we first compare its results to an exact solution and evaluate convergence rates, and then  compare it with other standard numerical
 schemes already mentioned in the previous section, that are the Godunov, upstream mobility, Lax-Friedrichs and FORCE schemes.  
 \subsection{Comparison with an exact solution}
In this section we compare the calculated and exact solutions of two Riemann problems. We consider the following functions 
\be  f(s,c)=s(4-s)/(1+c),  \quad a(c)=c. \label{flux1} \ee
Note that $f(0,c)=f(4,c)=0$ for all $c$ and the interval for $s$ is $[0,4]$ instead of $[0,1]$. This choice of $f$, which does not correspond to any physical reality, was done in order to try to have a large difference between the Godunov and the DFLU flux (see second experiment below).

In a first experiment the initial condition is
\be s(x,0) = \left\{ \begin{array}{lll}
2.5 &\mbox{if}& x<.5, \\ 1 &\mbox{if}& x>.5 \end{array} \right.  , \quad
c(x,0) = \left\{ \begin{array}{lll}
.5 &\mbox{if}& x<.5, \\ 0 &\mbox{if}& x>.5. \end{array} \right.  
\label{IC1} \ee
  These $f$ and initial data correspond to 
   the case 2a in sections \ref{Riemann} and \ref{comparison} where the DFLU flux coincides with the Godunov flux: $F^{DFLU}(s_L,s_R,c_L,c_R)=F^{G}(s_L,s_R,c_L,c_R)$ with 
   $s^*=1.236, A=2.587,\overline{s}=.394$.  The exact  solution of the Riemann problem at a time $t$ is given by
\be s(x,t) = \left\{ \begin{array}{lll}
2.5 &\! \mbox{if}&\! x<.5+\sigma_1 \,t\, \\ \frac{1}{2}(4-1.5(\frac{x-.5}{t}))&\!\mbox{if}&\! .5+\sigma_1 \,t <x < .5 +\sigma_c
\, t\\  \overline{s}=.394 &\!\mbox{if}&\! .5+\sigma_c \,t < x <.5+\sigma_2 \,t \\
1. & \!\mbox{if} &\! x > \sigma_2 t+.5
  \end{array} \right.  , \;\;
c(x,t) = \left\{ \begin{array}{lll}
.5 & \,\,\,  \!\!\! \mbox{if}& \,\,\, \!\!\!  x<.5+\sigma_c t, \\ 0. &\,\,\,  \!\!\! \mbox{if}& \,\,\,  \!\!\!  x>.5+ \sigma_c t. \end{array} \right.  
\ee
where $\sigma_1=f_s(s_L,c_L)=-2/3$,
 $\sigma_c=f_s(s^*,c_L)=\dfrac{f(s^*,c_L)}{s^*+\bar{a}_L(c_R)}=\dfrac{f(\overline{s},c_R)}{\overline{s}+\bar{a}_L(c_R)}=1.018$ and $ \sigma_2=\dfrac{f(\overline{s},c_R) - f(s_R,c_R)}{\overline{s}-s_R}=2.606$.

\begin{figure}[H]
 \includegraphics[width=5cm,angle=-90]{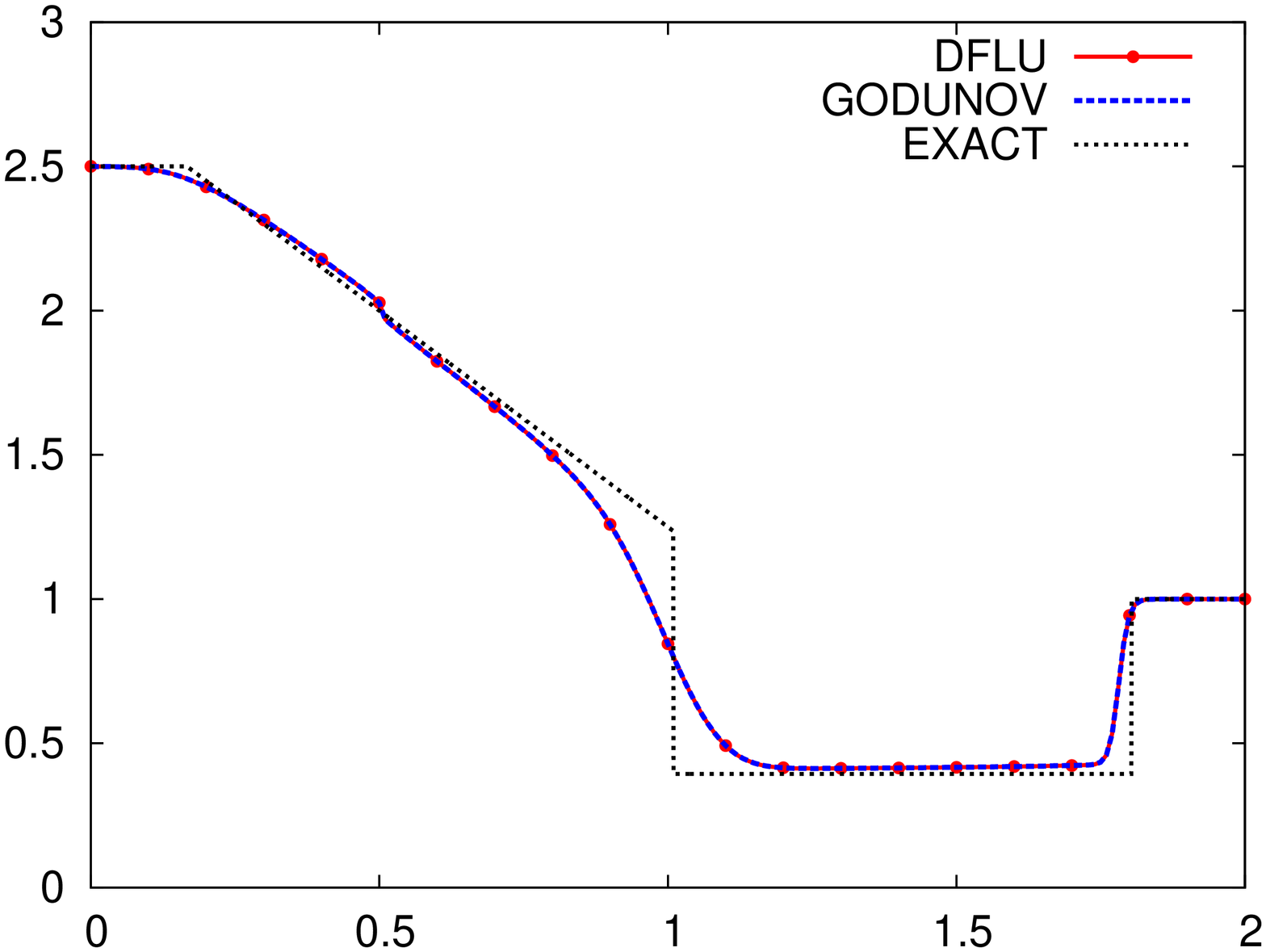} \hspace*{0.3cm}
 \includegraphics[width=5cm,angle=-90]{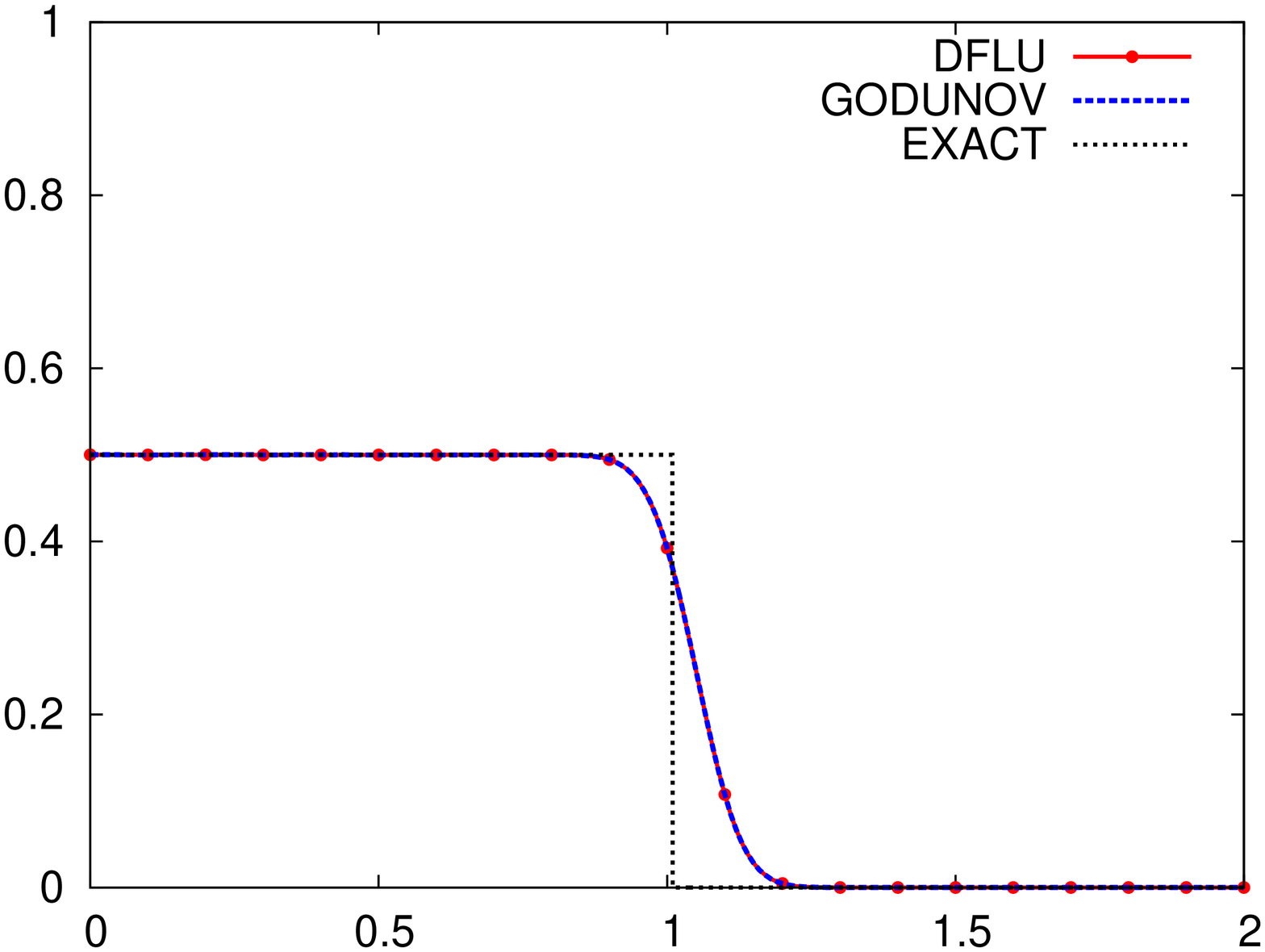}
\caption{Comparison with exact solution of Riemann problem (\ref{flux1}), (\ref{IC1}):  $s$ (left) and $c$ (right) at $t=.5$ for $h=1/100, \lambda=1/4$.} 
\label{fig2aexp1}
\end{figure}

\begin{figure}[H]
 \includegraphics[width=5cm,angle=-90]{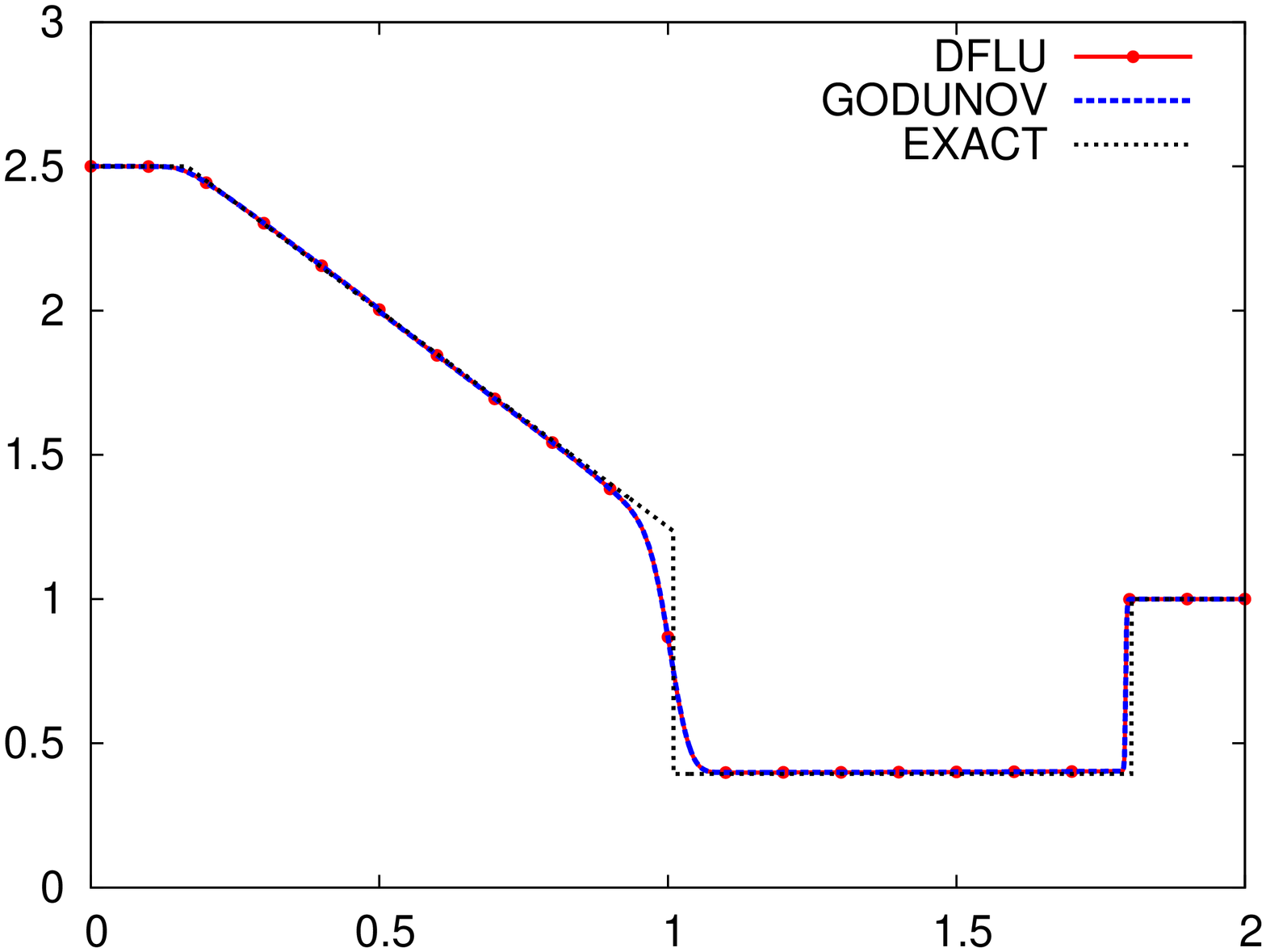} \hspace*{0.3cm}
 \includegraphics[width=5cm,angle=-90]{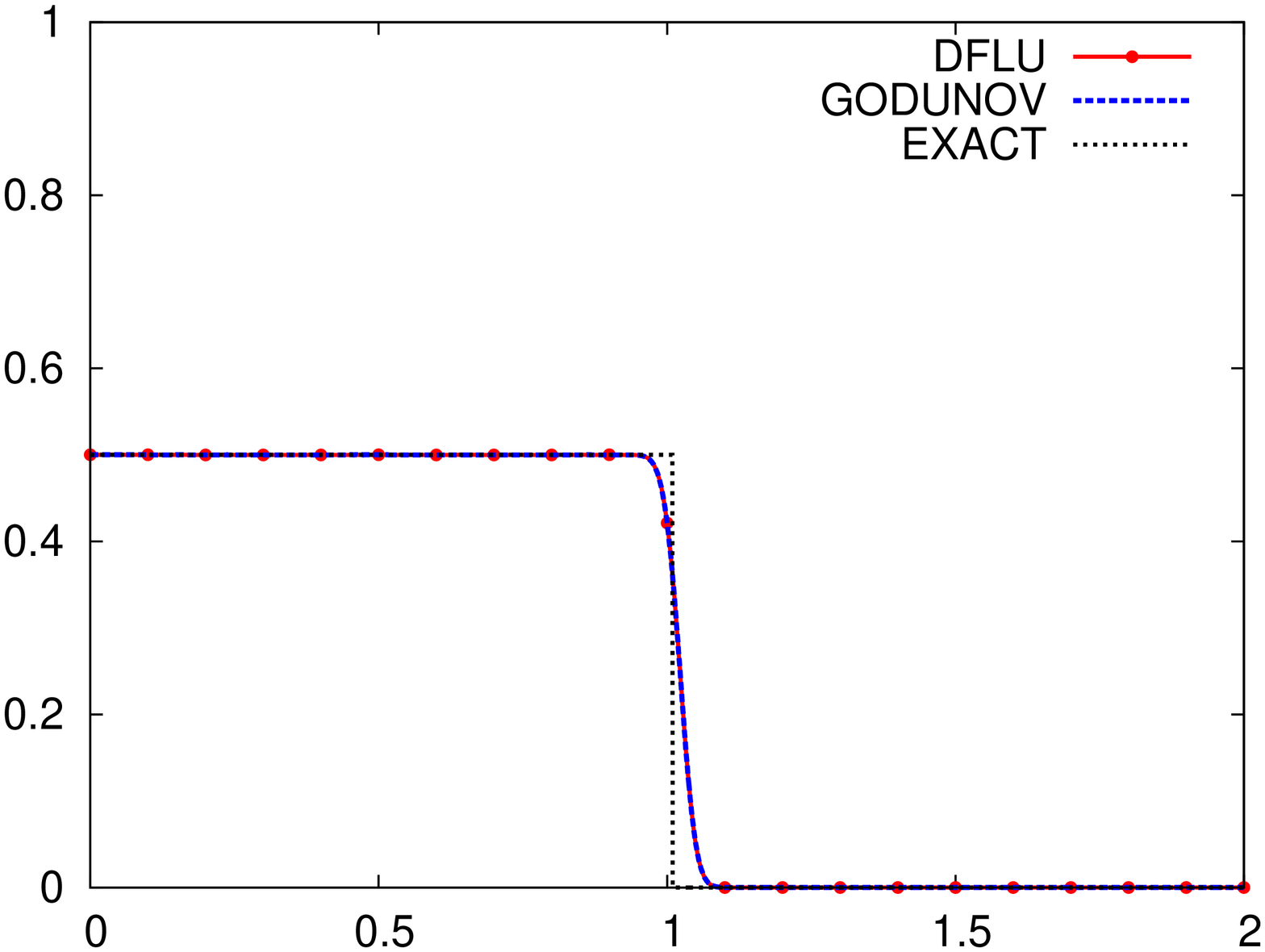}
\caption{Comparison with exact solution of Riemann problem (\ref{flux1}), (\ref{IC1}):  $s$ (left) and $c$ (right) at $t=.5$ for $h=1/800, \lambda=1/4$.} 
\label{fig2aexp11}
\end{figure}

 Figs. \ref{fig2aexp1} and \ref{fig2aexp11} verify that the DFLU and Godunov schemes give coinciding results.  
 As expected both  schemes are diffusive at $c$-shocks as well as at $s$-shocks but
 as the mesh size goes to zero calculated solutions are getting closer  
 to the exact solution (see Fig.\ref{fig2aexp11}).  Table \ref{table1} shows $L_1$ errors for $s$ and  $c$ and the convergence rate $\alpha$. Calculations are done with $\lambda=\frac{1}{4} (M=4)$, that is the largest time step allowed by the CFL condition.
\begin{table}[hbtp]
\begin{center}
\begin{tabular}{|l||c|c||c|c|}\hline
$h$& Godunov,$||s-s_h||_{L^1}$&$\alpha$&DFLU,$||s-s_h||_{L^1}$& $\alpha$ \\
\hline
1/50&.2373&&.2372&\\
\hline
1/100&0.15134&0.6489&0.1506&0.655\\
\hline
1/200& 9.6868 $\times 10^{-2}$&0.6437&9.6868 $\times 10^{-2}$&0.6366  \\
\hline
1/400&6.4228 $\times 10^{-2}$&0.5928&6.4228 $\times 10^{-2}$ &  0.5928\\
\hline
1/800&4.2198 $\times 10^{-2}$&0.606&4.2197 $\times 10^{-2}$ &0.606  \\
\hline
\end{tabular}\\[0.3cm]
\begin{tabular}{|l||c|c||c|c|}\hline
$h$& Godunov,$||c-c_h||_{L^1}$&$\alpha$&DFLU,$||c-c_h||_{L^1}$& $\alpha$\\
\hline
1/50& 6.3796 $\times 10^{-2}$&&6.3796 $\times 10^{-2}$&\\
\hline
1/100&4.1630 $\times 10^{-2}$&0.6158& 4.1630 $\times 10^{-2}$&0.6158 \\
\hline
1/200& 2.6669 $\times 10^{-2}$&0.6424&2.6669 $\times 10^{-2}$ &0.6424 \\
\hline
1/400&1.7398 $\times 10^{-2}$&0.6162&1.7398 $\times 10^{-2}$ &0.6162  \\
\hline
1/800&1.1522 $\times 10^{-2}$&0.5945&1.1522 $\times 10^{-2}$&0.5945  \\
\hline
\end{tabular}
\end{center}
\caption{ Riemann problem (\ref{flux1}), (\ref{IC1}): $L^1$-errors between exact and  calculated solutions at $t=.5$}
\label{table1}
\end{table}

   Now we want to have an experiment where the DFLU flux differs from the Godunov flux. Therefore we now consider the Riemann problem with initial data
\be s(x,0) = \left\{ \begin{array}{lll}
2.3 &\mbox{if}& x<.5, \\ 3.2 &\mbox{if}& x>.5, \end{array} \right.  , \quad
c(x,0) = \left\{ \begin{array}{lll}
.5 &\mbox{if}& x<.5, \\ 0 &\mbox{if}& x>.5. \end{array} \right.  
\label{IC2}
\ee
This initial data corresponds to case 2b of  sections \ref{Riemann} and \ref{comparison} with
 $c_R=0$, $s^*=1.236$.
     In this case, the  exact solution of  the Riemann problem at a time $t$ is given by
\[ s(x,t) = \left\{ \begin{array}{lll}
s_L=2.3 &\mbox{if}& x<.5+\sigma_s \,t \\
 \overline{s}=2.7536&\mbox{if}& .5+\sigma_s\,t <x < .5 +\sigma_c t,\\
s_R=3.2 &\mbox{if} & x > \sigma_c t+.5
  \end{array} \right.  , \quad
c(x,0) = \left\{ \begin{array}{lll}
.5 &\mbox{if}& x<.5+\sigma_c t, \\ 0. &\mbox{if}& x>.5+ \sigma_c t, \end{array} \right.  
\]
where $\sigma_s=\dfrac{f(s_L,c_L)-f(\overline{s},c_L)}{s_L-\overline{s}}=-.702$,
 and  $\sigma_c=\dfrac{f(s_R,c_R)}{s_R+\bar{a}_L(c_R)}=0.609$.       
 
\begin{figure}[H]
 \includegraphics[width=5cm,angle=-90]{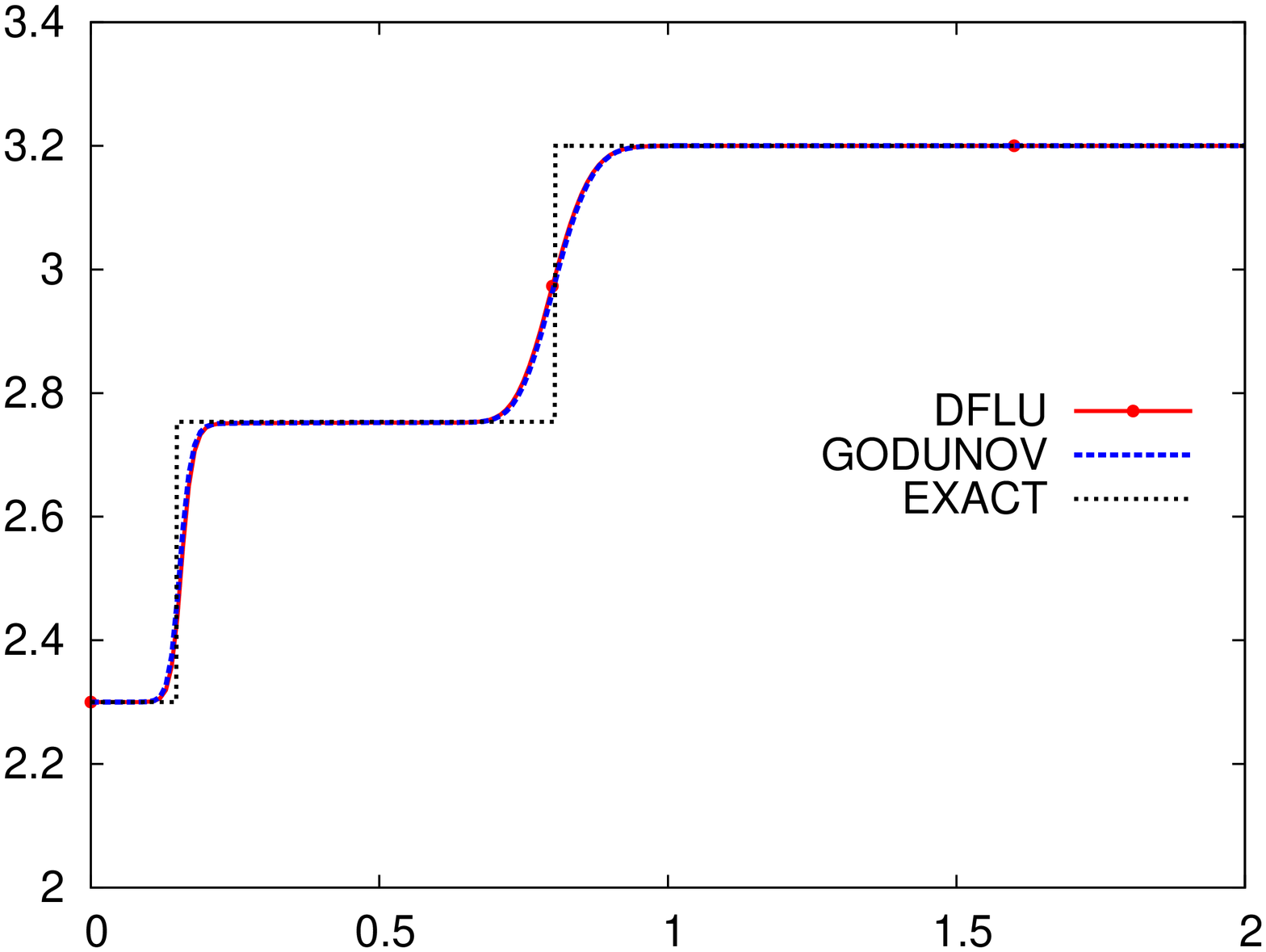}
 \includegraphics[width=5cm,angle=-90]{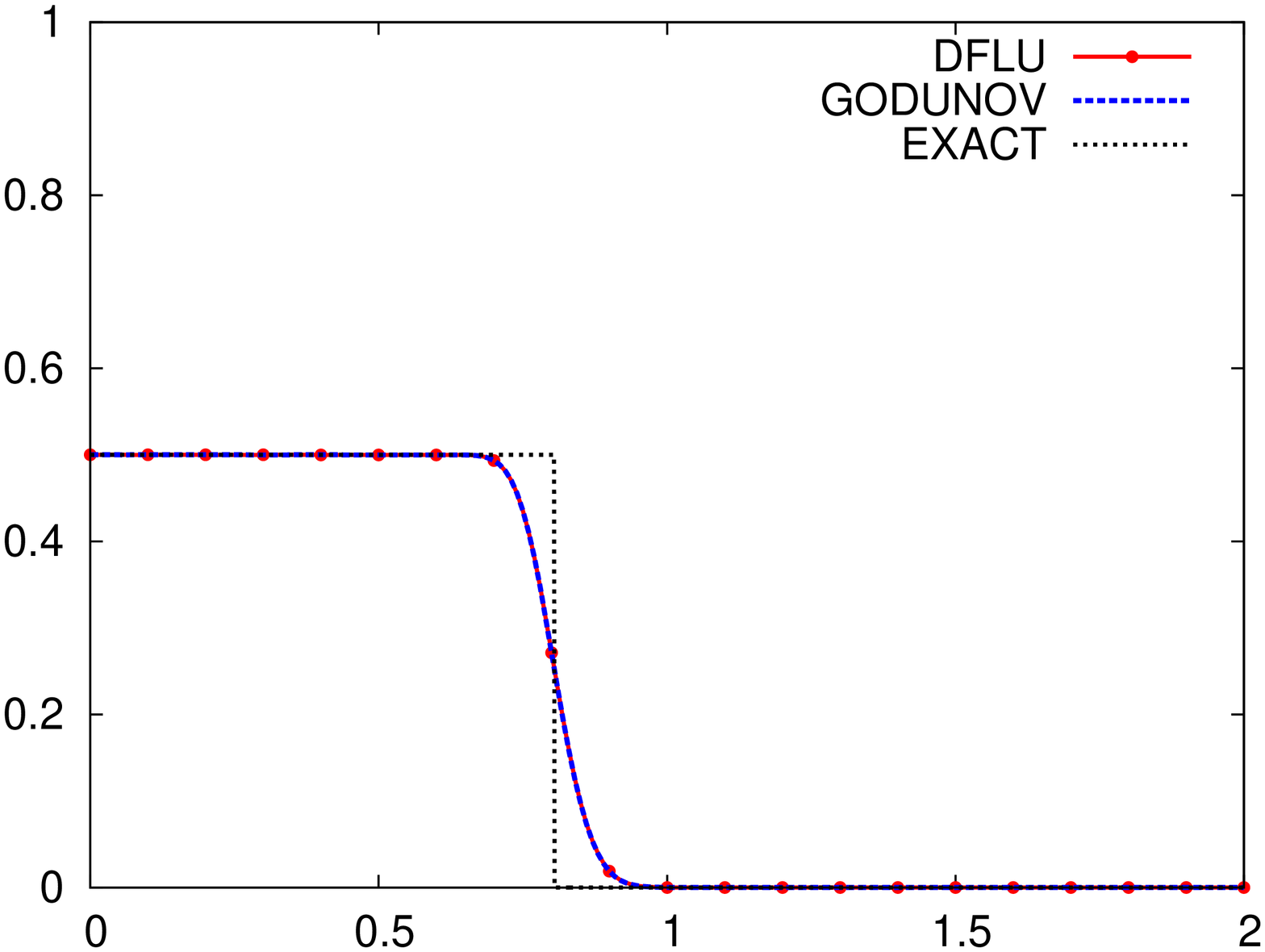}
\caption{Comparison with exact solution of Riemann problem (\ref{flux1}), (\ref{IC2}):  $s$ (left) and $c$ (right) at $t=.5$ for $h=1/100, \lambda=1/4$.} 
\label{fig2bexp2}
\end{figure}

\begin{figure}[H]
 \includegraphics[width=5cm,angle=-90]{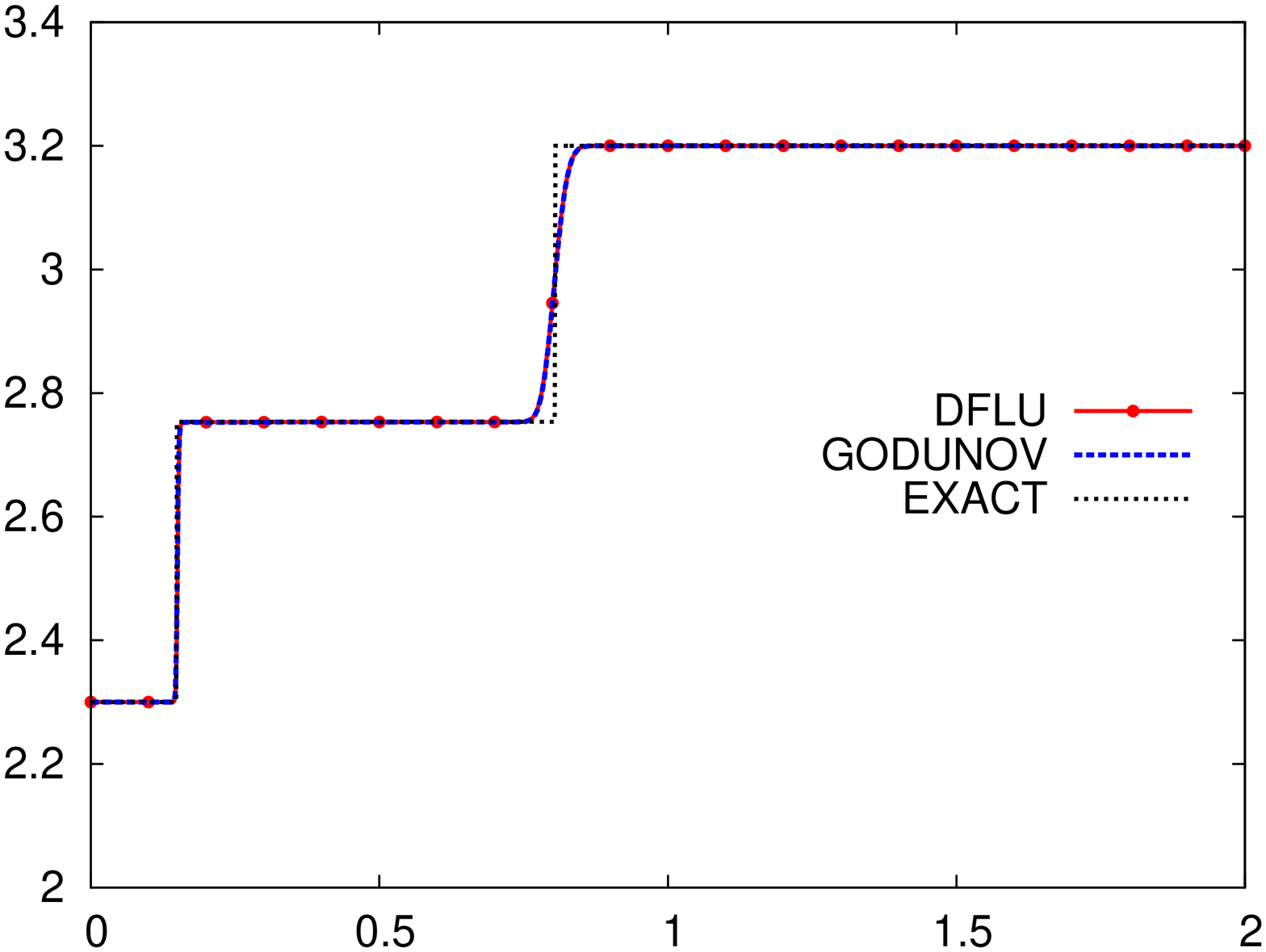}
 \includegraphics[width=5cm,angle=-90]{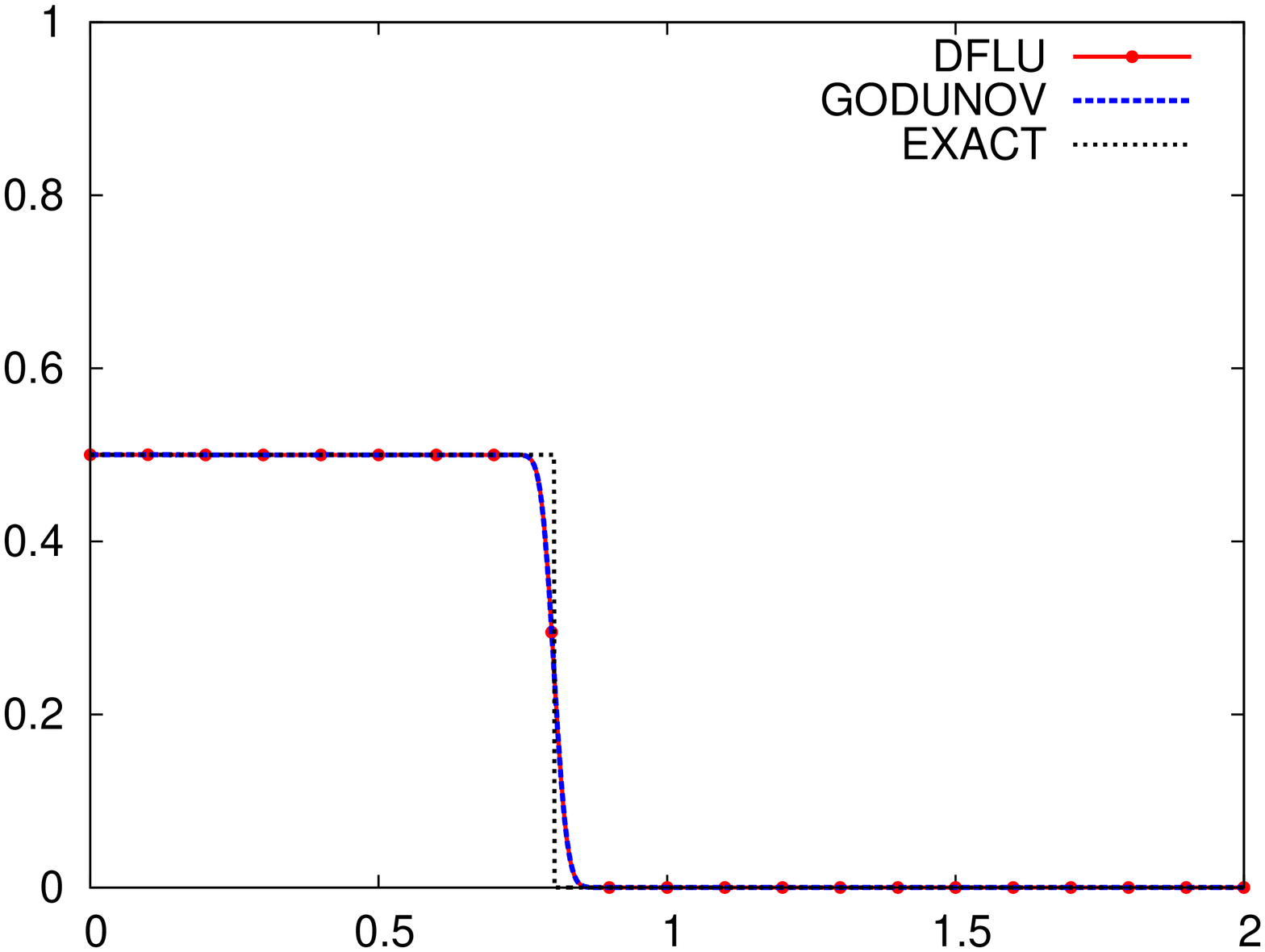}
\caption{ Comparison with exact solution of Riemann problem (\ref{flux1}), (\ref{IC2}):  $s$ (left) and $c$ (right) at $t=.5$ for $h=1/800, \lambda=1/4$.} 
\label{fig2bexp21}
\end{figure}

Figs. \ref{fig2bexp2} and \ref{fig2bexp21} show the comparison of the results obtained with the DFLU and Godunov fluxes with the exact solution. The solution obtained with the DFLU and Godunov flux are very close even if they do not coincide actually. Table \ref{table2} shows $L_1$ errors for $s$ and  $c$ and the convergence rate $\alpha$. Calculations are done with $\lambda=\frac{1}{4} (M=4)$, that is the largest time step allowed by the CFL condition.

\begin{table}[H]
\begin{center}
\begin{tabular}{|l||c|c||c|c|}\hline
$h$& Godunov,$||s-s_h||_{L^1}$&$\alpha$&DFLU,$||s-s_h||_{L^1}$&$\alpha$\\
\hline
1/50&0.10246&&0.10373&\\
\hline
1/100&5.7861 $\times 10^{-2}$ &0.8243 &5.8731 $\times 10^{-2}$&0.8206\\
\hline
1/200& 3.2849 $\times 10^{-2}$&0.81674&3.3259 $\times 10^{-2}$& 0.8203 \\
\hline
1/400&1.9152 $\times 10^{-2}$&0.7785&1.9353 $\times 10^{-2}$ & 0.7811 \\
\hline
1/800&1.1489 $\times 10^{-2}$&0.7370&1.1571 $\times 10^{-2}$ &  0.7420\\
\hline
\end{tabular}\\[0.3cm]
\begin{tabular}{|l||c|c||c|c|}\hline
$h$& Godunov,$||c-c_h||_{L^1}$&$\alpha$&DFLU,$||c-c_h||_{L^1}$& $\alpha$\\
\hline
1/50&4.8407 $\times 10^{-2}$&& 4.8486 $\times 10^{-2}$&\\
\hline
1/100&3.0161 $\times 10^{-2}$&0.6825& 3.0201 $\times 10^{-2}$&0.6829\\
\hline
1/200& 1.9307 $\times 10^{-2}$&0.6435&1.9328$\times 10^{-2}$ &0.6439  \\
\hline
1/400&1.2618 $\times 10^{-2}$&0.6136&1.2628 $\times 10^{-2}$ &0.6140  \\
\hline
1/800&8.4125$\times 10^{-3}$&0.5848&8.4173 $\times 10^{-3}$&0.5851   \\
\hline
\end{tabular}
\end{center}
\caption{Riemann problem (\ref{flux1}), (\ref{IC2}): $L^1$-errors between exact and calculated solutions at $t=.5$.}
\label{table2}
\end{table}

\subsection{Comparison of the DFLU, upstream mobility, FORCE and Lax-Friedrichs fluxes}
   In the previous section, we have seen that Godunov and DFLU fluxes give schemes with very close performances. In this section we compare the DFLU flux with the other fluxes that we mentioned in section \ref{finite difference} which are the
upstream mobility, FORCE and Lax-Friedrichs fluxes. We take now 
\begin{equation} \begin{array}{l} 
f(s,c)=  \var_1 = \dfrac{\lambda_1(s,c)}{\lambda_1(s,c) + \lambda_2(s,c)} [ \var + (g_1-g_2)\lambda_2(s,c) ],\\
\lambda_1(s,c)=\dfrac{s^2}{.5+c}, \lambda_2(s,c)=(1-s)^2, \,\,g_1=2, g_2=1,  \varphi=0,\\
a(c)=.25c. 
\end{array} \label{f2} \end{equation}
In all following experiments the discretization is such that $\Delta t=1/125$ and $h=1/100$.

\noi {\bf Remark:} Even for a total Darcy velocity $\var \neq 0$, the DFLU scheme works. For the DFLU scheme to work, what one needs is $f(0,c)=c_1$ for all $c \in I$ and $f(1,c)=c_2$ for all $c \in I$, for some constants $c_1$ and $c_2$.

\begin{figure}[H]
\includegraphics[width=7cm]{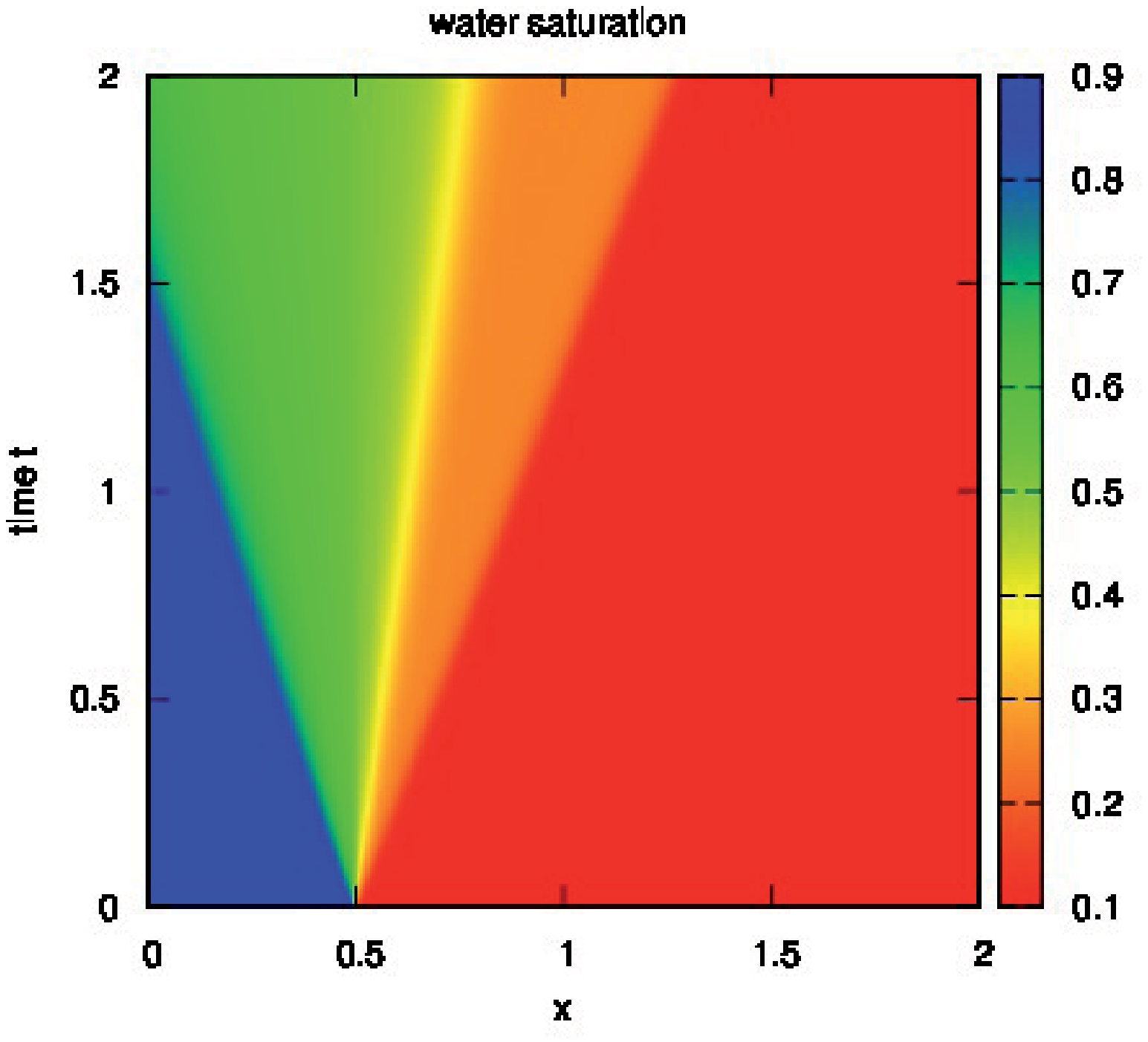} \hspace*{0.5cm}
\includegraphics[width=7cm]{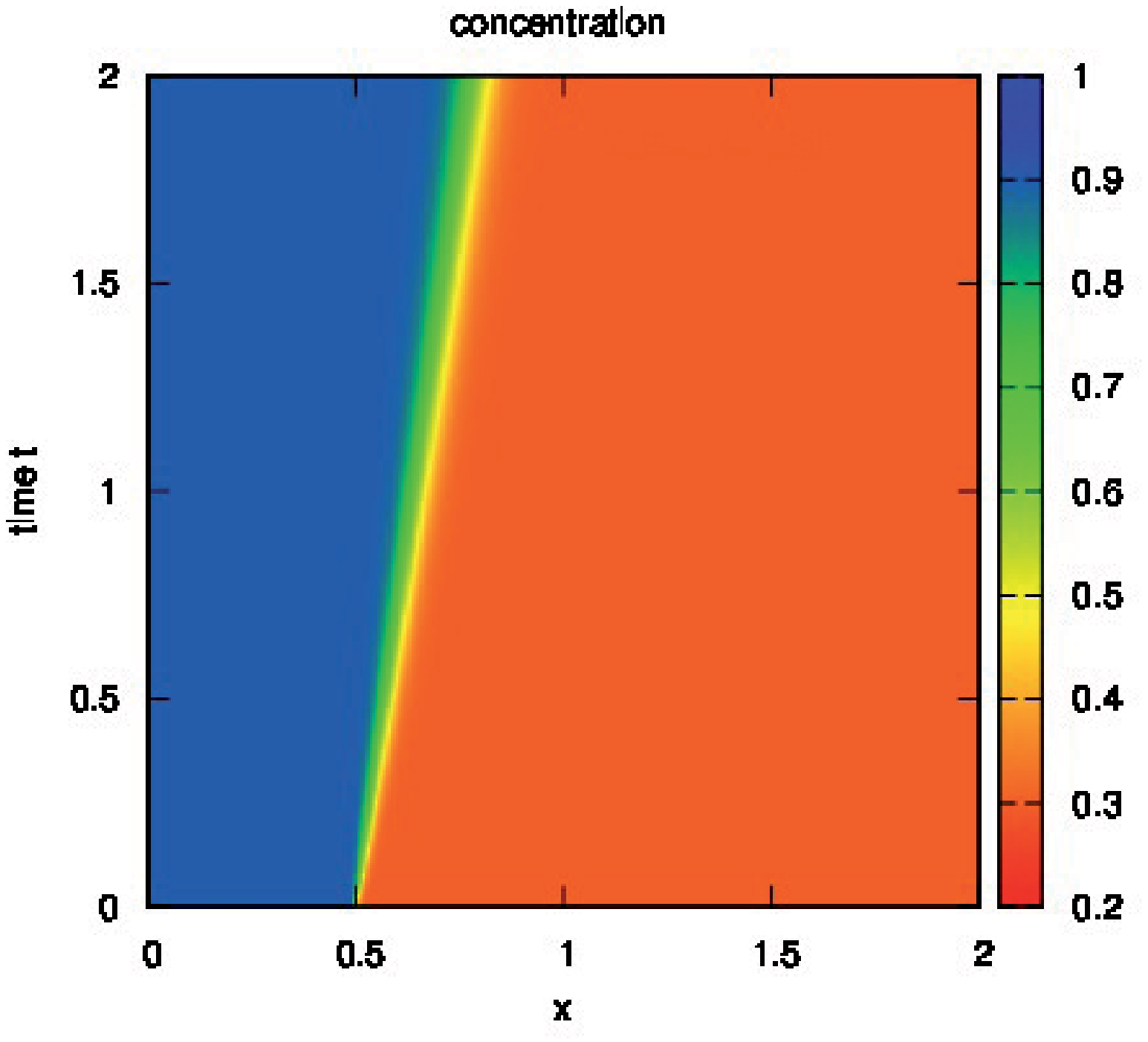} 
\caption{$s$ (left) and $c$ (right) 2D plot for data (\ref{f2}), (\ref{ivpb})and (\ref{bvivpb}).}
\label{fig-ivp2D}
\end{figure}

\begin{figure}[H]
 \includegraphics[width=6.5cm,height=6.5cm,angle=-90]{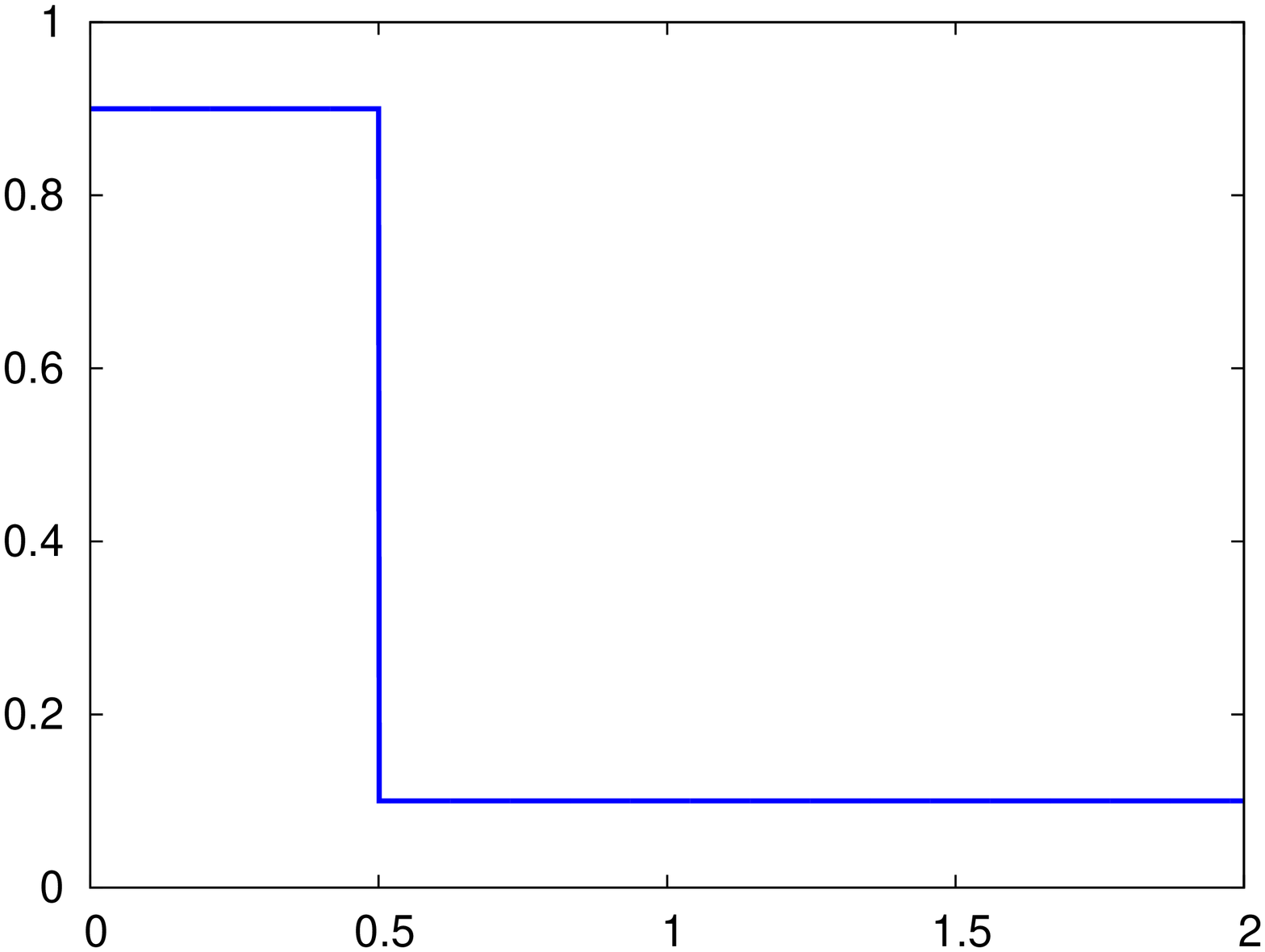}
 \includegraphics[width=6.5cm,height=6.5cm,angle=-90]{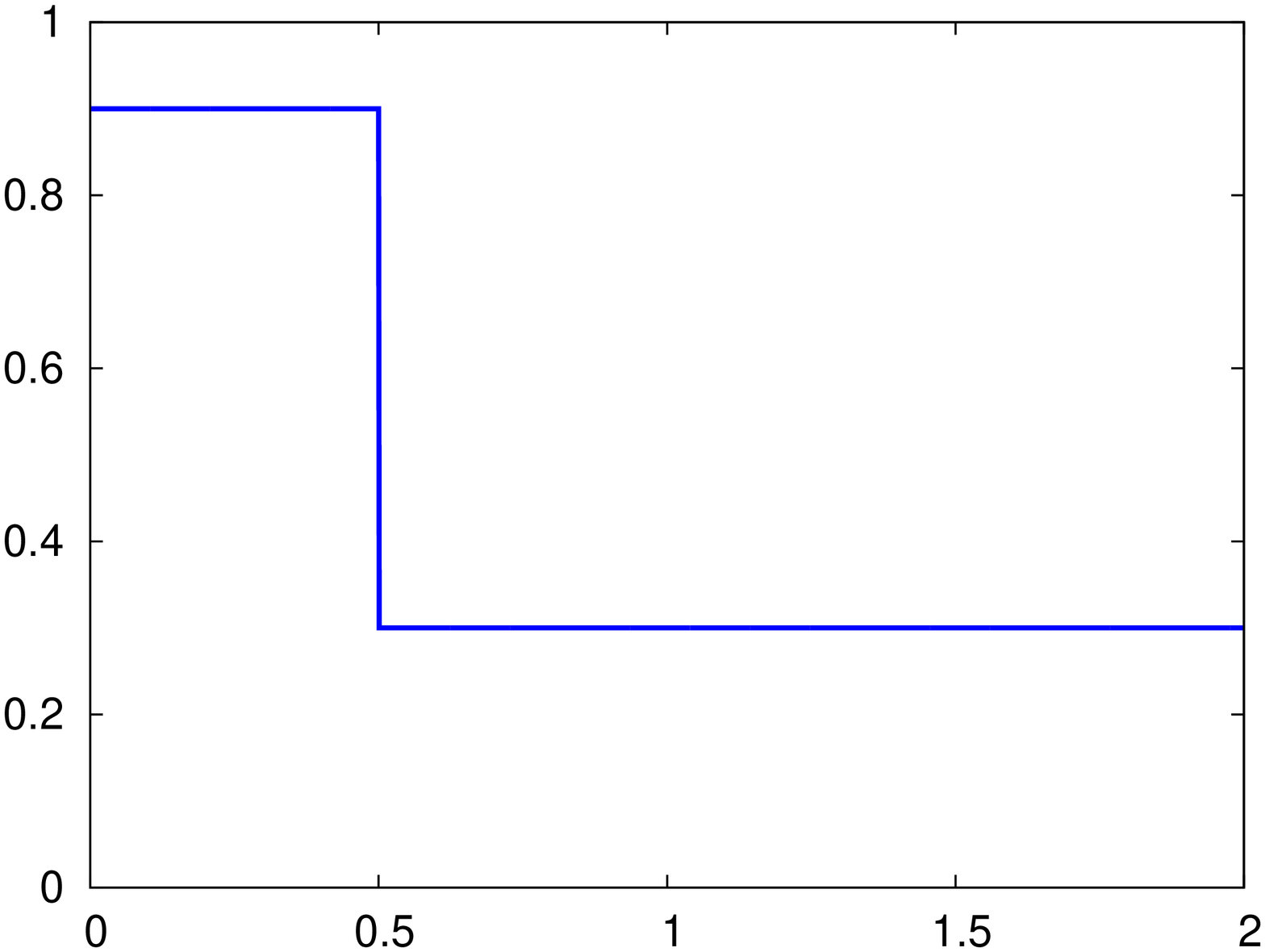} \\
 \includegraphics[width=6.5cm,height=6.5cm,angle=-90]{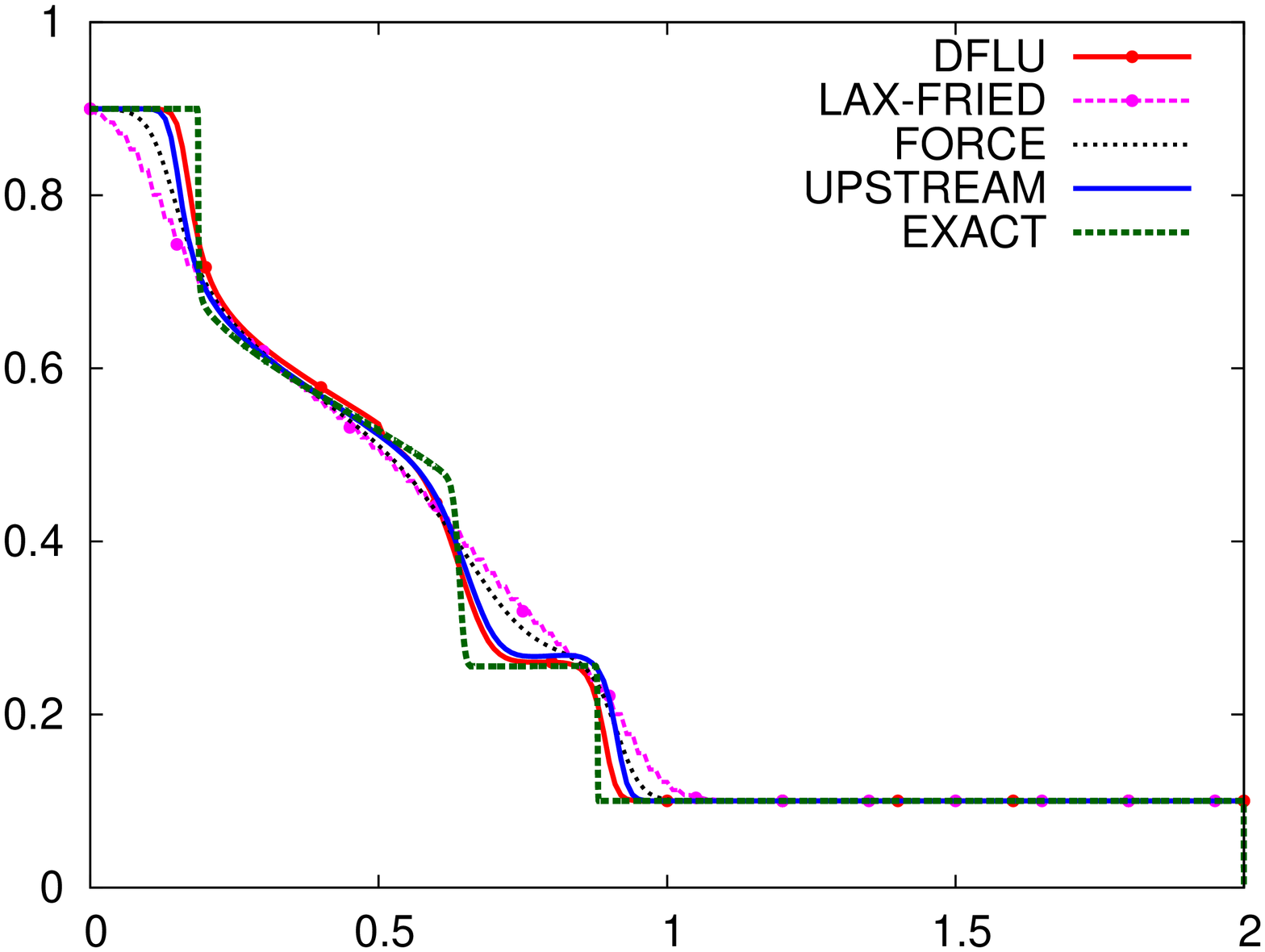}
 \includegraphics[width=6.5cm,height=6.5cm,angle=-90]{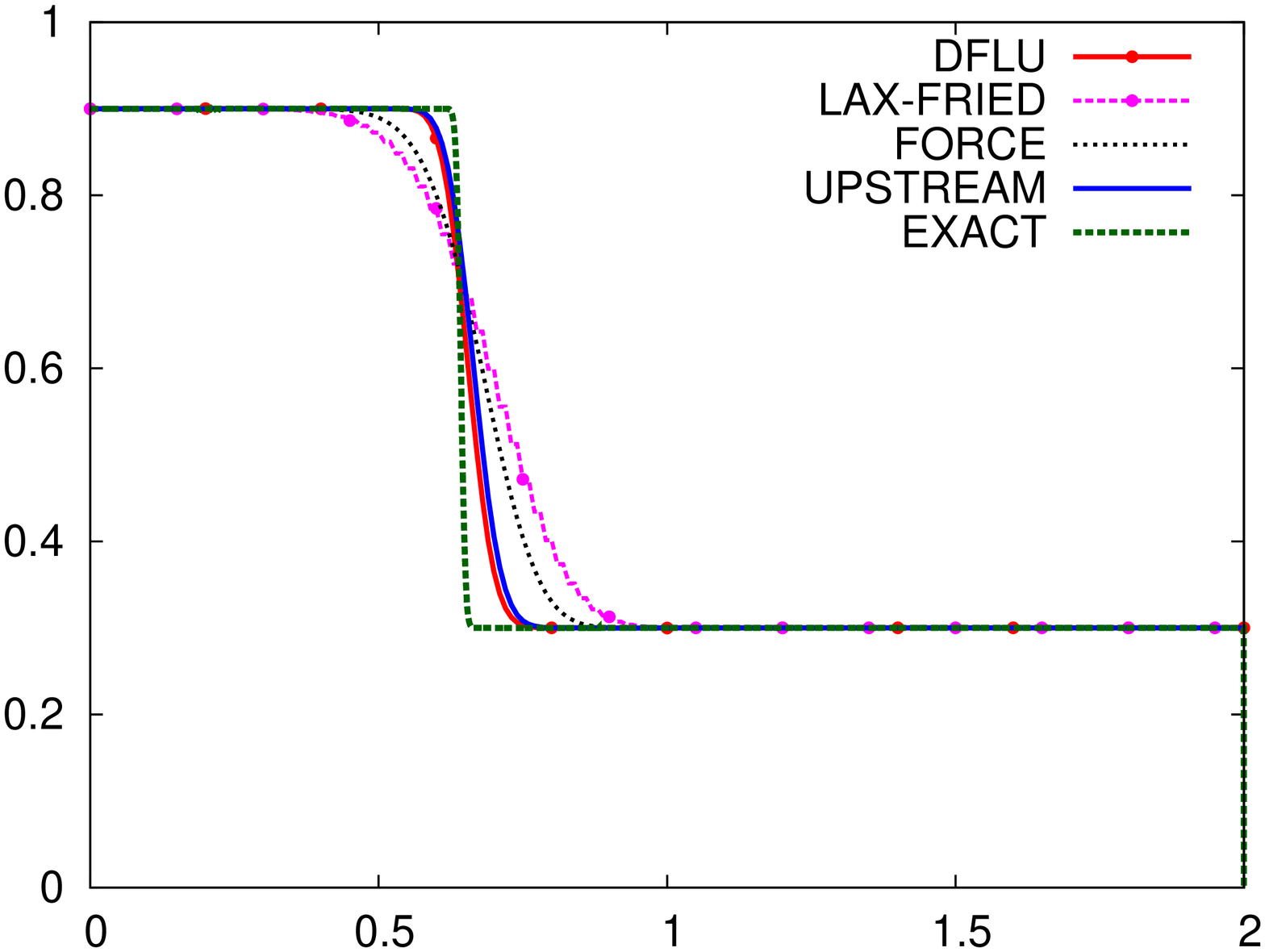} \\
 \includegraphics[width=6.5cm,height=6.5cm,angle=-90]{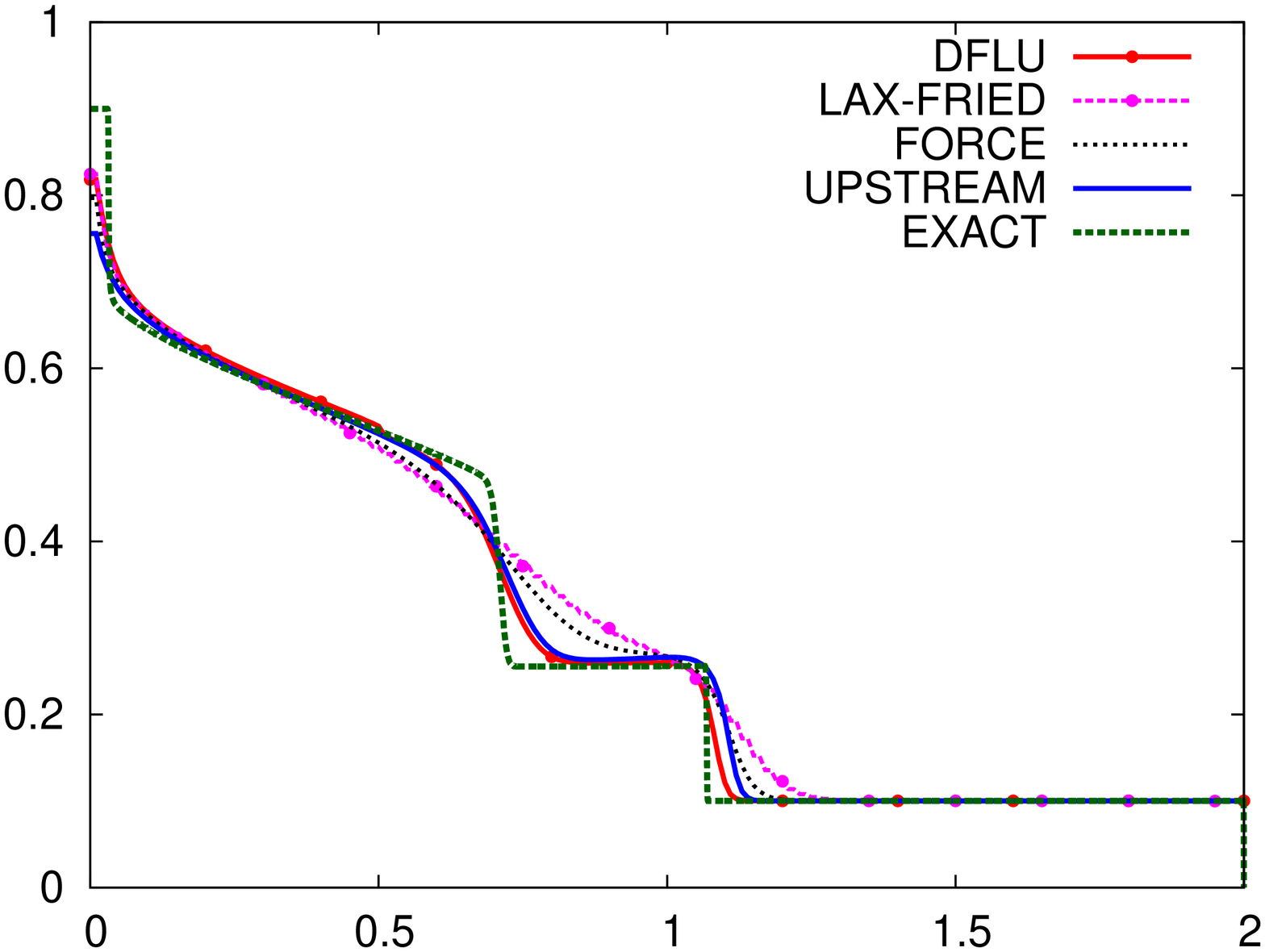}
 \includegraphics[width=6.5cm,height=6.5cm,angle=-90]{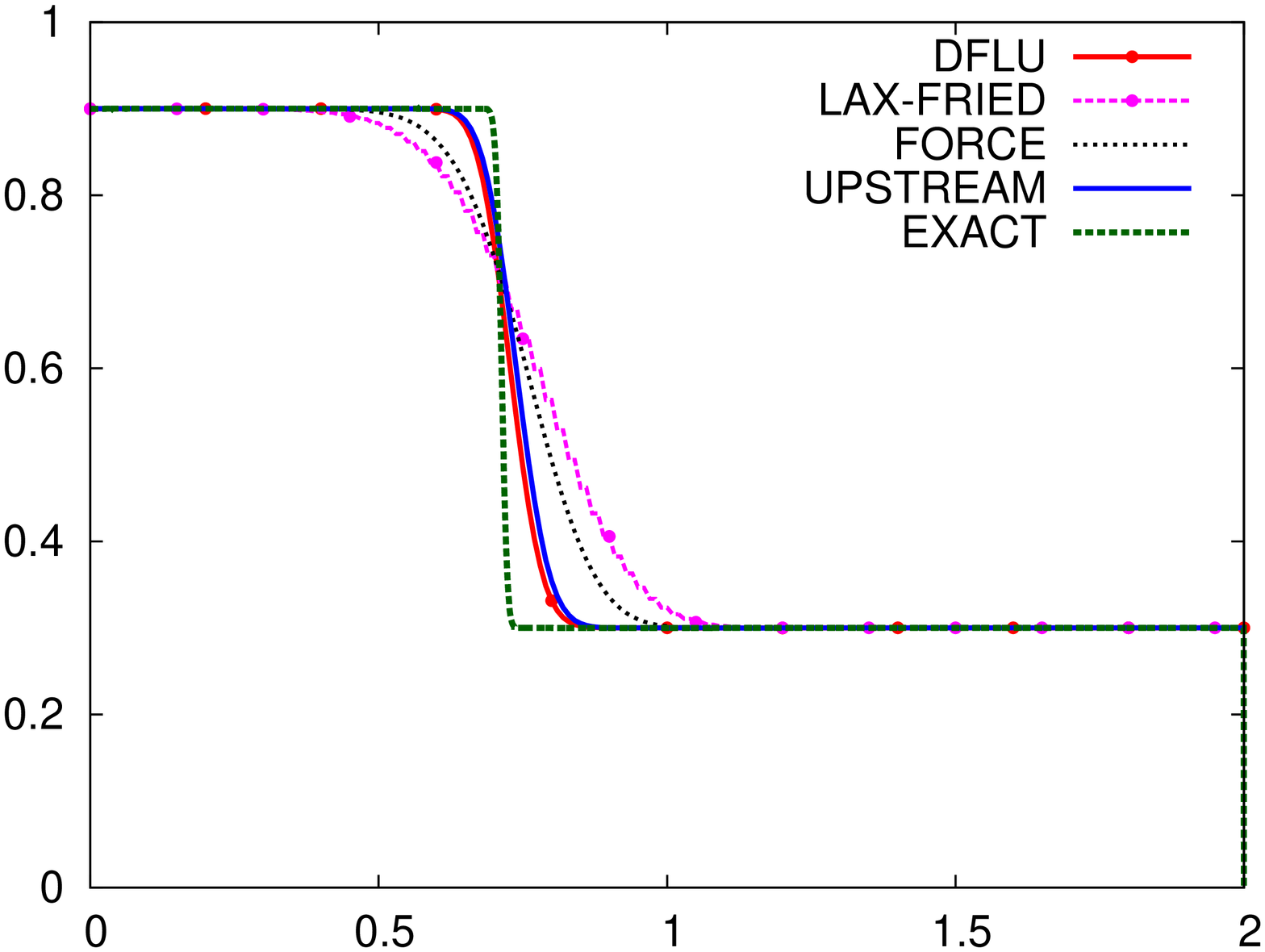}
\caption{$s$ (left) and $c$ (right) calculated at t=0., t=1. and t=1.5 for data (\ref{f2}), (\ref{ivpb}) and  (\ref{bvivpb}).}
\label{fig-154ivp}
\end{figure}

We first consider a pure initial value problem.
Initial condition (see top of Fig. \ref{fig-154ivp}) is given by   
\begin{equation} s(x,0) = \left\{ \begin{array}{lll}
.9 &\mbox{if}& x<.5, \\ .1 &\mbox{if}& x>.5 \end{array} \right.  , \quad
c(x,0) = \left\{ \begin{array}{lll}
.9 &\mbox{if}& x<.5, \\ .3 &\mbox{if}& x>.5 \end{array} \right.  .
\label{ivpb}
\end{equation}
With this initial condition we have $F^{DFLU}(s_L,s_R,c_L,c_R)=F^{G}(s_L,s_R,c_R,c_L)$ with
$ s_L=.9,s_R=.1,c_L=1.$ and $c_R=.3$. Boundary data are such that\\
\begin{equation} s(0,t)  = .9, \; s(2,t) = .1, \quad c(0,t)  = .9, \; c(2,t) = .3 \,\,\,\,\forall\,\,t \geq 0.
\label{bvivpb} 
\end{equation}
In Fig.\ref{fig-ivp2D}, a two dimensional plot in space and time for saturation and consentration is presented for the DFLU flux and in Fig. \ref{fig-154ivp} comparison
of the DFLU with other fluxes are given at time levels $t=1$ and $t=1.5$. They show that, as expected, the DFLU flux, which is the closest to a Godunov scheme, performs better than the other schemes. The upstream mobility flux, which is an upwind scheme, performs better than the  two central difference schemes, the FORCE and Lax-Friedrichs schemes. Here, in Fig.\ref{fig-bvp} and in Fig.\ref{figbvpc0u} reference(exact)  solution is calculated from DFLU with finer meshes for the comparison of various schemes

\newpage

To confirm these first observations we consider now a boundary value problem. We just changed the boundary functions, so instead of boundary conditions (\ref{ivpb}) we consider now a problem with closed boundaries, that is fluxes are zero at the boundary:
\be f \equiv 0 \,\mbox{ at } \,\,x=0 \mbox{ and }  x=2 \,\,\mbox{ for all } \, t \geq 0  .
\label{bvpbf0} \ee

They show  that, as expected, the DFLU scheme, which is the closest to a Godunov scheme, performs better than the  upstream mobility, the FORCE or the Lax-Friedrichs schemes.

The purpose of the last experiment whose results are shown in Fig. \ref{figbvpc0u} is to show the effect of polymer flooding. In this experiment we remove polymer flooding  and take  $c\equiv 0$ at all time. By comparing with the solution shown in Fig. \ref{fig-bvp} bottom left we observe that as expected the saturation front is moving faster since there is no retardation due to the increase of viscosity of the wetting fluid caused by the polymer injection. We also observe that the structure of the solution is less complex. In the absence of concentration FORCE scheme is closer to Upstream Mobility in that it has less diffusion, compare figures Fig.\ref{figbvpc0u} and Fig.\ref{fig-bvp}. In the presence of concentration,
it is diffusive particularly more at the points where the concentration $c$ is discontinuous. 

\begin{figure}[hb]
\includegraphics[width=7cm]{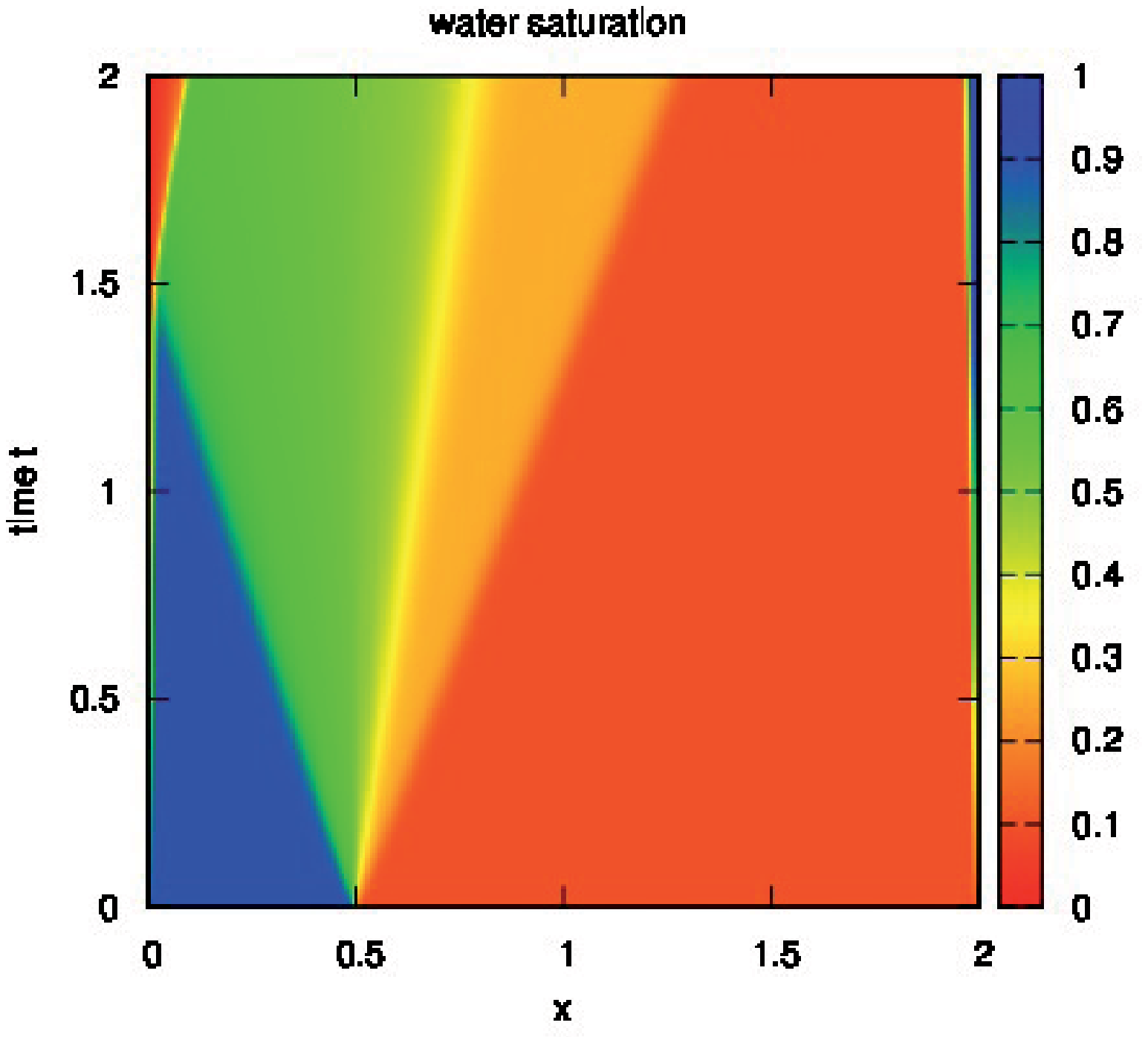} \hspace*{0.5cm}
\includegraphics[width=7cm]{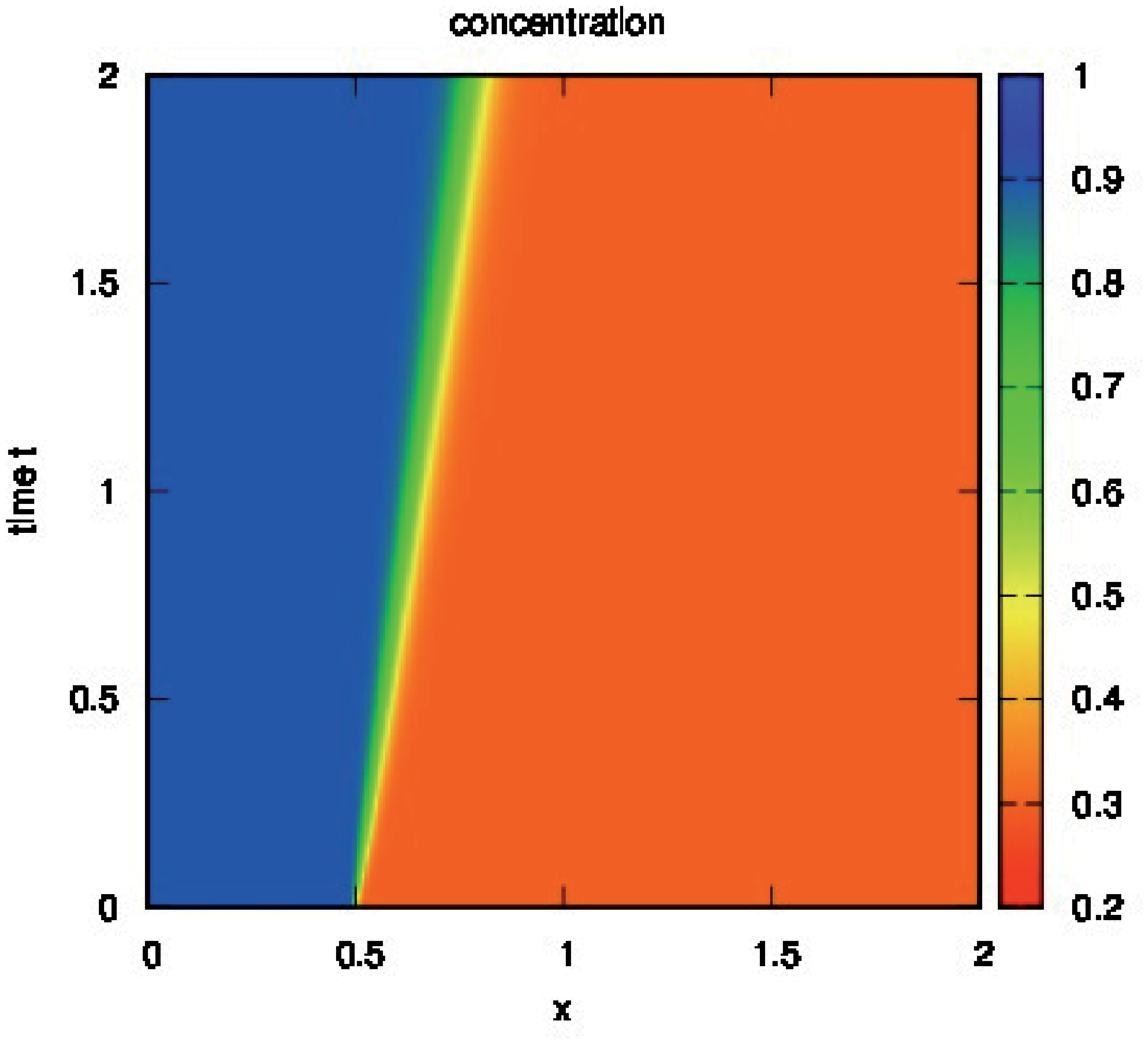} 
\caption{$s$ (left) and $c$ (right) 2D plot for data (\ref{f2}), (\ref{ivpb}) and  (\ref{bvpbf0}).}
\label{fig-bvp2D}
\end{figure}

\newpage

\begin{figure}[H]
\includegraphics[width=6.5cm,height=6.5cm,angle=-90]{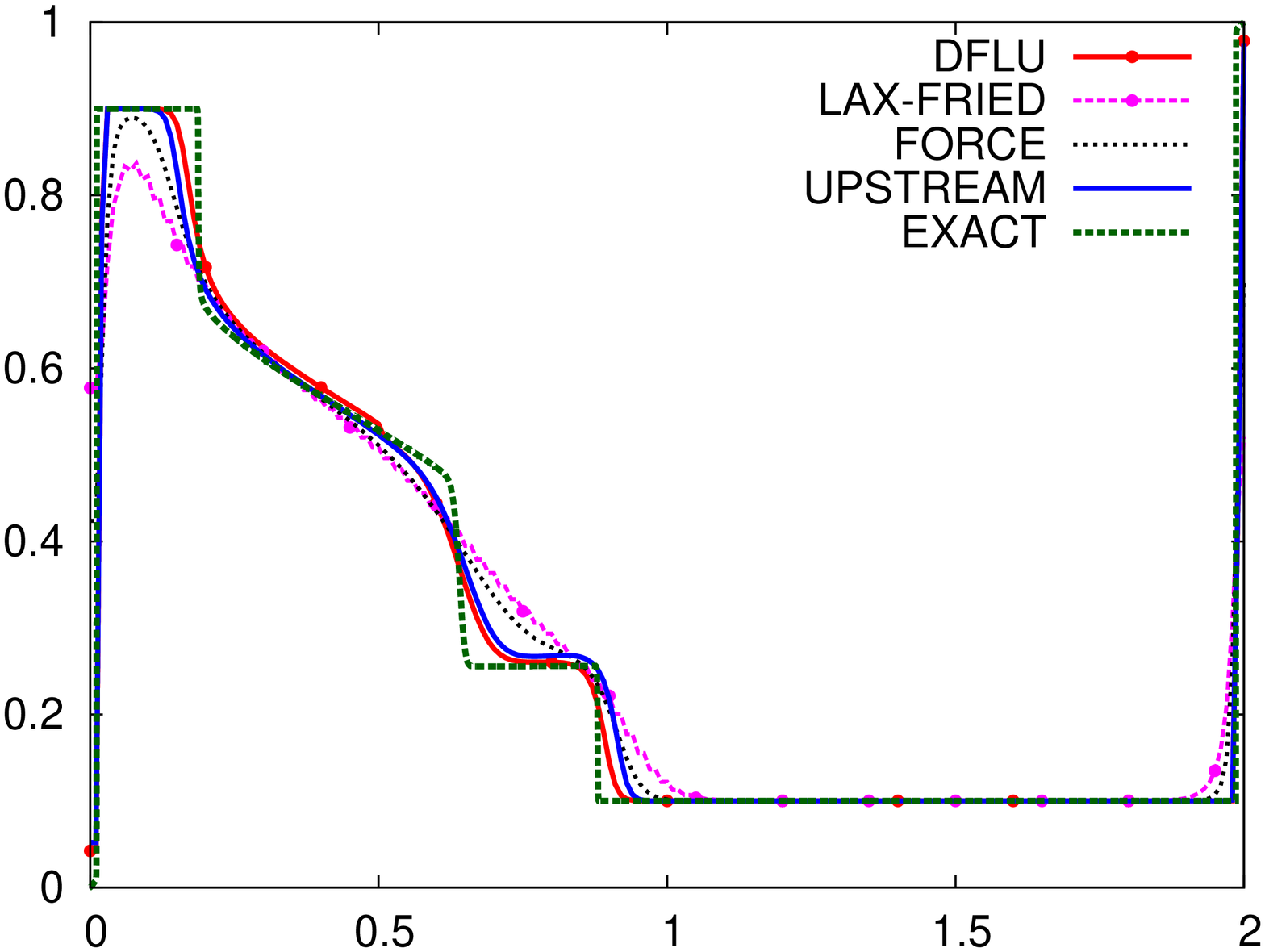} \hspace*{0.5cm}
\includegraphics[width=6.5cm,height=6.5cm,angle=-90]{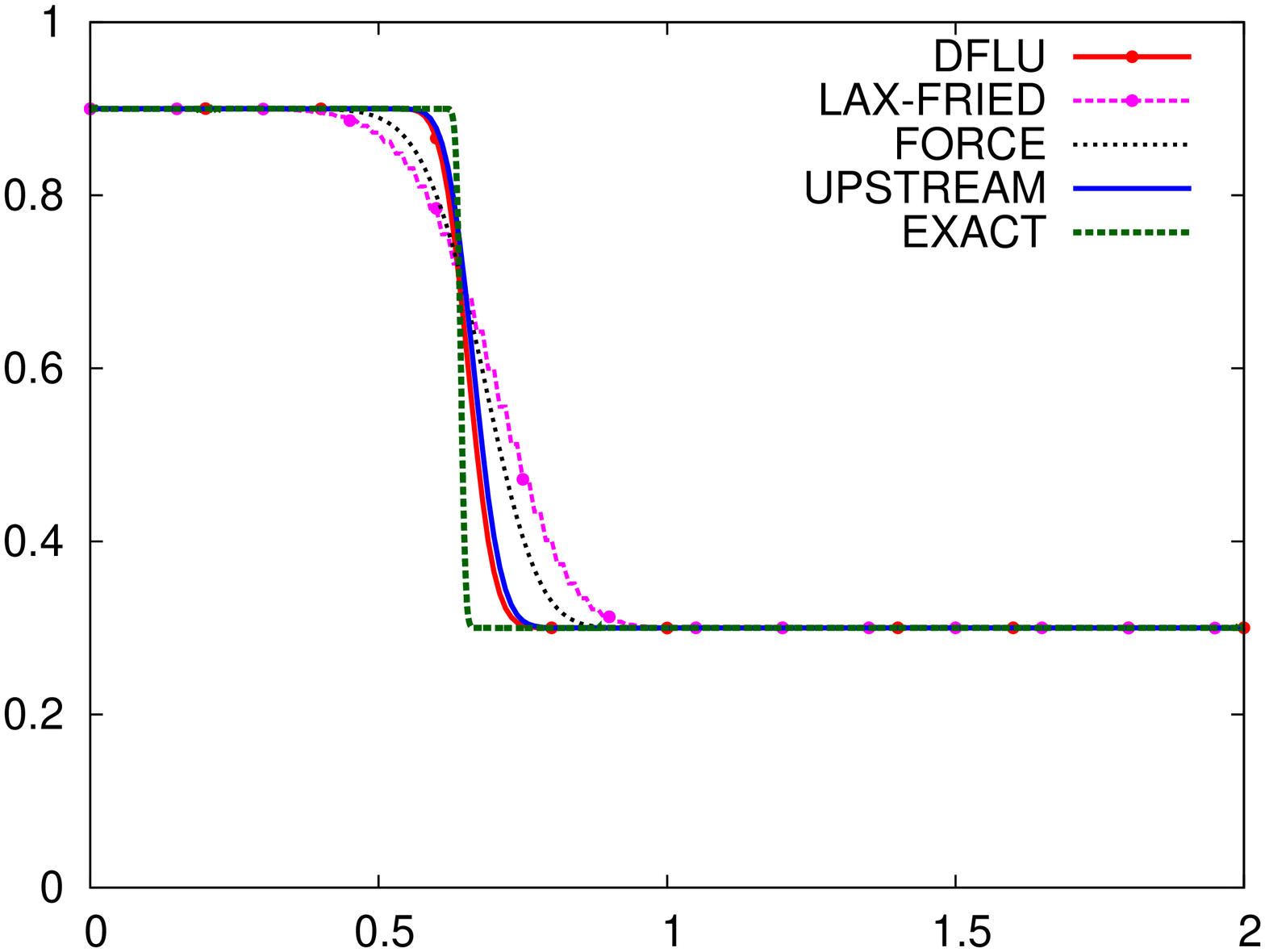}\\
\includegraphics[width=6.5cm,height=6.5cm,angle=-90]{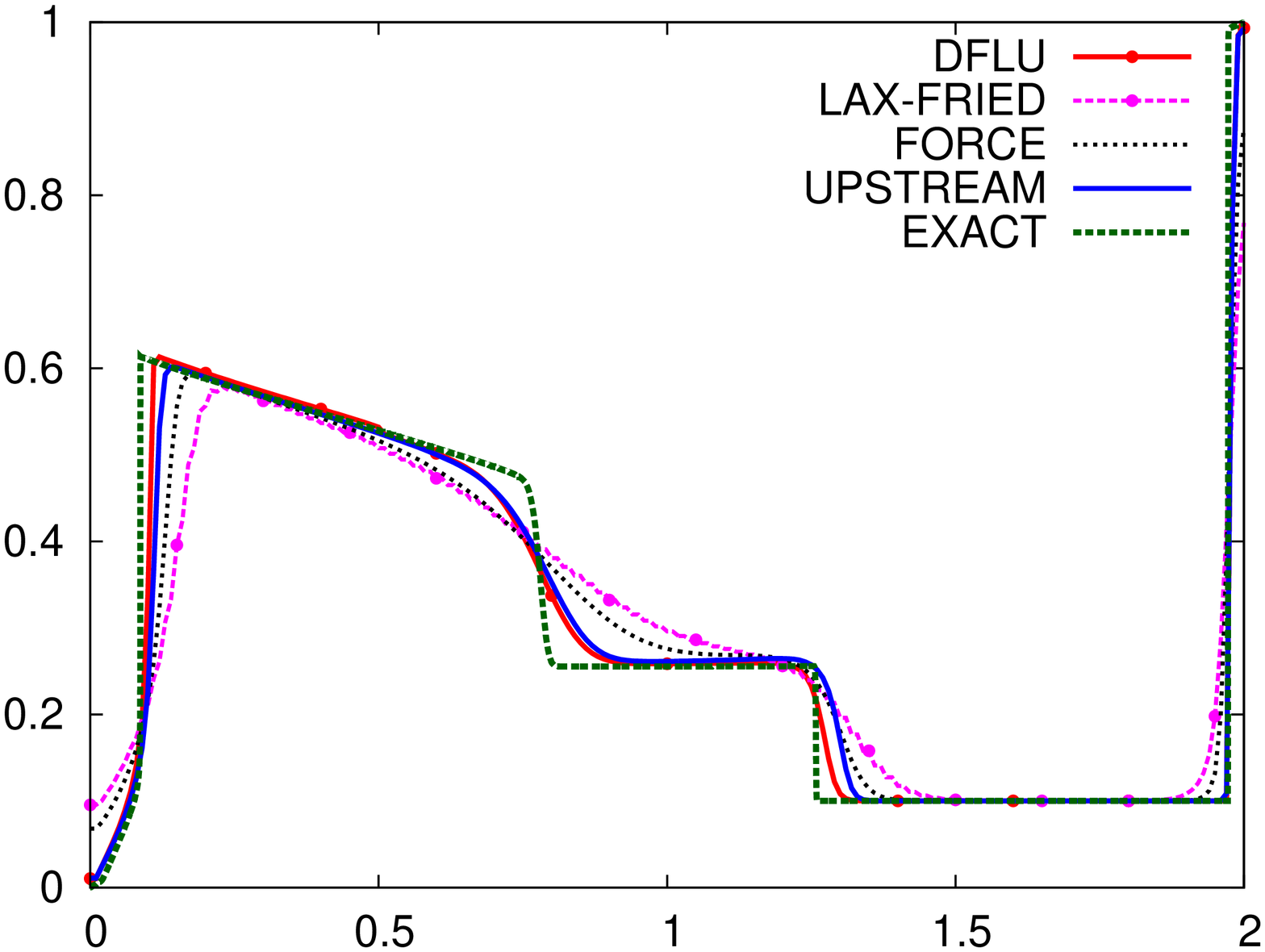} \hspace*{0.5cm}
\includegraphics[width=6.5cm,height=6.5cm,angle=-90]{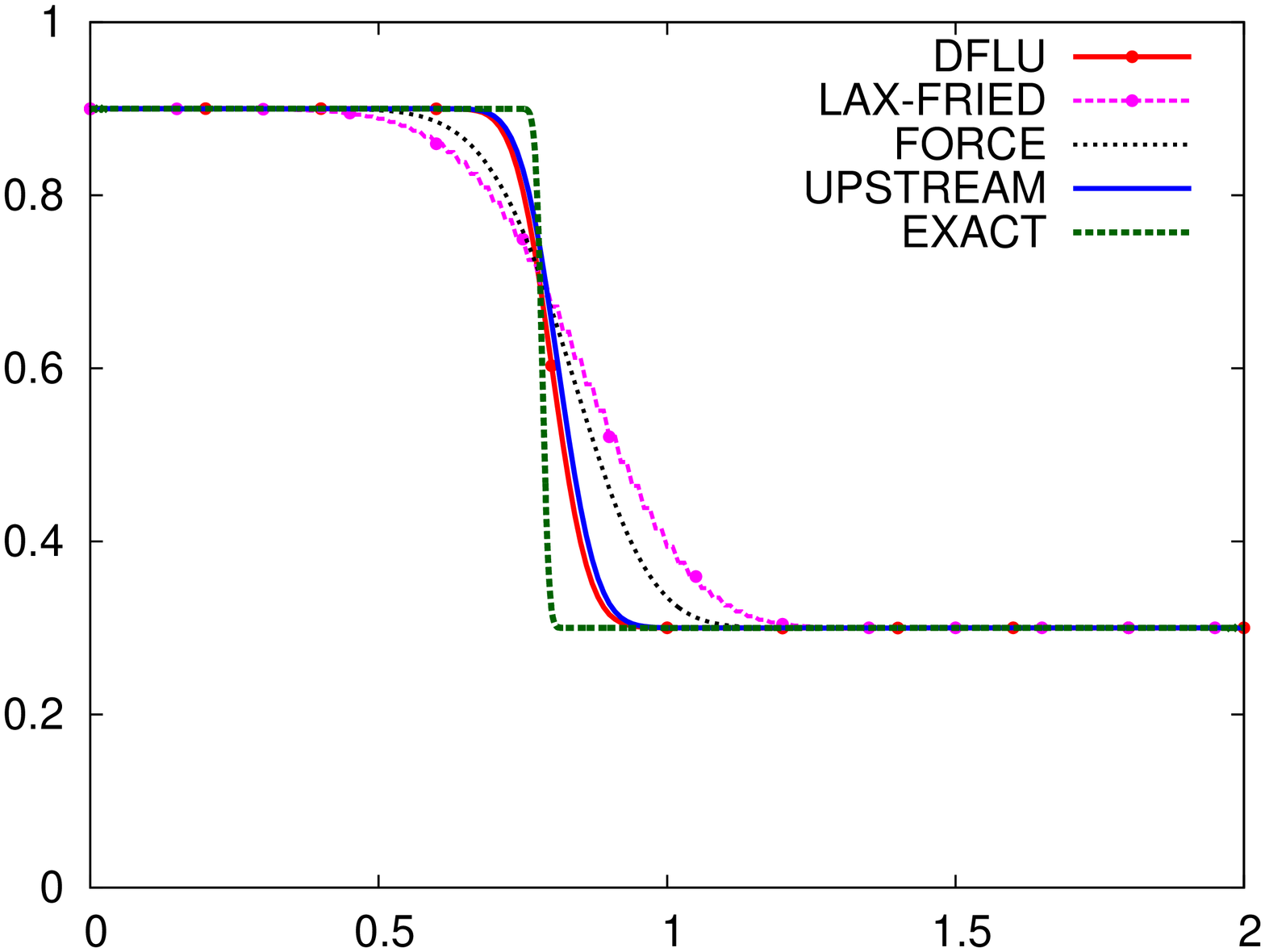}\\
\includegraphics[width=6.5cm,height=6.5cm,angle=-90]{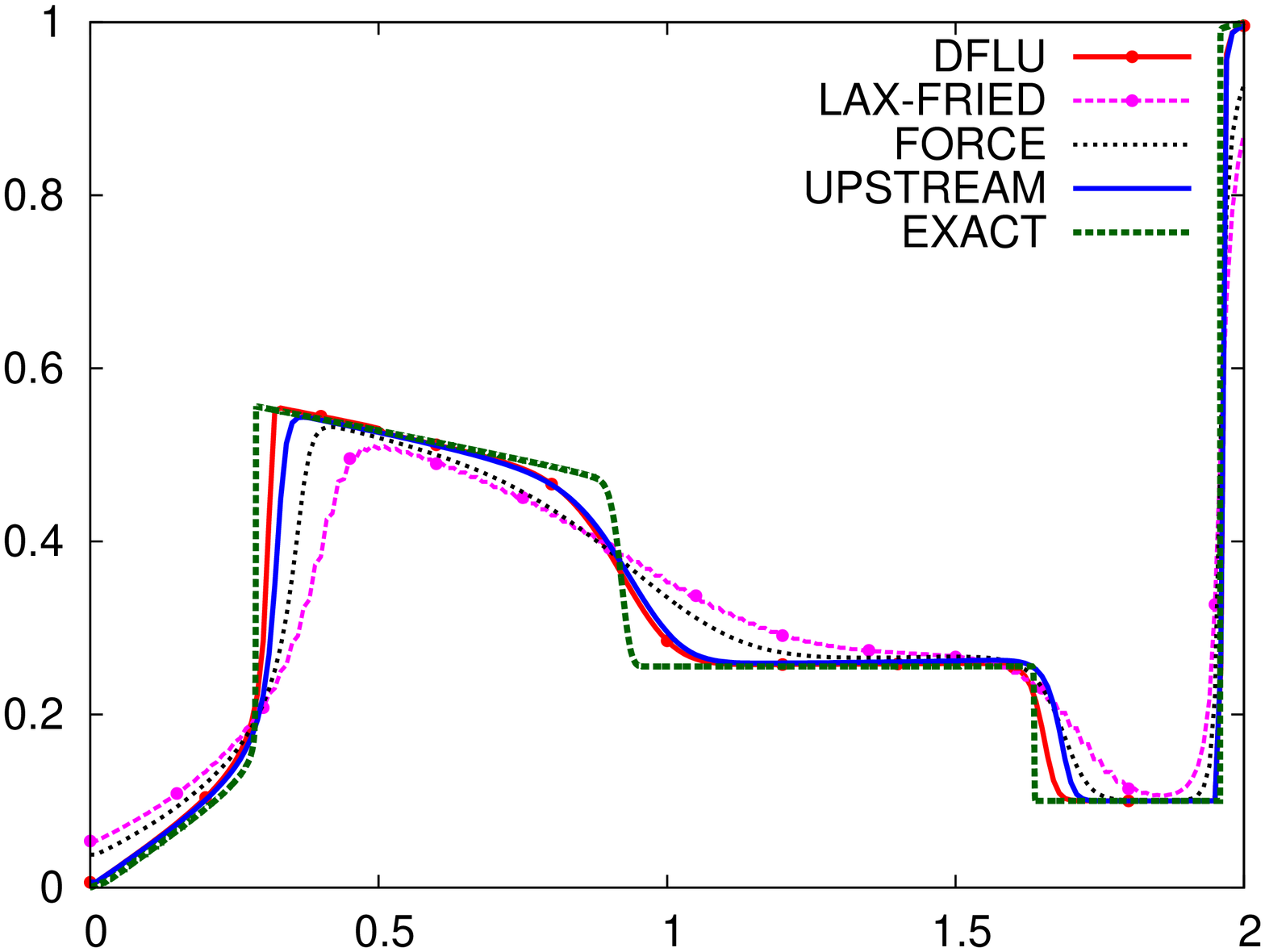} \hspace*{0.5cm}
\includegraphics[width=6.5cm,height=6.5cm,angle=-90]{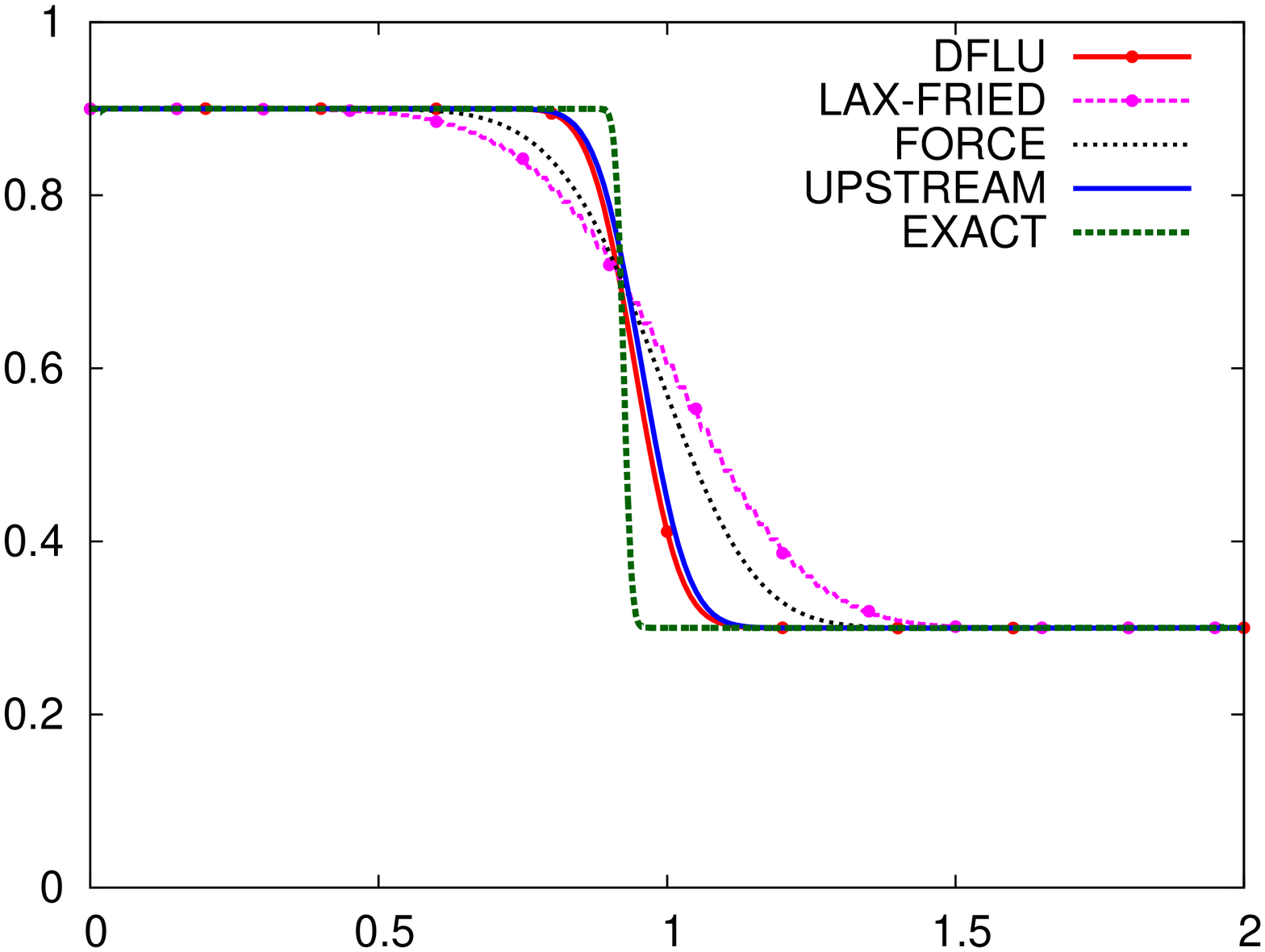}
\caption{$s$ (left) and $c$ (right) calculated at t=1., t=2. and t=3. for data (\ref{f2}), (\ref{ivpb}) and (\ref{bvpbf0}).}
\label{fig-bvp}
\end{figure}

\begin{figure}[h]
\includegraphics[width=6.5cm,height=6.5cm,angle=-90]{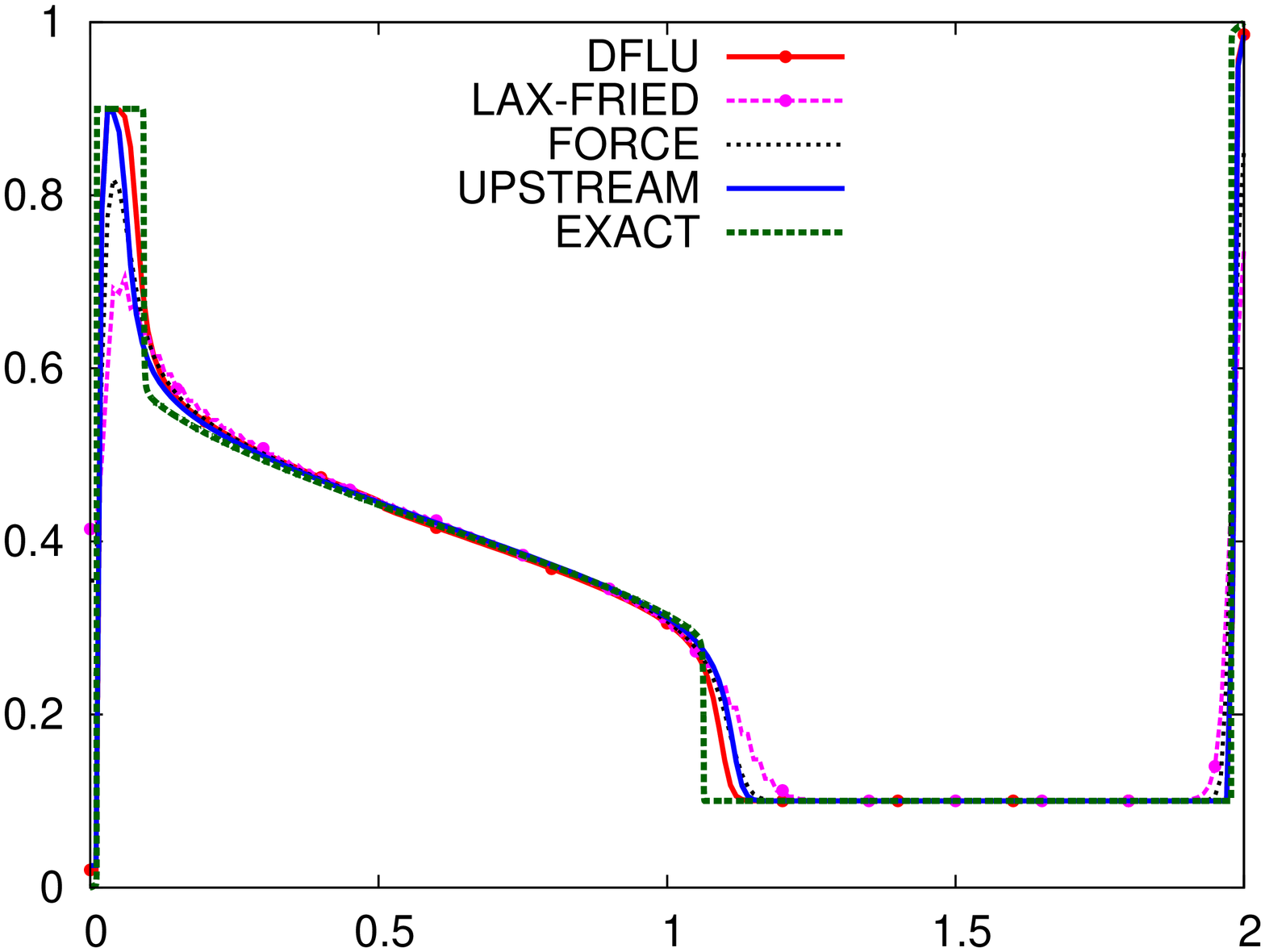} \hspace*{0.5cm}
\includegraphics[width=6.5cm,height=6.5cm,angle=-90]{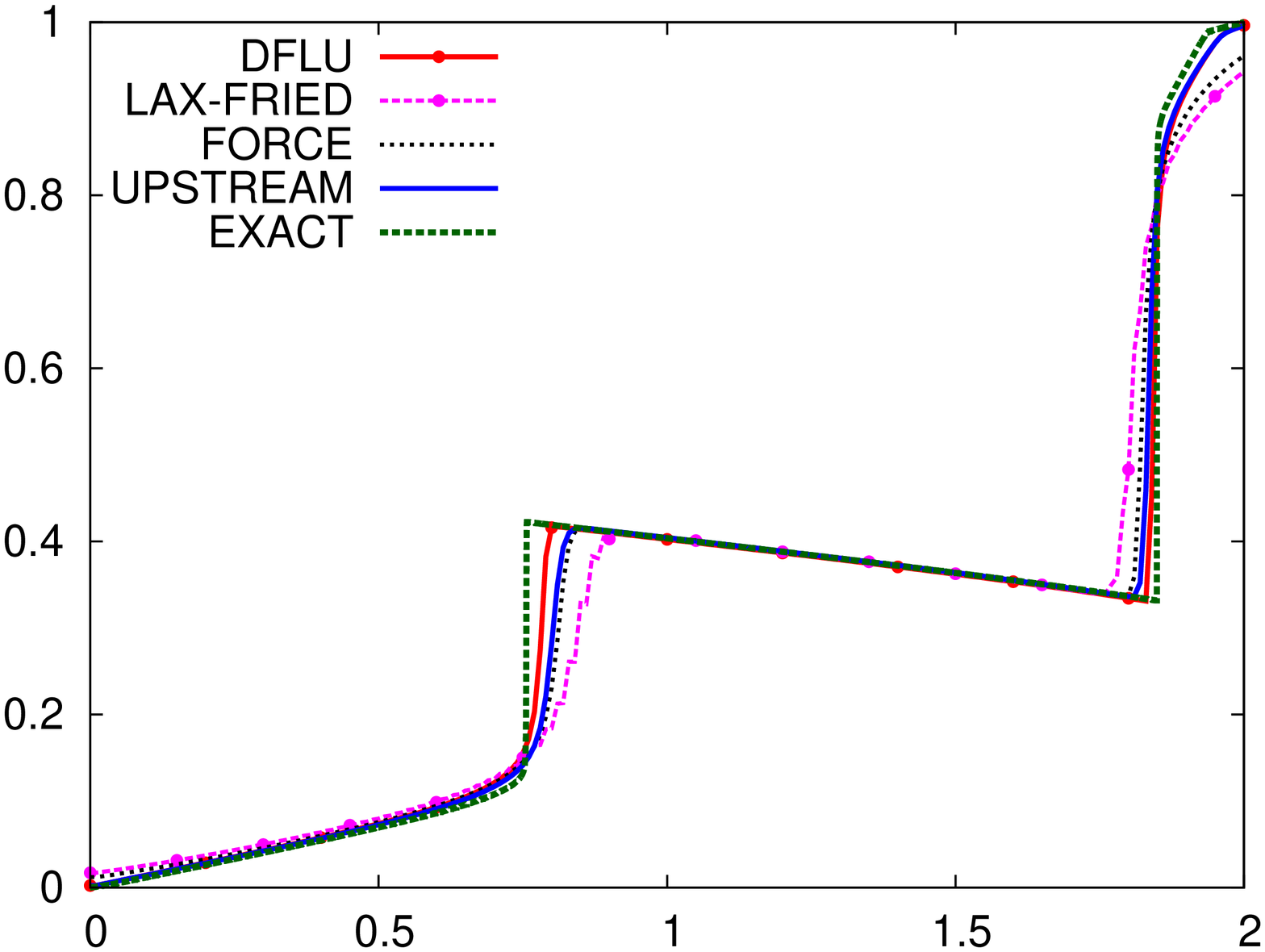}
\caption{$s$  calculated at t=1. and t=3. for same data as in Fig. \ref{fig-bvp} but without polymer injection.}
\label{figbvpc0u}
\end{figure}

\section{ \bf Polymer flood model with flux function discontinuous in the space variable}
\label{dispolymer}
     In this section, we extend the previous results to the case where the polymer flooding  model has a flux function discontinuous in the space variable: 

\be
\begin{array}{rrll}
 s_t + f(s, c,x)_x & =& 0 \\
 (sc+a(c))_t+(cf(s,c,x))_x &=&0
 \label{polydiseq1}
\end{array}
\ee
where $x \rightarrow f(s,c,x)$ is  discontinuous. For simplicity we assume that
$f$ has a single discontinuity at $x=0$.i.e.,
$$f(s,c,x)=H(x)f_l(s,c)+(1-H(x))f_r(s,c)$$
where $H$ is a Heaviside function and $f_l$ and $f_r$ as in section 1, satisfies the following conditions, for $p=l,r$
\begin{enumerate}
\item[(i)] $f_p(s,c) \ge 0, f_p(0,c)= f_p(1,c)=0$ for all $c \in I$.
\item[(ii)] The function $s \rightarrow f_p(s,c)$ has exactly one global
 maximum in $I$ with $\theta_p=$argmax$(f_p)$.
\item[(iii)] $\frac{\partial f_p}{\partial c}(s,c) < 0  \,\,\forall \,\,s \in (0,1)$ and for all $ c \in I$
 \end{enumerate}

  Equations of type (\ref{polydiseq1}) arise while dealing with polymer 
 flooding of oil reservoirs which are heterogeneous \cite{DarGlimMcbr88}.\\

\noi {\bf Remark}: Since $f$ is discontinuous at $x=0$, then the Rankine-Hugoniot condition for system (\ref{polydiseq1}) gives
$$
\begin{array}{rrll}
f_l(s^-,c^-)&=&f_r(s^+,c^+)\\
c^-f_l(s^-,c^-)&=&c^+f_r(s^+,c^+)
\end{array}
$$
where $(s^-,c^-)$ and $(s^+,c^+)$ denotes the left and right values of $(s,c)$ across the line $x=0$. This implies 
\begin{equation} c^-=c^+  \label{ccont} \end{equation} so
$c$ cannot have a discontinuity across the line $x=0$.\\

\begin{figure}
\begin{center}
 \begin{picture}(400,200)(0,-10)
{\includegraphics[height=6cm,width=12cm]{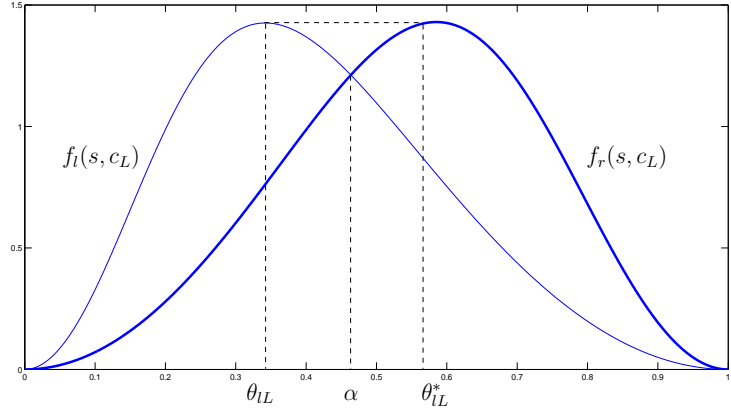}}
\end{picture}
\end{center}
\caption{ Flux functions $f_l(s,c_L),$ and $f_r(s,c_L)$ .}
\label{disfig}
 \end{figure}

    The solution to the Riemann problem corresponding to (\ref{polydiseq1}) is given in the Appendix. We now  present a numerical experiment to compare the DFLU, the upstream mobility, the FORCE and the Lax-Friedrichs fluxes in the case where the flux function $f$ is discontinuous in space:
$$f(s,c,x)=H(x)f_l(s,c)+(1-H(x))f_r(s,c)$$
where $H$ is the Heaviside function and $f_l$ and $f_r$ are given by
\be \begin{array}{l} 
f_r(s,c) = \dfrac{\lambda_1(s,c)}{\lambda_1(s,c) + \lambda_2(s,c)} [ \var + (g_1-g_2)\lambda_2(s,c) ],\\
f_l(s,c)= \dfrac{\mu_1(s,c)}{\mu_1(s,c) + \mu_2(s,c)} [ \var + (g_1-g_2)\mu_2(s,c) ] \\
a(c)=.25c
\end{array} \label{df2} \ee
where
$$ \begin{array}{l}
\lambda_1(s,c)=\dfrac{10 s^2}{.5+c}, \lambda_2(s,c)=20 (1-s)^2,\\
\mu_1(s,c)=\dfrac{50 s^2}{.5+c}, \mu_2(s,c)=5 (1-s)^2, \,\,g_1=2, g_2=1\,\,\, \mbox{and}\,\,\, \var=0
\end{array}  $$
(see Fig. \ref{disfig}),
with the initial condition
$$ s(x,0) = \left\{ \begin{array}{lll}
.9 &\mbox{if}& x<0, \\ .1 &\mbox{if}& x>0 \end{array} \right.  , \quad
c(x,0) = \left\{ \begin{array}{lll}
.9 &\mbox{if}& x<0, \\ .3 &\mbox{if}& x>0 \end{array} \right.  .
$$

Following   \cite{AdiJafGow04} the DFLU flux at the interface is given by
\be
\begin{array}{lcll}
 \bar{F}(s_{-1}^n,c_{-1}^n,s_1^n,c_1^n) 
             & = & \min\{ f_l(\min\{s_{-1}^n,\theta_{-1}^n\},c_{-1}^n),f_r(\max\{s_{1}^n,\theta_1^n\},c^n_1)\},
\end{array}
\ee
where
$ \theta_{-1}^n =\mbox{ argmax } f_l(\cdot,c_{-1}^n) $ and
$ \theta_1^n =\mbox{ argmax } f_r(\cdot,c_1^n) $.

\begin{figure}
\includegraphics[width=5.2cm, angle=-90]{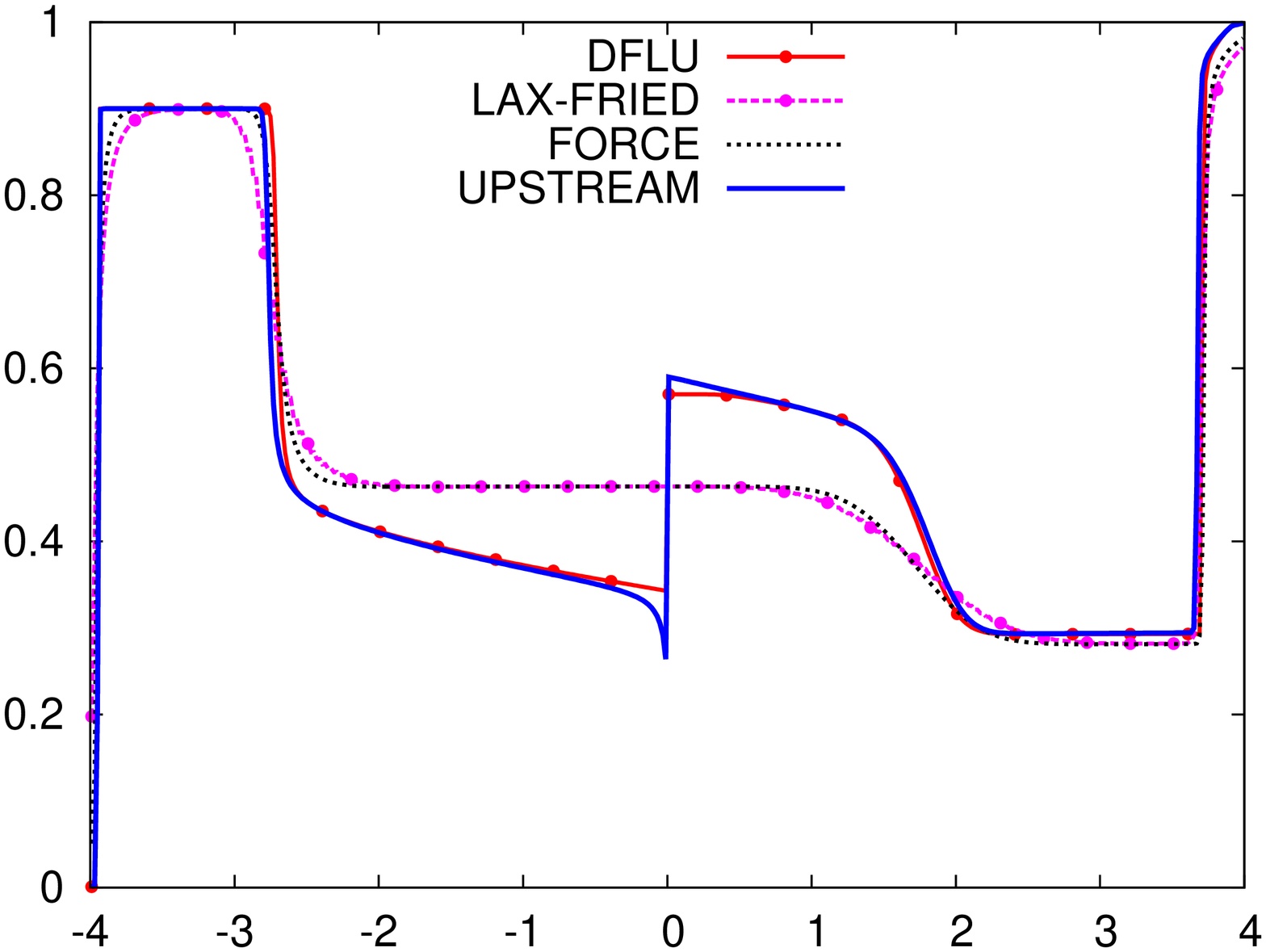}
\includegraphics[width=5.2cm, angle=-90 ]{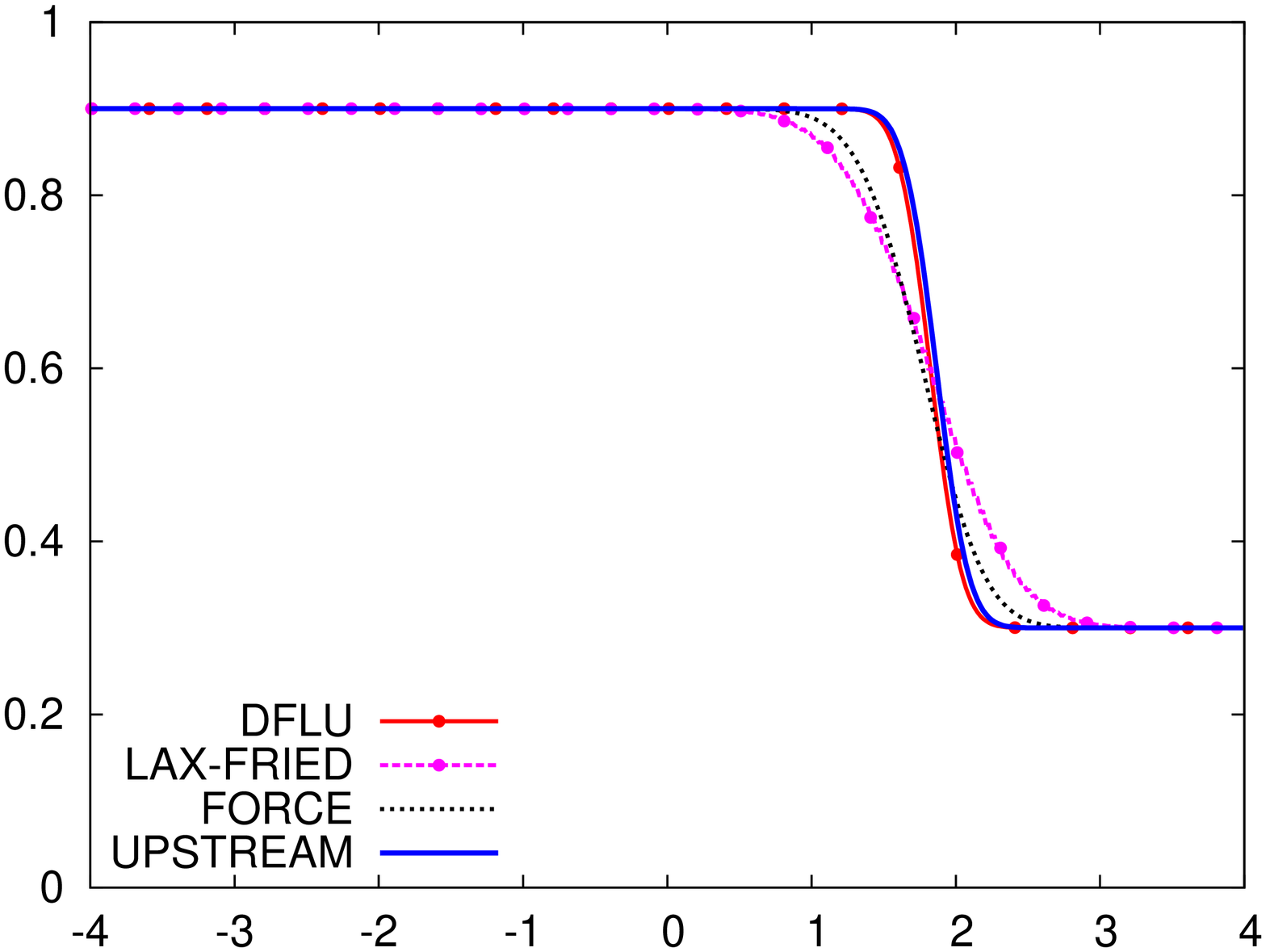}\\
\includegraphics[width=5.2cm, angle=-90]{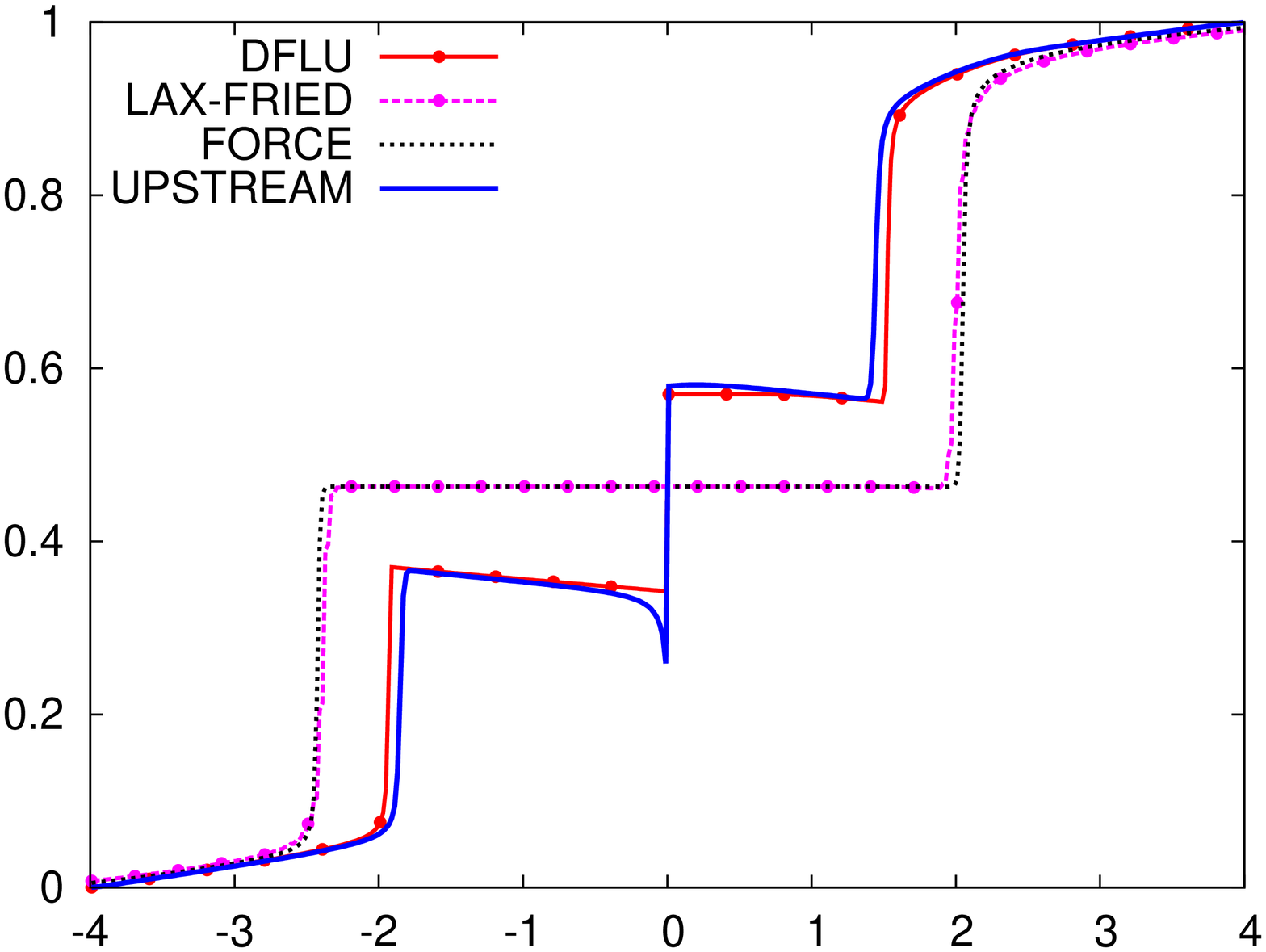}
\includegraphics[width=5.2cm, angle=-90]{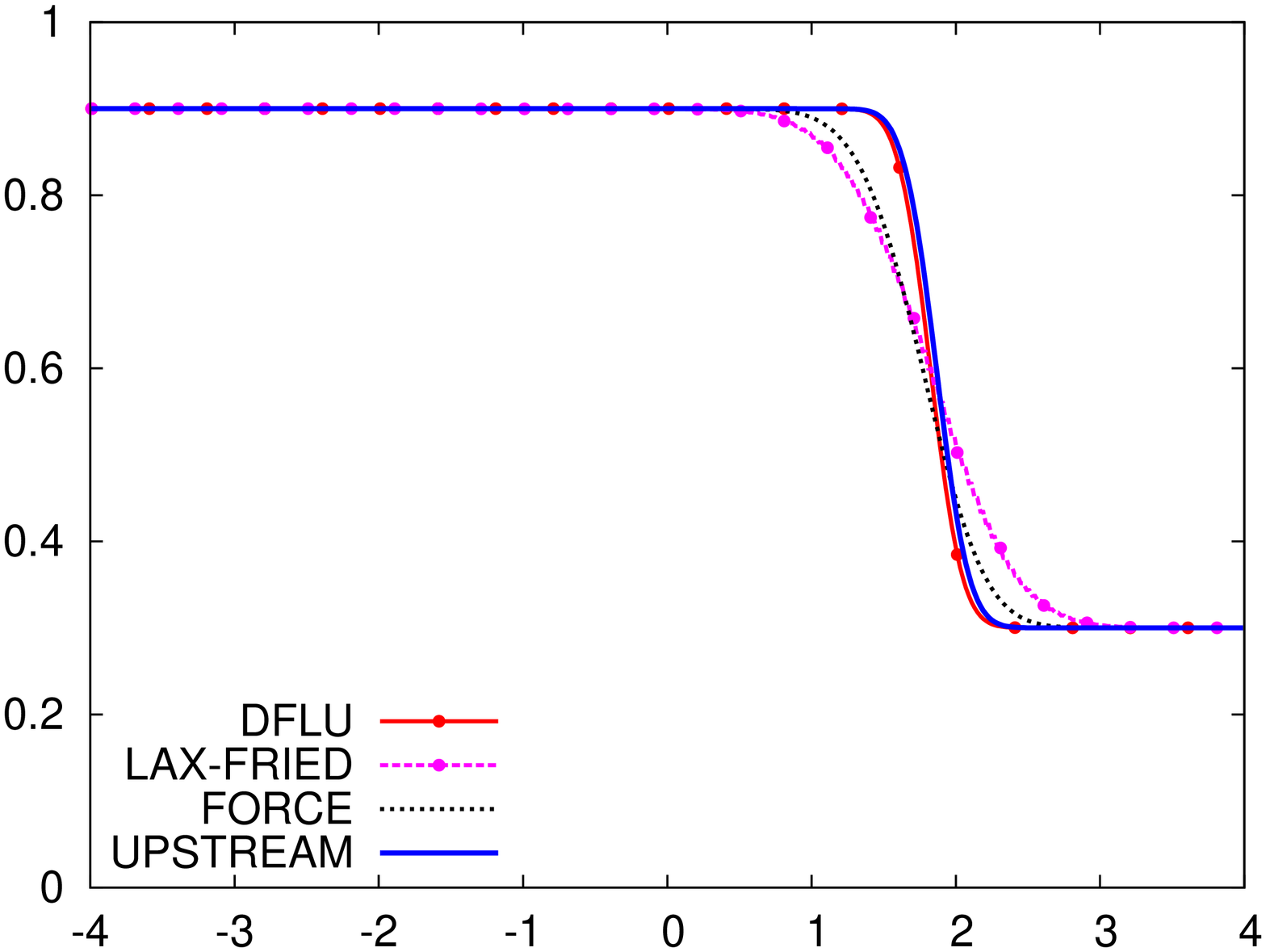}
 \caption{$s$ (left) and $c$ (right) at t=1 and t=2}
\label{discomparison}
\end{figure}

     Here we considered the case where the flux functions $ f_l(s,c_L)$ and $f_r(s,c_L)$ intersect at a point $\alpha$ where $ \frac{\partial f_l(s,c_L)}{\partial  s} < 0 $ and  $ \frac{\partial f_r(s,c_L)}{\partial s} > 0 $. At $\theta_{lL}$ and $\theta_{rL}$, $f_l(s,c_L)$ and  $f_r(s,c_L)$ attains their respective maxima.
Let  $\theta^*_{lL}$ be a point  such that  $f_l(\theta_{lL},c_L) = f_r(\theta^*_{lL},c_L)$. 
For the above $f_l(s,c_L)$ and $f_r(s,c_L)$, $\alpha=.464, s_L=.9,s_R=.1,c_L=.9,c_R=.3,\theta_{lL}=.342$ and $\theta_{lL}^*=.57$ (see Fig \ref{disfig}). This is an undercompressive intersection as in \cite{AdiJafGow04}. As the Lax-Friedrichs and the
 FORCE schemes are obtained from a linear parabolic regularisation, solutions  obtained from them differ from solutions obtained from  the upstream mobility  and the DFLU schemes  for an undercompressive initial data(see Fig.\ref{discomparison}). The Lax-Friedrichs and the FORCE schemes converge to the weak solution with a $(A,B)$ entropy condition \cite{AdiMisGow05} at the interface with $A=B=\alpha$ and the
 DFLU scheme and the upstream mobility flux schemes converge to the weak solution
 with a $(A,B)$ entropy condition at the interface $A=\theta_{lL},B=\theta_{lL}^*$.
In these numerical experiments here, the discretization is such that $\Delta t=1/600$ and $h=1/50$.

\section{Conclusion}
The DFLU flux defined in \cite{AdiJafGow04} for scalar conservation laws was used to construct a new scheme for a class of system of conservation laws such as systems modeling   polymer flooding in oil reservoir engineering. The resulting DFLU flux  is based on Godunov type flux for single conservation laws but with discontinuous coefficients. It is  easy to implement as it is not using detailed  information of eigenstructure of the full system. It is very close to the flux given by an exact Riemann solver and the corresponding finite volume scheme compares favorably to other schemes using the uptream mobility, the Lax-Friedrichs and the FORCE fluxes.  The extension to the case with a change of rock type is straightforward since the DFLU flux was built to solve this case. It will work even in cases where the upstream mobility fails \cite{SidJaf09}. Here we assumed, flux $f=f(s,c)$ is  not changing the sign which is equivalent to saying that second eigen value in (\ref{polyeq1}) is not allowed to change the sign. The sign changing case and the extension to system of polymer flooding in multidimensional case  will be taken up in a forth coming paper.  In a separate paper \cite{AdiJafGow09c} we show how to use the DFLU flux to solve Hamilton-Jacobi equations with a discontinuous Hamiltonian.\\

{ \bf Appendix. Riemann problem for a polymer flooding model with a discontinuous flux}:
  In this Appendix we briefly describe the construction of the solution to a Riemann problem associated to the system (\ref{polydiseq1}) with the initial condition
\be 
\tag{A-1}
s(x,0) = \left\{ \begin{array}{lll}
s_L &\mbox{if}& x<0, \\ s_R &\mbox{if}& x>0 \end{array} \right.  , \quad
c(x,0) = \left\{ \begin{array}{lll}
c_L &\mbox{if}& x<0, \\ c_R &\mbox{if}& x>0 \end{array} \right.  .
\label{riemanndispb}
\ee

\begin{figure}[H]
\begin{center}
 \begin{picture}(400,180)(0,-5)
{\includegraphics[height=6cm,width=12cm]{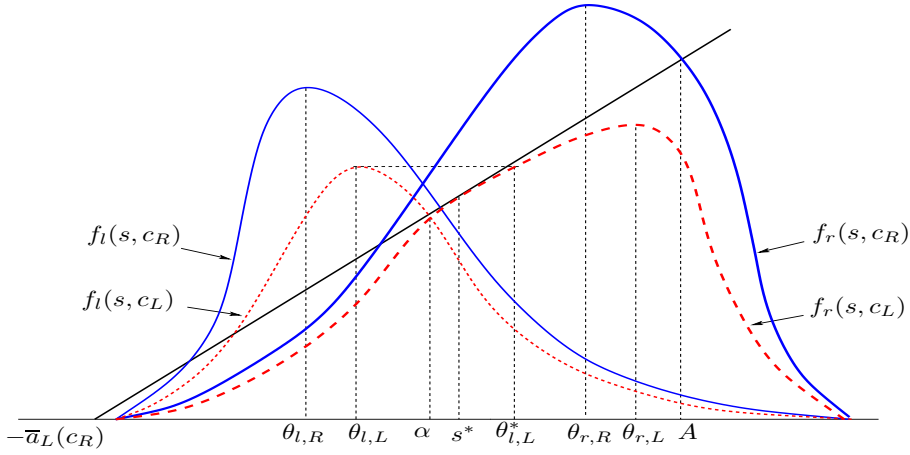}}
\end{picture}
\end{center}
\caption{ Flux functions $f_r(s,c_L), f_r(s,c_R), f_l(s,c_L)$ and $f_l(s,c_R)$ with $c_L > c_R$.}
\label{disrs0}
\end{figure}

 When $c_L > c_R$, the flux functions $f_l$ and $f_r$ satisfy
$ f_l(s,c_L) \leq f_l(s,c_R)$ and $ f_r(s,c_L) \leq f_r(s,c_R)$ for all $s$ in $ (0,1).$ Let $\theta_{lL},\theta_{lR},\theta_{rL}$ and $ \theta_{rR}$ be the points where $f_l(s,c_L), f_l(s,c_R), f_r(s,c_L)$ and $f_r(s,c_R)$ attain their maxima respectively(see Fig.\ref{disrs0}). As there is no discontinuity in $c=c(x,t)$ across the line $x=0$ (see equation (\ref{ccont})) and as
$\sigma$, the speed corresponding to to the $c-$shock, is strictly positive, in Riemann problems we have
$$c(0,t)=c_L\, \forall\,\, t>0.$$
Here we restrict ourselves to the case $c_L > c_R$. The case $c_L < c_R$ can be treated similarly. To study the Riemann problem, we split the problem ( \ref{polydiseq1})  into two problems, one for a scalar conservation law with a discontinuous flux and another for polymer flooding.

\noi{\bf Problem-I}:
\be
\tag{A-2}
\begin{array}{rrll}
 s_t + f_l(s, c_L)_x & =& 0 & \mbox{if}\,\,\, x >0 \\
 s_t + f_r(s, c_L)_x &=& 0  & \mbox{if}\,\,\, x < 0
 \label{problem1}
\end{array}
\ee
The Riemann problem for this equation can be solved as in  \cite{AdiGow03,AdiJafGow04}.\\
\noi{\bf Problem II}:
\be
\tag{A-3}
\begin{array}{rrll}
 s_t + f_r(s, c)_x & =& 0 \\
 (sc+a(c))_t+(cf_r(s,c))_x &=&0
 \label{problem2}
\end{array}
\ee
The Riemann problem for this system can be solved as in section 3.

  We assume without loss of generality that  $f_l(\theta_{lL},c_L) \leq f_r(\theta_{rL},c_L)$. Let $\theta^*_{lL}$ be a point  such that  $f_l(\theta_{lL},c_L) = f_r(\theta^*_{lL},c_L)$ and let $s^* \in (0,1)$ be a point where
$ \frac{\p}{\p s}f_r(s^*,c_L)=\dfrac{f_r(s^*,c_L)}{s^*+\bar{a}_L(c_R)}$, with $\bar{a}_L(c)$ defined as in section 3. Now draw a line through the points
$(-\bar{a}_L(c_R),0)$ and $(s^*,f_r(s^*,c_L))$ which intersects the curve $f_r(s,c_R)$ at a point $A \geq s^*$ (see Fig. \ref{disrs0}).

\begin{itemize}
\item {\bf Case 1:} $ s_L \geq \theta_{lL}$\\
 
 Draw a line through the points
$(-\bar{a}_L(c_R),0)$ and $(\min(s^*,\theta_{lL}^*),f_r(\min(s^*,\theta_{lL}^*),c_L))$ which intersects the curve $f_r(s,c_R)$ at a point $\bar{B} \geq s^*$. For example if $\theta_{lL}^* > s^* $ then $\bar{B}=A$.
 
\item{\bf Case 1a:} $s_R \leq \bar{B}$\\
 
\noi {\bf Step-1}: Let $s_1(x,t)$ be the solution of equation (\ref{problem1}) with initial condition
\be s(x,0) = \left\{ \begin{array}{lll}
s_L &\mbox{if}& x<0, \\ \theta^*_{lL} &\mbox{if}& x>0   
\end{array} \right.  . \nonumber
\label{riemanndispb11}
\ee \nonumber 
\noi {\bf Step-2}: Let $(s_2(x,t),c_2(x,t))$ be the solution of equations (\ref{problem2}) with initial condition
\be s(x,0) = \left\{ \begin{array}{lll}
\theta^*_{lL} &\mbox{if}& x<0, \\ s_R &\mbox{if}& x>0 \end{array} \right.  , \quad
c(x,0) = \left\{ \begin{array}{lll}
c_L &\mbox{if}& x<0, \\ c_R &\mbox{if}& x>0 \end{array} \right.  . \no
\label{riemanndispb12}
\ee \nonumber

   Then the solution to the Riemann problem  (\ref{polydiseq1}), (\ref{riemanndispb}) is given by
 
\be (s(x,t),c(x,t)) = \left\{ \begin{array}{lll}
(s_1(x,t),c_L) &\mbox{if}& x<0, \\ (s_2(x,t),c_2(x,t)) &\mbox{if}& x > 0 \end{array} \right. . \no
\ee \nonumber

\item {\bf Case 1b:} $ s_R > \bar{B}$ \\
 Draw a line through the points
$(-\bar{a}_L(c_R),0)$ and $(s_R,,f_r(s_R,c_R))$ which intersects the curve $f_r(s,c_L)$ at a point $\bar{s}$.

\noi {\bf Step-1:} Let $s_1(x,t)$ be the solution of equation (\ref{problem1}) with initial condition
\be s(x,0) = \left\{ \begin{array}{lll}
s_L &\mbox{if}& x<0, \\ \bar{s} &\mbox{if} &x>0
\end{array} \right.  . \no
\label{riemanndispb21}
\ee \nonumber
{\bf Step-2:} Let $(s_2(x,t),c_2(x,t))$ be the solution of equations (\ref{problem2}) with initial condition
\be s(x,0) = \left\{ \begin{array}{lll}
\bar{s} &\mbox{if}& x<0, \\ s_R &\mbox{if}& x>0 \end{array} \right.  , \quad
c(x,0) = \left\{ \begin{array}{lll}
c_L &\mbox{if}& x<0, \\ c_R &\mbox{if}& x>0 \end{array} \right.  . \no
\label{riemanndispb22}
\ee \nonumber

   Then the solution to the Riemann problem  (\ref{polydiseq1}), (\ref{riemanndispb}) is given by

\be (s(x,t),c(x,t)) = \left\{ \begin{array}{lll}
(s_1(x,t),c_L) &\mbox{if}& x<0, \\ (s_2(x,t),c_2(x,t)) &\mbox{if}& x > 0 \end{array} \right. \no
\ee \nonumber

\item {\bf Case 2:} $ s_L < \theta_{lL}$.

Let  $s_L^*$ be a point  such that  $f_r(s_L^*,c_L)=f_l(s_L,c_L)$ and $\frac{\p}{\p s}f_r(s,c_L)$ at $s=s_L^*  \geq 0$. Draw a line through the points
$(-\bar{a}_L(c_R),0)$ and $(\min(s^*,s_L^*),f_r(\min(s^*,s_L^*),c_L))$ which intersects the curve $f_r(s,c_R)$ at a point $\bar{B}$.

\item {\bf Case 2a:} $s_R \leq \bar{B}$

\noi Step-1: Let $s_1(x,t)$ be the solution of equation (\ref{problem1}) with initial condition
\be s(x,0) = \left\{ \begin{array}{lll}
s_L &\mbox{if}& x<0, \\ s_L^* &\mbox{if}& x>0 
\end{array} \right.  . \no
\label{riemanndispb31}
\ee \nonumber
Step-2: Let $(s_2(x,t),c_2(x,t))$ be the solution of equations (\ref{problem2}) with initial condition
\be s(x,0) = \left\{ \begin{array}{lll}
s_L^* &\mbox{if}& x<0, \\ s_R &\mbox{if}& x>0 \end{array} \right.  , \quad
c(x,0) = \left\{ \begin{array}{lll}
c_L &\mbox{if}& x<0, \\ c_R &\mbox{if}& x>0 \end{array} \right.  .\no
\label{riemanndispb32}
\ee \nonumber

   Then the solution to the Riemann problem  (\ref{polydiseq1}), (\ref{riemanndispb}) is given by

\be (s(x,t),c(x,t)) = \left\{ \begin{array}{lll}
(s_1(x,t),c_L) &\mbox{if}& x<0, \\ (s_2(x,t),c_2(x,t)) &\mbox{if}& x > 0 \end{array} \right. . \no
\ee \nonumber

\item {\bf Case 2b} $ s_R >  \bar{B}$.

 Draw a line through the points
$(-\bar{a}_L(c_R),0)$ and $(s_R,f_r(s_R,c_R))$ which intersects the curve $f_r(s,c_L)$ at a point $\bar{s}$.\\

\noi Step-1: Let $s_1(x,t)$ be the solution of equation (\ref{problem1}) with initial condition
\be s(x,0) = \left\{ \begin{array}{lll}
s_L &\mbox{if}& x<0, \\ \bar{s} &\mbox{if}& x>0 
\end{array} \right.  . \no
\label{riemanndispb41}
\ee \nonumber

\noi  Step-2: Let $(s_2(x,t),c_2(x,t))$ be the solution of equations (\ref{problem2}) with initial condition
\be s(x,0) = \left\{ \begin{array}{lll}
\bar{s} &\mbox{if}& x<0, \\ s_R &\mbox{if}& x>0 \end{array} \right.  , \quad
c(x,0) = \left\{ \begin{array}{lll}
c_L &\mbox{if}& x<0, \\ c_R &\mbox{if}& x>0 \end{array} \right.  . \no
\label{riemanndispb42}
\ee \nonumber

   Then the solution to the Riemann problem  (\ref{polydiseq1}), (\ref{riemanndispb}) is given by

\be (s(x,t),c(x,t)) = \left\{ \begin{array}{lll}
(s_1(x,t),c_L) &\mbox{if}& x<0, \\ (s_2(x,t),c_2(x,t)) &\mbox{if}& x > 0 \end{array} \right. \no
\ee \nonumber
\end{itemize}
{\bf Acknowledgements:}\\
Authors would like to thank anonymous referee for his valuble suggestions
in proving the convergence of approximated solution $\{s_i^n\}$ and
 Sudarshan Kumar for computing the results in Fig.\ref{discomparison}.



\end{document}